\documentclass[11pt,twoside]{article}

\usepackage{amsbsy,amsfonts,amsmath,amssymb,eucal,mathrsfs}
\usepackage[all]{xy}
\usepackage{pstricks,epsfig}
\usepackage{pst-plot}
\usepackage{psfrag}

\newtheorem{definition}{Definition}[section]
\newenvironment{defi}{\begin{definition} \rm}{\end{definition}}

\newtheorem{prop}[definition]{Proposition}

\newtheorem{lemm}[definition]{Lemma}
\newtheorem{fact}[definition]{Fact}

\newtheorem{coro}[definition]{Corollary}
\newtheorem{theo}[definition]{Theorem}
\newtheorem{notation}[definition]{Notation}
\newenvironment{nota}{\begin{notation} \rm}{\end{notation}}
\newtheorem{construction}[definition]{Construction}

\newtheorem{remark}[definition]{Remark}
\newenvironment{rema}{\begin{remark} \rm}{\end{remark}}
\newtheorem{remarks}[definition]{Remarks}

\newtheorem{example}[definition]{Example}
\newenvironment{exam}{\begin{example} \rm}{\end{example}}
\newtheorem{examples}[definition]{Examples}

\newtheorem{nothing}[definition]{$\!\!$}

\newenvironment{proo}{{\flushleft \it Proof.}}{\hfill $\square$
  \vspace{2mm}}
\newenvironment{proo-prop}{{\flushleft \it Proof of Proposition
    \ref{prop-f_4}.}}{\hfill $\square$ \vspace{2mm}}
\newenvironment{conj}{\begin{conjecture} \rm}{\end{conjecture}}
\newtheorem{conjecture}[definition]{Conjecture}

\newtheorem{definition*}{Definition}[section]
\newenvironment{defi*}{\begin{definition*} \rm}{\end{definition*}}
\newtheorem{definitions*}[definition*]{Definitions}
\newenvironment{defis*}{\begin{definitions*} \rm}{\end{definitions*}}
\newtheorem{prop*}[definition*]{Proposition}
\newtheorem{lemm*}[definition*]{Lemma}
\newtheorem{coro*}[definition*]{Corollary}
\newtheorem{theo*}[definition*]{Theorem}
\newtheorem{remark*}[definition*]{Remark}
\newenvironment{rema*}{\begin{remark*} \rm}{\end{remark*}}
\newtheorem{remarks*}[definition*]{Remarks}
\newenvironment{remas*}{\begin{remarks*} \rm}{\end{remarks*}}
\newtheorem{example*}[definition*]{Example}
\newenvironment{exam*}{\begin{example*} \rm}{\end{example*}}
\newtheorem{examples*}[definition*]{Examples}
\newenvironment{exams*}{\begin{examples*} \rm}{\end{examples*}}
\newtheorem{nothing*}[definition*]{$\!\!$}
\newenvironment{noth*}{\begin{nothing*} \rm}{\end{nothing*}}

\newtheorem{commentaire*}[definition*]{Commentaire}

\begin{document}

\def \JK {{\mathfrak{J}}}
\def \IK {{\mathfrak{I}}}
\def \Rh {{\widehat{R}}}
\def \Sh {{\widehat{S}}}
\def \supp {{\rm Supp}}
\def \codim {{\rm Codim}}
\def \b {{\beta}}
\def \T {{\Theta}}
\def \t {{\theta}}
\def \L {{\cal L}}
\def \sca #1#2{\langle #1,#2 \rangle}
\def\pt{\{{\rm pt}\}}
\def\x {{\underline{x}}}
\def\y {{\underline{y}}}
\def\aut{{\rm Aut}}
\def\ra{\rightarrow}
\def\s{\sigma}\def\OO{\mathbb O}\def\PP{\mathbb P}\def\QQ{\mathbb Q}
 \def\CC{\mathbb C} \def\ZZ{\mathbb Z}\def\JO{{\mathcal J}_3(\OO)}
\newcommand{\G}{\mathbb{G}}
\def\proof{\noindent {\it Proof.}\;}
\def\qed{\hfill $\square$}
\def \uh {{\widehat{u}}}
\def \vh {{\widehat{v}}}
\def \fh {{\widehat{f}}}
\def \wh {{\widehat{w}}}
\def \Wh {{{W_{{\rm aff}}}}}
\def \Wt {{\widetilde{W}_{{\rm aff}}}}
\def \Qt {{\widetilde{Q}}}
\def \Waff {{W_{{\rm aff}}}}
\def \Waffm {{W_{{\rm aff}}^-}}
\def \Wpaff {{{(W^P)}_{{\rm aff}}}}
\def \Wtpaff {{{(\widetilde{W}^P)}_{{\rm aff}}}}
\def \Wtaffm {{\widetilde{W}_{{\rm aff}}^-}}
\def \lh {{\widehat{\lambda}}}
\def \pit {{\widetilde{\pi}}}
\def \lt {{{\lambda}}}
\def \xh {{\widehat{x}}}
\def \yh {{\widehat{y}}}
\def \a {\alpha}
\def \b {\beta}
\def \l {\lambda}
\def \t {\theta}
%%%%%ICI ATTENTION \def \T {\theta}

%%%%%%%%% Commandes pe %%%%%%%%%%%%%%

\newcommand{\expxy}{\exp_{x \to y}}
\newcommand{\drat}{d_{\rm rat}}
\newcommand{\dmax}{d_{\rm max}}
\newcommand{\zl}{Z(x,L_x,y,L_y)}

% alg{\~A}{\AA}{!`}bres norm{\~A}{\copyright}es et anneaux usuels

\newcommand{\tiff}{if and only if }
\newcommand{\N}{\mathbb{N}}
\newcommand{\A}{{\mathbb{A}_{\rm Aff}}}
\newcommand{\Ah}{{\mathbb{A}_{\rm Aff}}}
\newcommand{\At}{{\widetilde{\mathbb{A}}_{\rm Aff}}}
\newcommand{\Ht}{{{H}^T_*(\Omega K^{\ad})}}
\renewcommand{\H}{{\rm Hi}}
\newcommand{\Ih}{{I_{\rm Aff}}}
\newcommand{\psit}{{\widetilde{\psi}}}
\newcommand{\xit}{{\widetilde{\xi}}}
\newcommand{\Jt}{{\widetilde{J}}}
\newcommand{\Zt}{{\widetilde{Z}}}
\newcommand{\Xt}{{\widetilde{X}}}
\newcommand{\at}{{\widetilde{A}}}
\newcommand{\Z}{\mathbb Z}
\newcommand{\R}{\mathbb{R}}
\newcommand{\Q}{\mathbb{Q}}
\newcommand{\C}{\mathbb{C}}
\renewcommand{\O}{\mathbb{O}}
\newcommand{\F}{\mathbb{F}}
\newcommand{\p}{\mathbb{P}}
\newcommand{\co}{{\cal O}}
\newcommand{\pos}{{\bf P}}

\renewcommand{\a}{{\alpha}}
\newcommand{\az}{\a_\Z}
\newcommand{\ak}{\a_k}

\newcommand{\rc}{\R_\C}
\newcommand{\cc}{\C_\C}
\newcommand{\hc}{\H_\C}
\newcommand{\oc}{\O_\C}

\newcommand{\rk}{\R_k}
\newcommand{\ck}{\C_k}
\newcommand{\hk}{\H_k}
\newcommand{\ok}{\O_k}

\newcommand{\rz}{\R_Z}
\newcommand{\cz}{\C_Z}
\newcommand{\hz}{\H_Z}
\newcommand{\oz}{\O_Z}

\newcommand{\RR}{\R_R}
\newcommand{\CR}{\C_R}
\newcommand{\HR}{\H_R}
\newcommand{\OR}{\O_R}

\newcommand{\re}{\mathtt{Re}}

\newcommand{\matttr}[9]{
\left (
\begin{array}{ccc}
{} \hspace{-.2cm} #1 & {} \hspace{-.2cm} #2 & {} \hspace{-.2cm} #3 \\
{} \hspace{-.2cm} #4 & {} \hspace{-.2cm} #5 & {} \hspace{-.2cm} #6 \\
{} \hspace{-.2cm} #7 & {} \hspace{-.2cm} #8 & {} \hspace{-.2cm} #9
\end{array}
\hspace{-.15cm}
\right )   }

% alg{\~A}{\AA}{!`}bre

\newcommand{\dual}{{\bf v}}
\newcommand{\com}{\mathtt{Com}}
\newcommand{\rg}{\mathtt{rg}}
\newcommand{\pu}{{\mathbb{P}^1}}
\newcommand{\scal}[1]{\langle #1 \rangle}
\newcommand{\MK}[2]{{\overline{{\rm M}}_{#1}(#2)}}
\newcommand{\mor}[2]{{{\rm Mor}_{#1}(\pu,#2)}}

\newcommand{\fg}{\mathfrak g}
\newcommand{\fgad}{{\mathfrak g}^{\rm ad}}
\renewcommand{\fh}{\mathfrak h}
\newcommand{\fu}{\mathfrak u}
\newcommand{\fz}{\mathfrak z}
\newcommand{\fn}{\mathfrak n}
\newcommand{\fe}{\mathfrak e}
\newcommand{\fp}{\mathfrak p}
\newcommand{\ft}{\mathfrak t}
\newcommand{\fl}{\mathfrak l}
\newcommand{\fq}{\mathfrak q}
\newcommand{\fsl}{\mathfrak {sl}}
\newcommand{\fgl}{\mathfrak {gl}}
\newcommand{\fso}{\mathfrak {so}}
\newcommand{\fsp}{\mathfrak {sp}}
\newcommand{\ff}{\mathfrak {f}}

\newcommand{\ad}{{\rm ad}}
\newcommand{\jad}{{j^\ad}}
\newcommand{\id}{{\rm id}}

%%%% Poids et racines   %%%%

\newcommand{\dynkinadeux}[2]
{
$#1$
\setlength{\unitlength}{1.2pt}
\hspace{-3mm}
\begin{picture}(12,3)
\put(0,3){\line(1,0){10}}
\end{picture}
\hspace{-2.4mm}
$#2$
}

\newcommand{\mdynkinadeux}[2]
{
\mbox{\dynkinadeux{#1}{#2}}
}

\newcommand{\dynkingdeux}[2]
{
$#1$
\setlength{\unitlength}{1.2pt}
\hspace{-3mm}
\begin{picture}(12,3)
\put(1,.8){$<$}
\multiput(0,1.5)(0,1.5){3}{\line(1,0){10}}
\end{picture}
\hspace{-2.4mm}
$#2$
}

\newcommand{\poidsesix}[6]
{
\hspace{-.12cm}
\left (
\begin{array}{ccccc}
{} \hspace{-.2cm} #1 & {} \hspace{-.3cm} #2 & {} \hspace{-.3cm} #3 &
{} \hspace{-.3cm} #4 & {} \hspace{-.3cm} #5 \vspace{-.13cm}\\
\hspace{-.2cm} & \hspace{-.3cm} & {} \hspace{-.3cm} #6 &
{} \hspace{-.3cm} & {} \hspace{-.3cm}
\end{array}
\hspace{-.2cm}
\right )      }

\newcommand{\copoidsesix}[6]{
\hspace{-.12cm}
\left |
\begin{array}{ccccc}
{} \hspace{-.2cm} #1 & {} \hspace{-.3cm} #2 & {} \hspace{-.3cm} #3 &
{} \hspace{-.3cm} #4 & {} \hspace{-.3cm} #5 \vspace{-.13cm}\\
\hspace{-.2cm} & \hspace{-.3cm} & {} \hspace{-.3cm} #6 &
{} \hspace{-.3cm} & {} \hspace{-.3cm}
\end{array}
\hspace{-.2cm}
\right |      }

\newcommand{\poidsesept}[7]{
\hspace{-.12cm}
\left (
\begin{array}{cccccc}
{} \hspace{-.2cm} #1 & {} \hspace{-.3cm} #2 & {} \hspace{-.3cm} #3 &
{} \hspace{-.3cm} #4 & {} \hspace{-.3cm} #5 & {} \hspace{-.3cm} #6
\vspace{-.13cm}\\
\hspace{-.2cm} & \hspace{-.3cm} & {} \hspace{-.3cm} #7 &
{} \hspace{-.3cm} & {} \hspace{-.3cm}
\end{array}
\hspace{-.2cm}
\right )      }

\newcommand{\copoidsesept}[7]{
\hspace{-.12cm}
\left |
\begin{array}{cccccc}
{} \hspace{-.2cm} #1 & {} \hspace{-.3cm} #2 & {} \hspace{-.3cm} #3 &
{} \hspace{-.3cm} #4 & {} \hspace{-.3cm} #5 & {} \hspace{-.3cm} #6
\vspace{-.13cm}\\
\hspace{-.2cm} & \hspace{-.3cm} & {} \hspace{-.3cm} #7 &
{} \hspace{-.3cm} & {} \hspace{-.3cm}
\end{array}
\hspace{-.2cm}
\right |      }

\newcommand{\poidsehuit}[8]{
\hspace{-.12cm}
\left (
\begin{array}{cccccc}
{} \hspace{-.2cm} #1 & {} \hspace{-.3cm} #2 & {} \hspace{-.3cm} #3 &
{} \hspace{-.3cm} #4 & {} \hspace{-.3cm} #5 & {} \hspace{-.3cm} #6 &
{} \hspace{-.3cm} #7   \vspace{-.13cm}\\
\hspace{-.2cm} & \hspace{-.3cm} & {} \hspace{-.3cm} #8 &
{} \hspace{-.3cm} & {} \hspace{-.3cm}
\end{array}
\hspace{-.2cm}
\right )      }

\newcommand{\copoidsehuit}[8]{
\hspace{-.12cm}
\left |
\begin{array}{cccccc}
{} \hspace{-.2cm} #1 & {} \hspace{-.3cm} #2 & {} \hspace{-.3cm} #3 &
{} \hspace{-.3cm} #4 & {} \hspace{-.3cm} #5 & {} \hspace{-.3cm} #6 &
{} \hspace{-.3cm} #7  \vspace{-.13cm}\\
\hspace{-.2cm} & \hspace{-.3cm} & {} \hspace{-.3cm} #8 &
{} \hspace{-.3cm} & {} \hspace{-.3cm}
\end{array}
\hspace{-.2cm}
\right |      }

\newcommand{\im}{\mathtt{Im}}

%%%% Caligraphic letters %%%%%%%%%%%%%%

\def\cA{{\cal A}} \def\cC{{\cal C}} \def\cD{{\cal D}} \def\cE{{\cal E}}
\def\cF{{\cal F}} \def\cG{{\cal G}} \def\cH{{\cal H}} \def\cI{{\cal I}}
\def\cK{{\cal K}} \def\cL{{\cal L}} \def\cM{{\cal M}} \def\cN{{\cal N}}
\def\cO{{\cal O}}
\def\cP{{\cal P}} \def\cQ{{\cal Q}} \def\cT{{\cal T}} \def\cU{{\cal U}}
\def\cV{{\cal V}} \def\cX{{\cal X}} \def\cY{{\cal Y}} \def\cZ{{\cal Z}}

\def \g {{\gamma}}

\def \pp {\odot}
\def \tr {{}^t}
\def \ct {{}_Tc}
\def \lcom {$\Lambda$-(co)minuscule }
\def \lmin {$\Lambda$-minuscule }
\renewcommand{\ok}[2]{{#1} \cdot {#2} = {#1} \pp {#2}}
\newcommand{\notok}[2]{{#1} \cdot {#2} \not = {#1} \pp {#2}}
\newcommand{\oknu}[3]{c_{{#1},{#2}}^{{#3}} = t_{{#1},{#2}}^{{#3}} \cdot
m_{{#1},{#2}}^{{#3}}}
\newcommand{\notoknu}[3]{c_{{#1},{#2}}^{{#3}} \not = t_{{#1},{#2}}^{{#3}} \cdot
m_{{#1},{#2}}^{{#3}}}

 \title{Towards a Littlewood-Richardson rule \\
for Kac-Moody homogeneous spaces}
 \author{P.-E. Chaput, N. Perrin}

\maketitle

\begin{abstract}

We prove a general combinatorial formula yielding the
intersection number $c_{u,v}^w$ of three particular
$\Lambda$-minuscule Schubert
classes in any Kac-Moody homogeneous space,
generalising the Littlewood-Richardson rule.
The combinatorics are
based on jeu de taquin rectification in a poset defined by the
heap of $w$.

\end{abstract}

 {\def\thefootnote{\relax}
 \footnote{ \hspace{-6.8mm}
 Key words: Littlewood-Richardson rule, Schubert calculus, Kac-Moody
homogeneous spaces, jeu de taquin. \\
 Mathematics Subject Classification: 14M15, 14N35}
 }

\section{Introduction}

Schubert calculus is an old important problem. Its main focus is the
computation of the structure constants (the Littlewood-Richardson
coefficients) in the cup product of Schubert classes in the cohomology 
of a homogeneous space. Schubert calculus is now well understood 
in many aspects (see for example \cite{borel}, \cite{demazure}, 
\cite{BGG}, \cite{duan})
but several problems remain open. In particular a
combinatorial formula for the Littlewood-Richardson coefficients 
is not known in general. The most striking example of such a formula
is the celebrated Littlewood-Richardson rule computing these
coefficients for Grassmannian using jeu de taquin (see Section
\ref{section-taquin}). This rule was conjectured by D.E. Littlewood and 
A.R. Richardson in \cite{LR} and proved by M.P. Sch{\"u}tzenberger in 
\cite{schutz}. For a historical account, the reader may consult 
\cite{vanlee}. Generalisation to minuscule and cominuscule homogeneous 
spaces of classical types were proved by D. Worley \cite{worley} and 
P. Pragacz \cite{pragacz}. Recently, this rule has
been extended to exceptional minuscule homogeneous spaces by H. Thomas and A. 
Yong \cite{TY}.

In this paper, we largely extend their rule to any homogeneous space
$X$ for certain cohomology classes called $\Lambda$-minuscule classes (see
Definition \ref{defi-minuscule}). For $X$ minuscule, any cohomology
class is $\Lambda$-minuscule. We even prove this rule in many cases
where the space $X$ is homogeneous under a Kac-Moody group.

Let us be more precise and introduce some notation. Let $G$ be a
Kac-Moody group and let $P$ be a parabolic subgroup of $G$. Let $X$ be
the homogeneous space $G/P$. A basis of the cohomology group
$H^*(X,\Z)$ is indexed by the set of minimal length representative
$W^P$ of the quotient $W/W_P$ where $W$ is the Weyl group of $G$ and
$W_P$ the Weyl group of $P$. Let us denote with $\s^w$ the Schubert
class corresponding to $w\in W^P$.
The Littlewood-Richardson coefficients are the contants
$c_{u,v}^w$ defined for $u$ and $v$ in $W^P$ by the formula:
$$\s^u\cup\s^v=\sum_{w\in W^P}c_{u,v}^w\s^w.$$
Let us denote with $\Lambda$ the dominant weight associated to
$P$. Following Dale Peterson, we define special elements in $W^P$
called $\Lambda$-minuscule (see Definition
\ref{defi-minuscule}). These elements have the nice property of being
fully commutative: they admit a unique reduced expression up to
commuting relations. In particular, they have a well defined heap which 
is a colored poset, the colors being simple roots (see
Definition \ref{defi-heap}, this was first introduced by X. G. Viennot in 
\cite{viennot}, we use J. Stembridge's definition in
\cite{stembridge}, heaps were reintroduced in \cite{small} as Schubert
quivers). One of the major points we shall use here to define our
combinatorial rule is the fact proved by R. Proctor
\cite{proctor_taquin} that these heaps do have the jeu de taquin
property (see Section \ref{section-taquin}). In particular, given two
elements $u$ and $v$ in $W$ smaller than a $\Lambda$-minuscule element
$w$, we define combinatorially using jeu de taquin an integer
$t_{u,v}^w$ (see Proposition \ref{prop_independant_t}). We make the
following conjecture:

\begin{conj}
\label{conj1}
For $w$ a $\Lambda$-minuscule element and $u$ and $v$ in $W$ smaller
than $w$, we have the equality $c_{u,v}^w=t_{u,v}^w$.
\end{conj}

Following \cite{TY}, we extend these considerations to
$\Lambda$-cominuscule elements (see Definition \ref{defi-minuscule})
defined using $\Lambda$-minuscule elements in the Langlands dual
group. 
Let $S(\Lambda)$ denote the set of roots $\alpha$ such that
$\scal{\Lambda,\alpha^\vee} > 0$. Let $w$ be a $\Lambda$-cominuscule
element, if $w=s_{\alpha_1} \cdots s_{\alpha_l}$ is a reduced
expression we define
$$
\displaystyle{
m(w) := \prod_{\stackrel{\stackrel{i \in [1,l],\alpha \in S(\Lambda) ,}
{(\alpha , \alpha ) > ( \alpha_i,\alpha_i ) ,\ i \geq (\alpha,1) }}{}} 
\frac{ ( \alpha , \alpha ) } { ( \alpha_i,\alpha_i ) }\ \ ,
}
$$
were $( \cdot , \cdot )$ is any $W$-invariant scalar product and $(\a,1)$ 
is the minimal element of the heap colored by $\a$. This 
only depends on $w$ and not on the choice of a reduced expression.
Let $u$ and $v$ in $W$ smaller than $w$, we denote with $m_{u,v}^w$
the number $m(w)/(m(u) \cdot m(v))$. If $w$ is $\Lambda$-minuscule,
the same definition gives $m_{u,v}^w=1$, by Lemma
\ref{lemm-racines-courtes}. We extend the previous
conjecture as follows:

\begin{conj}
\label{conj2}
For $w$ a $\Lambda$-cominuscule element and $u$ and $v$ in $W$ smaller
than $w$, we have the equality $c_{u,v}^w=m_{u,v}^w t_{u,v}^w$.
\end{conj}

Our inspiration in the work of H. Thomas and A. Yong is very clear
with these conjectures. The first evidences for them are the
Littlewood-Richarson rule (\emph{i.e.} Conjecture \ref{conj1} is true
for $X$ a Grassmannian) and the result of H. Thomas and A. Yong
\cite{TY} proving that conjectures \ref{conj1} and \ref{conj2} are
true for $X$ a minuscule or a cominuscule homogeneous space. Our main
result is a proof of these conjectures in many cases including all
finite dimensional homogeneous spaces $X$. Indeed, we define for $w$ a
$\Lambda$-minuscule or $\Lambda$-cominuscule element of the Weyl group
the condition of being slant-finite-dimensional (see Definition
\ref{defi-slant}). This includes all $\Lambda$-minuscule or
$\Lambda$-cominuscule elements in the Weyl group $W$ of a finite
dimensional group $G$. Our main result is the following:

\begin{theo}
\label{main_theo}
Let $G/P$ be a Kac-Moody homogeneous space where $P$ corresponds
to the dominant weight $\Lambda$. Let $u,v,w \in W$ be
$\Lambda$-(co)minuscule. Assume that $w$ is slant-finite-dimensional.
Then we have $c_{u,v}^w = m_{u,v}^wt_{u,v}^w$.
\end{theo}

Let us observe here that we restrict the statement to slant-finite 
dimensional elements essentially for technical reasons: this simplifies 
a lot the combinatorics involved and allows us to find easily generators 
of the cohomology algebra.

The strategy of proof is very similar to the one of H. Thomas and
A. Yong but we add two powerful ingredients: first we prove \emph{a priori}
that jeu de taquin numbers $t_{u,v}^w$ as well as modified jeu de
taquin numbers $m_{u,v}^wt_{u,v}^w$ define a commutative and
associative algebra (see Subsection
\ref{subsection-algebra-taquin}). As an example of the strength of this
fact, we will reprove that in classical
(co)minuscule homogeneous spaces the modified jeu de taquin coefficients are
equal to the intersection numbers, assuming that only very
few intersection numbers are known. For example, to reprove the case
of Grassmannians we only need to assume that we know the cohomology
ring of the 4-dimensional Grassmannian $G(2,4)$: see
Lemma \ref{lemm-an}. We believe that
this was not possible without this fact only with the arguments of
H. Thomas and A. Yong. Our main use of this result relies on the
fact that we will only need to prove Conjectures \ref{conj1}
and \ref{conj2} for a
system of generators of the cohomology.

Another powerful tool is the decomposition of any $\Lambda$-minuscule
element into a product of so-called slant-irreducible elements and the
classification, by Proctor and Stembridge, of the irreducible ones.
We are thus able to reduce the proof of Theorem
\ref{main_theo} to the classical cases plus a finite number of
exceptional ones: see Subsection \ref{subsection-decomposable}.

To prove theorem \ref{main_theo} we need
two more ingredients already contained in \cite{TY}: the fact that our
rule is compatible with the Chevalley formula and a Kac-Moody
recursion which enables to boil the computation of certain
Littlewood-Richardson coefficients down to the computation of other
Littlewood-Richardson coefficients in a smaller group. This idea of
recursion was contained in the work of H. Thomas and A. Yong
\cite{TY}, however we had to adapt their proof 
in the general Kac-Moody situation.
This is done in Subsection \ref{subsection-recursion}.

Before describing in more details the sections in this article, let us
remark that, even if $\Lambda$-(co)minuscule elements may be rare in
certain homogeneous spaces, our result can be applied to compute an
explicit presentation of the cohomology ring of adjoint varieties
and thus to compute
all their Littlewood-Richardson coefficients. This will be done
in a subsequent work \cite{CP}.

In Section \ref{section-taquin}, we define $\Lambda$-minuscule and
$\Lambda$-cominuscule elements and the combinatorial invariants
$t_{u,v}^w$ and $m_{u,v}^w$. We state our main conjecture. We prove
that this conjecture is compatible with the Chevalley formula and
define an associative and commutative algebra using these
combinatorial invariants. 
We also define the notion of Bruhat
recursion and prove
that the Littlewood-Richardson coefficients
$c_{u,v}^w$ satisfy Bruhat recursion. In Section
\ref{section-general-argument}, we define the notion of
slant-finite-dimensional elements and state our main result. We explain
our strategy to prove Theorem \ref{main_theo}. We prove several lemmas
implying that the two products (the cup product and the combinatorial
product) are equal. In Section \ref{section-calc}, we prove by a case
by case analysis that Theorem \ref{main_theo} holds for simply laced
Kac-Moody groups. In type $A$, Lemma \ref{lemm-an} gives a very short
proof (using the fact that our combinatorial product is commutative
and associative) of the classical Littlewood-Richardson rule. In
Section \ref{section-pliage}, we explain how, using foldings, we can
deduce Theorem \ref{main_theo} in the non simply laced cases, using 
the simply laced case. We will need in particular to make involved
computations to deal with a single coefficient in one case related
to $F_4$.

\vskip 0.3 cm

{\bf Acknowledgement:} We would like to thank G{\'e}rald Gaudens
and Antoine Touz{\'e} for discussions about topology of infinite 
dimensional spaces. Both authors are thankful to the Max-Planck
Institut for providing ideal research conditions, and
Pierre-Emmanuel Chaput thanks the Nantes university for giving
a special funding.

\vskip 0.3 cm

{\bf Convention:} We work over an algebraically closed field of
characteristic zero. We will use several times the notation in
\cite{bou} especially for labelling the simple roots of a semisimple
Lie algebra.

\tableofcontents

\newcommand{\wv}{{\bf v}}
\newcommand{\ww}{{\bf w}}

\section{Jeu de taquin}

\label{section-taquin}

\subsection{The jeu de taquin property}

\label{subsection-taquin}

Jeu de taquin is a combinatorial game encoding all Schubert
intersection numbers for (co)mi\-nuscu\-le varieties, as it was shown
by H. Thomas and A. Yong in
\cite{TY}. For the convenience of the reader we recall their
definition of the jeu de taquin. Let $P$ be a poset which we assume to
be {\it bounded below}, meaning that for any $x \in P$ the set $\{ y :
y \leq x \}$ is finite. Elements of $P$ will be called boxes.
Recall that a subset $\lt$ of a poset $P$ is an order ideal if for $x\in \lt$
and $y\in P$ we have the implication $(y\leq x\Rightarrow y\in\lt)$.
We denote with $I(P)$ the set of finite order ideals of $P$. For
$\lambda \subset \nu$ two finite order ideals in $P$ we
denote with $\nu / \lambda$ the pair $(\lt,\nu)$. Any such pair
is called a skew shape. A standard tableau $T$ of skew shape $\nu /
\lambda$ is an increasing bijective map $(\nu-\lambda )\to 
[1,d]$, where $d$ is the cardinal of the set theoretic difference 
$(\nu-\l)$. 

Consider $x \in \lt$ and maximal in $\lt$ among the 
elements that are below some element of $(\nu - \l)$. We associate another 
standard tableau $j_x(T)$ (of a different skew shape) arising
from $T$ : let $y$ be the box of $(\nu- \l)$ with the smallest
label, among those that cover $x$. Move the label of the box $y$ to $x$, 
leaving $y$ vacant. Look for the smallest label of
$(\nu-  \l)$ that covers $y$ and repeat the process. The tableau $j_x(T)$ is 
outputted when no more such moves are possible. A
rectification of $T$ is the result of an iteration of jeu de taquin slides 
until we terminate at
a standard tableau which shape is an order ideal. By the assumption that $P$ 
is bounded below this
will occur after a finite number of slides.

According to Proctor \cite{proctor_taquin}, we will say that $P$ has the {\it 
jeu de taquin property} if the rectification of any tableau does not depend 
on the choices of the empty boxes used to perform
jeu de taquin slides.

\subsection{Jeu de taquin poset associated with a 
$\Lambda$-(co)minuscule element}
\label{subsection-taquin-poset}

Let us first recall some results of Proctor and Stembridge.
Let $A$ be a symmetrisable matrix, $G$ be the associated
symmetrisable Kac-Moody group and let $(\varpi_i)_{i \in I}$ be the
set of fundamental weights.
Let $W$ be the Weyl group of $A$ with generators denoted with 
$s_i$. Note that $W$ acts on the root system $R(A)$ of $A$, and since the Weyl
group of the dual root system $R(\tr\!\! A)$ is isomorphic with $W$ in a canonical 
way,
$W$ also acts on $R(\tr\!\! A)$. The fundamental weights of $R(\tr\!\! A)$ will be
denoted with $\varpi_i^\vee$.
According to Dale Peterson
\cite[p.273]{proctor_minuscule} we give the following definition:

\begin{defi}
\label{defi-minuscule}
Let $\Lambda = \sum_{i} \Lambda_i \varpi_i$ be a dominant weight.
\begin{itemize}
\item
An element $w \in W$ is $\Lambda$-minuscule if there exists a
reduced decomposition
$w = s_{i_1} \cdots s_{i_l}$ such that for any $k \in [1,l]$ we have
$s_{i_k} s_{i_{k+1}}\cdots s_{i_l}(\Lambda) 
= s_{i_{k+1}}\cdots s_{i_l}(\Lambda) - \alpha_{i_k}$.
\item
$w$ is $\Lambda$-cominuscule if $w$ is $(\sum \Lambda_i 
\varpi_i^\vee)$-minuscule.
\item
We will write that $w$ is $\Lambda$-(co)minuscule when we mean that $w$ is
either $\Lambda$-minuscule or $\Lambda$-cominuscule. We denote with $W_m$ the
set of all $\Lambda$-(co)minuscule elements of $W$.
\item $w$ is fully commutative if all the reduced expressions of $w$
  can be deduced one from the other using commutation relations.
\end{itemize}
\end{defi}

By \cite[Proposition 2.1]{stembridge}, any $\Lambda$-minuscule element is fully
commutative. Since the property of being fully commutative depends on $W$ only,
and not on the underlying root system, $\Lambda$-cominuscule elements are
also fully commutative. Moreover Stembridge shows that if the above
condition, defining $\Lambda$-minuscule elements, holds for one reduced
expression $w = s_{i_1} \cdots s_{i_l}$, then it holds for any reduced 
expression of $w$.

For the convenience of the reader we recall the definition of the heap of $w$
given by Stembridge \cite{stembridge} (except that we reverse the order):

\begin{defi}
\label{defi-heap}
Let $w \in W$ be fully commutative and
let $w = s_{i_1} \cdots s_{i_l}$ be a reduced expression. The heap 
$H(w)$ of $w$
is the set $[1,l]$ ordered by the transitive closure of the
relations ``$p$ is smaller than $q$'' if $p > q$ and $s_{i_p}$ and $s_{i_q}$ 
do not commute.
\end{defi}
As Stembridge explains, the full
commutativity implies that the heap is well-defined up to isomorphisms of posets.
Moreover he shows the following (he shows this for $\Lambda$-minuscule 
elements, the statement for $\Lambda$-cominuscule elements follows because 
the statement only depends on the Weyl group):

\begin{prop}
\label{prop-intervalle}
Let $w$ be $\Lambda$-(co)minuscule. There is an order-preserving bijection
between the set of order ideals of $H(w)$ and the Bruhat interval $[e,w]$.
\end{prop}
The bijection maps an ideal 
$\lambda = \{ n_1 , \ldots , n_k \}$ to the element
$u=s_{n_1} \cdots s_{n_k}$.

\begin{prop}
\label{prop_q_taquin}
Let $w \in W$ be $\Lambda$-(co)minuscule. The poset $H(w)$
has the jeu de taquin property.
\end{prop}
\begin{proo}
If $w$ is $\Lambda$-minuscule, by \cite[Corollary 4.3]{stembridge},
$H(w)$ is a $d$-complete 
poset (the precise definition of $d$-completeness is given in
\cite[Section 3]{proctor_minuscule}). By \cite[Theorem 5.1]{proctor_taquin},
any $d$-complete poset has the jeu de taquin property, proving the proposition.
Since the definition of the heap $H(w)$ does not involve the root system, the 
same property holds for $w$ a $\Lambda$-cominuscule element.
\end{proo}

\begin{prop}
\label{prop_independant_t}
Let $w \in W$ be $\Lambda$-(co)minuscule and let $\lambda,\mu,\nu$ be
order ideals in $H(w)$. 
Then the number of tableaux of shape
$\nu / \lt$ which rectify on a standard tableau $U$ of shape $\mu$ does
not depend on the 
given standard tableau $U$ of shape $\mu$. Denote with 
$t_{\lambda,\mu}^\nu(W)$ this number: we have
$t_{\lambda,\mu}^\nu(W) = t_{\mu,\lambda}^\nu(W)$. 
\end{prop}
When $W$ will be clear from the context, the notation $t_{\lambda,\mu}^\nu(W)$
will be simplified to $t_{\lambda,\mu}^\nu$.
\begin{proo}
In \cite[Section 4]{TY}, the authors study properties of the jeu de taquin 
on so-called (co)minuscule posets, which are a very special class of posets 
with the jeu de taquin property. In fact they use two main properties of 
these posets, namely the jeu de taquin property and the fact that there
is a decreasing involution on these posets. However, this involution is 
used only for results
involving the Poincar{\'e} duality. As one readily checks, Proposition 4.2(b-c), 
Theorem 4.4,
its Corollary 4.5 and the
first equality of Corollary 4.6 are still true for any poset enjoying the 
jeu de taquin property.
The last two statements are the two claims of the proposition.
\end{proo}

\begin{rema}
As the proof shows, a similar result holds for any poset having the 
jeu de taquin property.
\end{rema}

We now prove an easy combinatorial lemma for
$\Lambda$-(co)minuscule elements.

\begin{lemm}
\label{lemm-racines-courtes}
Let $\Lambda$ be a fundamental weight with corresponding simple root 
$\alpha_\Lambda$.
Let $w = s_{\alpha_1} \cdots s_{\alpha_l}$ a
reduced expression of an element in $W$. Let $i \in [1,l]$. If
$w$ is $\Lambda$-minuscule then the root $\alpha_i$ cannot be
shorter than $\alpha_\Lambda$, and if $w$ is $\Lambda$-cominuscule 
then $\alpha_i$
cannot be longer than $\alpha_\Lambda$.
\end{lemm}
\begin{proo}
It is enough to consider the case when $w$ is $\Lambda$-minuscule.
Write $w = s_{\alpha_1} \cdots s_{\alpha_l}$
and assume on the contrary that there exists an integer $i$ such that
$(\alpha_i,\alpha_i) < (\alpha_\Lambda,\alpha_\Lambda)$. Let then $i_0$ be the
maximal such integer.
Since $\scal{\Lambda,\alpha_{i_0}} = 0$
(in fact $\Lambda$ is fundamental and $\alpha_{i_0} \not = \alpha_\Lambda$), 
we have 
$1 = \scal{s_{i_0+1} \cdots s_l(\Lambda) , \alpha_{i_0}^\vee} = - \sum_{i>i_0}
\scal{\alpha_i,\alpha_{i_0}^\vee}$, so there exists
$i>i_0$ such that $\scal{\alpha_i,\alpha_{i_0}^\vee} < 0$. Since 
$\alpha_{i_0}$ is shorter
than $\alpha_i$ we have $\scal{\alpha_i,\alpha_{i_0}^\vee} < -1$. 
Furthermore, for 
any $j>i_0$, we have the inequalities $(\a_{i_0},\a_{i_0})< 
(\a_{\Lambda},\a_{\Lambda}) \leq (\a_{j},\a_{j})$ thus 
$\a_j\neq\a_{i_0}$ and $\scal{\a_{j},a_{i_0}^\vee}\leq0$. This 
contradicts the above equality $\sum_{i>i_0} \scal{\alpha_i,\alpha_{i_0}^\vee}=-1$.
\end{proo}

\begin{rema}
Let $w$ be a $\Lambda$-(co)minuscule element and let $D$ be the subdiagram of 
the Dynkin diagram made of simple roots appearing in a reduced expression 
of $w$. Let $A$ be the generalised Cartan matrix associated to $D$, then with
arguments similar to those in the previous lemma one can show that:
for any couple $i<j$, if $a_{i,j}\neq0$, then one of the equalities
$a_{i,j}=-1$ or $a_{j,i}=-1$ holds.
\end{rema}

\vskip .5cm

We now recall some notation of \cite{proctor-classification}
and \cite{stembridge}, and introduce some new ones. If $D$ is a marked
diagram and $d \in D$, then we say that $(D,d)$ is a marked diagram.
A {\bf $D$-colored poset} is the data of a poset $P$ and a map $c:P \to D$
satisfying the condition: if $s_{c(i)}s_{c(j)}\neq s_{c(j)}s_{c(i)}$, then
$i\leq j$ or $j\leq i$ in $P$. 
To such a poset is associated an element $w$ of the
Weyl group of $D$ defined by $w = \prod_{p \in P} s_{c(p)}$, where the order
in this product is any order compatible with the partial order in $P$.
We say that $P$ is {\bf $d$-(co)minuscule} if $w$ is $\Lambda$-(co)minuscule 
for $\Lambda$ the fundamental weight corresponding to $d$. In the sequel, 
we shall assume that the element $w$ corresponding to the poset $P$ is
$\Lambda$-(co)minuscule. 

If $P$ is a $D$-colored poset with coloring function $c:P \to D$,
$\alpha \in D$ and $i$ is an integer, we denote with 
$(\alpha,i) \in P$ the unique element $p$, if it exists, such that
$c(p)=\alpha$ and such that
$\# \{ q \leq p : c(q)=\alpha \} = i$. In particular,
for each $\alpha$ in $c(P)$, $(\a,1) \in P$ is the minimal element colored
by $\a$. The set of all elements of the form $(\a,1)$ is an ideal in
$P$ called the {\bf rooted tree} of $P$ and denoted with $T$. The map
$\a \mapsto (\a,1)$ establishes a bijection from $c(P)$ to $T$ which
is a poset, thus
yielding a partial order on $c(P)$. We say
that $P$ is {\bf slant-irreducible} if each color in $c(P)$ which is non
maximal with respect to this order is the color of at least two elements
in $P$. In \cite{proctor-classification} and \cite{stembridge-classification},
the $D$-colored slant-irreducible $d$-minuscule posets are classified for any
marked Dynkin diagram $(D,d)$.

If $(p_i)_{i\in[1,k]}$ are elements of a poset $P$, we denote with
$\scal{(p_i)_{i\in[1,k]}}$ the ideal generated by $(p_i)_{i\in[1,k]}$.

\subsection{Conjecture on a general Littlewood-Richardson rule}

We now are in position to state a conjecture relating the Schubert 
calculus and the jeu de taquin.
Let  $\Lambda$ be a dominant weight. Let $X=G/P$ be the homogeneous space
corresponding to $\Lambda$, $W_P$ be the Weyl group of $P$, and
$W^P$ the set of minimum length representatives of the coset $W/W_P$.
Let $(\s^w)_{w \in W^P}$ denote the basis of the cohomology of $G/P$
dual to the Schubert basis in homology (see \cite[Proposition
11.3.2]{kumar}). We denote with $c_{u,v}^w$ the integer coefficients such that
$\s^u \cup \s^v = \sum c_{u,v}^w \s^w$.
Note the following:

\begin{fact}
\label{fait-min-length}
If $w \in W$ is $\Lambda$-(co)minuscule then $w \in W^P$.
\end{fact}
\begin{proo}
We may assume that $w$ is $\Lambda$-minuscule.
Write a length additive expression $w = v p$ with $v \in W^P$
and $p \in W_P$. By
\cite[Proposition 2.1]{stembridge} any reduced expression of a 
$\Lambda$-minuscule
element satisfies the condition of Definition~\ref{defi-minuscule}, thus
$p(\Lambda) = \Lambda$ implies $p=e$; thus $w \in W^P$.
\end{proo}

On the other hand, let $w \in W$ be $\Lambda$-(co)minuscule and $u,v 
\in W$ be less or equal to $w$. To
$u$ and $v$ we can associate order ideals $\lambda(u),\l(v)$
of the poset $H(w)$ of $w$ by Proposition~\ref{prop-intervalle}.
Recall the definition of $t_{\l(u),\l(v)}^{H(w)}$ in Proposition
\ref{prop_independant_t}; this
number will be denoted just with $t_{u,v}^w$.

Let $S(\Lambda)$ denote the set of roots $\alpha$ such that
$\scal{\Lambda,\alpha^\vee} > 0$. 
If $u=s_{\alpha_1} \cdots s_{\alpha_l}$ is
a reduced expression we define
$$
\displaystyle{
m(u) := \prod_{\stackrel{\stackrel{i \in [1,l],\alpha \in S(\Lambda) ,}
{(\alpha , \alpha ) > ( \alpha_i,\alpha_i ) ,\ i \geq (\alpha,1) }}{}} 
\frac{ (\alpha , \alpha) } { ( \alpha_i,\alpha_i ) }\ \ ,
}
$$
were $( \cdot , \cdot )$ is any $W$-invariant scalar product.
Let $u,v \leq w \in W$, we denote with $m_{u,v}^w$ the number
$m(w)/(m(u) \cdot m(v))$.

\begin{conj}
\label{main_conj}
Let $w\in W$ be $\Lambda$-(co)minuscule and $u,v \in W$
with $u,v \leq w$.
Then the Schubert intersection number $c_{u,v}^w$ is equal to the jeu de taquin combinatorial
number $m_{u,v}^w \cdot t_{u,v}^w$.
\end{conj}

By \cite{TY} this conjecture holds for $G/P$ a (co)minuscule homogeneous space
and Theorem~\ref{main-theo} proves it when $G/P$ is a finite dimensional
homogeneous space.
Our strategy of proof is essentially the same as in \cite{TY}:
we argue that the numbers
$c_{u,v}^w$ and $m_{u,v}^w \cdot t_{u,v}^w$ both satisfy
some recursive identities
(this holds for any $G/P$), and then we check in the particular case
of finite dimensional varieties that these
identities together with a few number of equalities
$c_{u,v}^w = m_{u,v}^w \cdot t_{u,v}^w$ imply the theorem.
The recursive identities are:
\begin{itemize}
\item
The numbers $m_{u,v}^w \cdot t_{u,v}^w$ satisfy the same
identity as the identity on the numbers
$c_{u,v}^w$ implied by the Chevalley formula:
see Subsection~\ref{subsection-chevalley}.
\item A Kac-Moody recursion which is a general procedure drawing down
  the computation of some numbers $c_{u,v}^w$ (resp. $t_{u,v}^w$) for
  $G/P$ to the computation of the similar numbers for a quotient $H/Q$
  with $H$ a Levi subgroup of $G$: see 
Subsection~\ref{subsection-recursion}.
\item
Jeu de taquin defines a natural algebra with basis indexed by all 
$\Lambda$-(co)minuscule elements
which is commutative and associative (and will turn out to be, once the 
theorem is proved, isomorphic with a quotient of $H^*(G/P)$): see Subsection
\ref{subsection-algebra-taquin}.
\end{itemize}

The last point was not used in \cite{TY}. We will see that
it simplifies a lot our argument, since
it implies that to prove the theorem it is enough to show
some Pieri formulas.
The statement corresponding to the Chevalley formula is well-known;
we prove the two other fundamental results in the general
context of Kac-Moody groups.

\subsection{Chevalley formula in the (co)minuscule case}

\label{subsection-chevalley}

Let $w \in W$ and $i \in I$ such that $l(s_{\alpha_i} w) = l(w) + 1$.
We denote with $m(w,i)$ the integer
$(\alpha_\Lambda , \alpha_\Lambda) / (\alpha_i,\alpha_i)$ if
$(\alpha_\Lambda , \alpha_\Lambda) > (\alpha_i,\alpha_i)$  and $m(w,i) = 1$
otherwise.

\begin{prop}
\label{prop-chevalley}
If $s_{\alpha_i}  w$ is length additive and $\Lambda$-(co)minuscule,
then the coefficient of the class 
$\s^{s_{\alpha_i} w}$ in the product $\s^w \cup \s^{s_{\alpha_\Lambda}}$
is $m(w,i)$.
\end{prop}
Thus, Conjecture~\ref{main_conj} is true when $u$ or $v$ has length one.
\begin{proo}
Recall the Chevalley formula
$$
\displaystyle{
\s^{s_{\alpha_\Lambda}} \cup \s^w = \sum_{\alpha:\ 
l(s_{\alpha} w) = l(w) + 1}
\scal{w(\Lambda),\alpha^\vee} \s^{s_{\alpha} w}.
}
$$
This follows from \cite[Theorem 11.1.7(i) and Remark
11.3.18]{kumar}. We only want to compute the coefficient of $\s^{s_\a w}$ in
$\s^{s_{\a_\Lambda}}\cup\s^w$ for $s_\a w$ a $\Lambda$-(co)minuscule element 
thus
we may in the sequel assume that $\a$ is simple (this comes from the
fact that weak and strong Bruhat order coincide for
$\Lambda$-(co)minuscule elements). 

Assume first that $s_{\alpha} w$ is $\Lambda$-minuscule. This means by
definition that $\scal{w(\Lambda) , \alpha^\vee} = 1$. Thus we only
have to prove that $(\alpha_\Lambda , \alpha_\Lambda) \leq
(\alpha,\alpha)$. This follows from Lemma
\ref{lemm-racines-courtes}.

Assume now that $s_{\alpha} w$ is $\Lambda$-cominuscule.
This means that $\scal{\alpha , w(\Lambda^\vee)} = 1$,
and therefore $\scal{w^{-1}(\alpha) , \Lambda^\vee} = 1$. By the
following Lemma 
\ref{lemm-chevalley} we have $\scal{\Lambda , w^{-1}(\alpha^\vee)} =
(\alpha_\Lambda , \alpha_\Lambda) / (\alpha,\alpha) $. Since
$s_{\alpha} w$ is $\Lambda$-cominuscule, by Lemma~\ref{lemm-racines-courtes}
the root $\alpha$ cannot be longer than $\alpha_\Lambda$ so this integer is
$m(w,i)$ and the proposition is proved.
\end{proo}

\begin{lemm}
\label{lemm-chevalley}
Let $\alpha,\beta$ be simple roots and $w \in W$. Then
$$
\scal { w(\alpha) , \varpi_\beta^\vee } \cdot ( \beta , \beta )
=
\scal { \varpi_\beta , w(\alpha^\vee) } \cdot ( \alpha , \alpha ).
$$
\end{lemm}
\begin{proo}
We prove this by induction on the length of $w$. If $w=e$, then both members of the equality equal
$(\alpha , \alpha)$ if $\alpha = \beta$ and 0 otherwise. Assume that
$$
\scal { w(\alpha) , \varpi_\beta^\vee } \cdot ( \beta , \beta )
=
\scal { \varpi_\beta , w(\alpha^\vee) } \cdot ( \alpha , \alpha ).
$$
and let $\gamma$ be a simple root. Since
$\scal { \varpi_\beta , \gamma^\vee }$ (resp. $\scal { \gamma , 
\varpi_\beta^\vee }$) is
by definition the coefficient of $\beta^\vee$ (resp. $\beta$) in 
$\gamma^\vee$ (resp. $\gamma$),
these coefficients are 1 if $\gamma = \beta$ and 0 otherwise. If 
$\gamma \not = \beta$, then
$\scal { s_\gamma w(\alpha) , \varpi_\beta^\vee } = \scal { w(\alpha) , 
\varpi_\beta^\vee }$
and $\scal { \varpi_\beta , s_\gamma w(\alpha^\vee) } = \scal { 
\varpi_\beta , w(\alpha^\vee) }$,
so the lemma is still true for $s_\gamma w$.
Moreover
$\scal { s_\beta w(\alpha) , \varpi_\beta^\vee } = \scal { w(\alpha) , 
\varpi_\beta^\vee } - \scal {w(\alpha) , \beta^\vee}$
and $\scal { \varpi_\beta , s_\beta w(\alpha^\vee) } = \scal { 
\varpi_\beta , w(\alpha^\vee) } - \scal {\beta , w(\alpha^\vee)}$.
Since $\scal {w(\alpha) , \beta^\vee} \cdot ( \beta , \beta ) = 
\scal {\beta , w(\alpha^\vee)} \cdot ( \alpha , \alpha ) 
= ( \alpha , \beta )$, the lemma is again true for $s_\beta \cdot w$.
\end{proo}

\subsection{Recursions}

\label{subsection-recursion}

Let us now introduce the notion of recursion, which is our essential
inductive argument, and was introduced in \cite{TY}.

\subsubsection{Homogeneous subspaces}

Let $G_1\subset G_2$ be an inclusion of Kac-Moody groups defined by 
an inclusion of their Dynkin diagrams (in particular we have an inclusion 
of the maximal torus $T_1$ of $G_1$ in the maximal torus $T_2$ of $G_2$). 
Let $\Lambda_2$ be a dominant weight for $G_2$ and $\Lambda_1$ its 
restriction to $T_1$.
We have an inclusion of the corresponding
Weyl groups $W_1 \subset W_2$ and of the homogeneous spaces
$G_1/P_1 \subset G_2/P_2$ where $P_i$ is associated to $\Lambda_i$ for 
$i\in\{1,2\}$.  

\begin{prop}
\label{coro-sous-groupe}
With the above notation, let $u$, $v$ and $w$ be elements in $W_1$ such that 
$u,v \leq w$. Assume that $w$ is
$\Lambda_1$-(co)minuscule. We have
$c_{u,v}^w(G_1/P_1) = c_{u,v}^w(G_2/P_2)$. 
Moreover we have
$t_{u,v}^w(W_1)m_{u,v}^w(W_1) = t_{u,v}^w(W_2)m_{u,v}^w(W_2)$.
\end{prop}
\begin{proo}
The claim for the coefficients $t$ and $m$ follows from the fact that 
the heap of $w$ does not depend on whether we consider $w$ as an element
of $W_1$ or $W_2$.

Let $i : G_1/P_1 \to G_2/P_2$ denote the natural inclusion. Observe
that $w$ (and thus also $u$ and $v$) is $\Lambda_2$-(co)minuscule.
To prove the proposition it is
enough to use the fact $i^*$ preserves the cup product: in fact,
we have the equality $i^*(\sigma^u(G_2/P_2)) = \sigma^u(G_1/P_1)$,
and thus the equality $i^*(\sigma^u(G_2/P_2) \cup \sigma^v(G_2/P_2)) =
\sigma^u(G_1/P_1) \cup \sigma^v(G_1/P_1)$ holds. Expanding these
products with the coefficients $c_{u,v}^w$ yields the result.
\end{proo}

Using this proposition, we see that that the coefficients $c_{u,v}^w(G/P)$
resp. $t_{u,v}^w(W)$ do not depend on $G/P$ resp. $W$, allowing us to 
simplify the notation into $c_{u,v}^w$ resp. $t_{u,v}^w$.

\begin{coro}
\label{coro-picard-superieur}
If Conjecture~\ref{main_conj}
holds when $P$ is a maximal parabolic subgroup,
then it holds in general.
\end{coro}
\begin{proo}
Let $u,v,w \in W$ and assume $w$ is $\Lambda$-(co)minuscule. Write
$\Lambda = \sum \Lambda_i \varpi_i$, with $\varpi_i$ the fundamental
weights. Let $S(\Lambda) \subset S$ be the set of indices
$i$ such that $\Lambda_i > 0$.
By \cite[Proposition page 65]{proctor-classification}
we can write $w$ as a commutative product
$w = \prod_{i \in S(\Lambda)} w_i$ where the supports of all the $w_i$'s
are disjoint and $i \in Supp(w_i)$. In the same way we write $u =
\prod_{i \in S(\Lambda)} u_i$ and $v = \prod_{i \in S(\Lambda)} v_i$. 
It follows that $m(w)=\prod m(w_i)$, that $m_{u,v}^w=\prod
m_{u_i,v_i}^{w_i}$ and that $t_{u,v}^w
=\prod {t}_{u_i,v_i}^{w_i}$. 
Moreover by Proposition~\ref{coro-sous-groupe} we have
$c_{u,v}^w = \prod c_{u_i,v_i}^{w_i}$.
Thus assuming that
$c_{u_i,v_i}^{w_i} = m_{u_i,v_i}^{w_i}\cdot t_{u_i,v_i}^{w_i}$ we get
$c_{u,v}^w = m_{u,v}^w \cdot t_{u,v}^w$.
\end{proo}

\subsubsection{Bruhat and taquin recursions}

\begin{defi}
\label{defi-recursion}
Let $x \in W$ be a $\Lambda$-(co)minuscule element.
\begin{itemize}
\item
Let $S(x) \subset S$ defined by $\a \in S(x)$ if and only if
$\scal{\alpha^\vee,x(\Lambda)} \geq 0$.
\item
Let $H_x \subset G$ be generated by the subgroups $SL_2(\alpha)$ of $G$ for
$\a \in S(x)$.
\item
Let $Q_x \subset H_x$ be the stabilisor of $[x]$ in $H_x$.
\item
Let $W_x \subset W$ be generated by the simple reflections $s_\a$ for $\a \in
S(x)$.
\item
We denote with $W_x \cdot x \subset W$ the subset of all elements of
the form $yx$ for some $y \in W_x$.
\end{itemize}
\end{defi}

\begin{fact}
\label{fait-px-parabolique}
$Q_x$ is a parabolic subgroup of $H_x$.
\end{fact}
\begin{proo}
Let $\alpha$ be a positive root of $H_x$. We can write 
$$\alpha = \sum_{i \in S(x)} n_i \alpha_i,$$
with $n_i \geq 0$. By definition of $S(x)$ it follows that
$\scal{\alpha^\vee,x(\Lambda)} \geq 0$. Since the set of weights
of the $SL_2(\alpha)$-representation
generated by the weight line $L_x$ of weight $x$ is the interval
$[x(\Lambda),s_\alpha(x(\Lambda))]$, it therefore contains weights of the 
form $x(\Lambda)-n\alpha$ with $n \geq 0$. Thus $\fg_{\alpha}$ acts 
trivially on $L_x$.
\end{proo}

Let $x$ be a $\Lambda$-(co)minuscule element and let $H(x)$ be its heap. We 
define the peaks of $H(x)$ to be the maximal elements in $H(x)$ with respect 
to the partial order (see \cite{small} for more combinatorics on these peaks 
and some geometric interpretations). Denote with ${\rm Peak}(x)$ the
set of peaks in $H(x)$. Recall that we denote with $c:H(x)\to D$ the
coloration of the heap.

\begin{prop}
We have $S(x)=S\setminus c({\rm Peak}(x))$.
\end{prop}

\begin{proo}
  Remark that it is enough to prove this statement for
  $\Lambda$-minuscule elements, the corresponding statement for 
$\Lambda$-cominuscule elements will
  follow by taking the dual root system.

Take $x=s_{\b_1}\cdots s_{\b_n}$ a reduced expression for $x$. We have
for any index $i\in[1,n-1]$ the equality $s_{\b_i}\cdots
s_{\b_n}(\Lambda)= s_{\b_{i+1}}\cdots s_{\b_n}(\Lambda)-\b_i$. If
$\a\in c({\rm  Peak}(x))$ we may assume that $\b_1=\a$ and we have
$s_{\a}(s_{\b_2}\cdots s_{\b_n}(\Lambda))= x(\Lambda)=s_{\b_2}\cdots
s_{\b_n}(\Lambda)-\a$. We get
$$\sca{\a^\vee}{x(\Lambda)}=\sca{\a^\vee}{s_{\b_2}\cdots
  s_{\b_n}(\Lambda)}-\sca{\a^\vee}{\a}=1-2=-1$$
therefore $c({\rm Peak}(x))$ does not meet $S(x)$. 

Now consider a
simple root $\a$ not in $c({\rm Peak}(x))$ and keep the reduced expression 
$x=s_{\b_1}\cdots s_{\b_n}$ for $x$. We have 
$$\scal{\a^\vee,x(\Lambda)}=\scal{\a^\vee,\Lambda}-\sum_{i=1}^n
\scal{\a^\vee,\beta_i}.$$
If $\a$ is not in the support of $x$, then for all $i$ we have 
$\scal{\a^\vee,\beta_i}\leq0$ thus $\scal{\a^\vee,x(\Lambda)}\geq0$ and 
$\a\in S(x)$. If $\a$ is in the support of $x$, then there exists an index 
$j$ with $\a=\b_j$. Because $\a$ is not in $c({\rm Peak}(x))$, there
exists an index $k<j$ such that for all $i\in[1,k-1]$ we have
$\sca{\b_i}{\a^\vee}=0$ and $\sca{\b_k}{\a^\vee}<0$. We may even assume 
that $\b_{k+1}=\a$ and we have $s_{\a}(s_{\b_{k+2}}\cdots s_{\b_n}(\Lambda))=
s_{\b_{k+2}}\cdots s_{\b_n}(\Lambda)-\a$ thus the equality
$$\sca{\a^\vee}{x(\Lambda)}=\sca{\a^\vee}{s_{\b_k}\cdots
  s_{\b_n}(\Lambda)} = \sca{\a^\vee}{s_{\b_{k+2}}\cdots
    s_{\b_n}(\Lambda)-\a-\b_k} = 1-2-\sca{\a^\vee}{\b_k}$$
holds. As $\sca{\a^\vee}{\b_k}<0$ we get
$\sca{\a^\vee}{x(\Lambda)}\geq0$. 
\end{proo}

Let $w$ be a $\Lambda$-(co)minuscule element with $w\geq x$ and denote
with $H(w)$ its heap.

\begin{coro}
\label{recu-carquois}
The element $w$ is in $W_x\cdot x$ as soon as $c(H(w)-H(x))\cap 
c({\rm Peak}(x))=\emptyset$.
\end{coro}

\begin{rema}
\label{cas-part}
In particular, we shall consider the special case of
recursion when $x$ is such that $c({\rm Peak}(x))$ consists of a
unique simple root: see Lemma~\ref{lemm-plus-grand}.
\end{rema}

\begin{defi}
Let $x \in W$. We say that $x$ is a Bruhat recursion
resp. a taquin recursion if for all $u,w \in W_x \cdot x$ with $u \leq w$
and $w$ a $\Lambda$-(co)minuscule element, and for all $v \leq w$,
the following holds:
$$
\begin{array}{c}
\displaystyle{c_{u,v}^w (G/P) = \sum_{s \in [e,wx^{-1}]}
  c_{ux^{-1},s}^{wx^{-1}} (H_x/Q_x) 
\cdot c_{x,v}^{sx} (G/P).}\\
\displaystyle{\mbox{ resp. } t_{u,v}^w (W)m_{u,v}^w (W) = 
\sum_{s \in [e,wx^{-1}]}
t_{ux^{-1},s}^{wx^{-1}} (W_x)m_{ux^{-1},s}^{wx^{-1}} (W_x)
\cdot t_{x,v}^{sx} (W)m_{x,v}^{sx} (W)\ .}
\end{array}
$$
\end{defi}

By \cite[Proposition 2.1]{stembridge} if $x \leq w$ and $w$ is 
$\Lambda$-(co)minuscule,
then $x$ is also $\Lambda$-(co)minuscule; thus if $x$ is not
$\Lambda$-(co)minuscule then the above statement is empty.

\begin{prop}
\label{recu}
Let $x \in W$ be $\Lambda$-(co)minuscule. Then $x$ is a taquin recursion.
\end{prop}
\begin{proo}
We start with the same formula involving only the taquin terms:
$$t_{u,v}^w (W)=
\sum_{s \in [e,wx^{-1}]}
t_{ux^{-1},s}^{wx^{-1}} (W_x)
\cdot t_{x,v}^{sx} (W).
$$
This formula was proved by Thomas and Yong in the more restrictive 
setting of cominuscule recursion (see \cite[Theorem 5.5]{TY}). Their proof
adapts here verbatim. 

We need to include the $m_{u,v}^w$ terms. For $u$ a $\Lambda$-minuscule 
element, we have, 
by Lemma~\ref{lemm-racines-courtes}, the equality $m(u)=1$ and the result 
follows.
For $u$ a $\Lambda$-cominuscule element, we may by Lemma 
\ref{lemm-racines-courtes} 
rewrite $m(u)$ as follows:
$$m(u)=\prod_{a\in H(u)}\frac{(\a_\Lambda,\a_\Lambda)}{({c(a)},{c(a)})}.$$
In particular we 
get for $m_{u,v}^w(W)$ an expression independent of $\a_\Lambda$ and thus 
independent of $W$. It only depends on the heaps of $u$, $v$ and $w$:
\begin{equation}
  \label{eq:equa_m}
  m_{u,v}^w(W)=\frac{\displaystyle{\prod_{a\in H(u)}({c(a)},
{c(a)})\prod_{a\in H(v)}({c(a)},
{c(a)})}}{\displaystyle{\prod_{a\in H(w)}({c(a)},
{c(a)})}}
\end{equation}
Now we remark that for $u'\in W_x$ with $u=u'x$, the heap $H(u)$ of $u$ is 
the union of the heaps $H(x)$ and $H(u')$.
In particular this 
gives $m(u)=m(ux^{-1}) m(x)$ so $m_{u,v}^w=m_{ux^{-1},s}^{wx^{-1}} 
m_{x,v}^{sx}$ and
the result follows.
\end{proo}

\subsubsection{A $\Lambda$-(co)minuscule element defines a Bruhat recursion}

Let us first prove a result on the length of elements of the form $wx$.

\begin{lemm}
\label{nicolas}
  Let $w\in (W_x)^{Q_x}$.

(\i) We have $wx\in W^P$.

(\i\i) We have $l(wx)=l(w)+l(x)$.
\end{lemm}

\begin{proo}
  Let us prove this result for a $\Lambda$-minuscule element
  first. The result for a $\Lambda$-cominuscule element follows
  since all these properties depend only on the Weyl group and thus
  not on the orientations of the arrows in the Dynkin diagram. Recall
  the characterisation
$$W^P=\{w\in W\ /\ w(\a)>0 \textrm{ for $\a>0$ with }
\sca{\Lambda}{\a^\vee}=0\}.$$
Recall also that for $u\in W^P$ we have $l(u)=\vert {\rm Inv}(u)\vert$ where
${\rm Inv}(u)$ is the set of inversions in $u$ defined by:
$${\rm Inv}(u)=\{\a>0\ /\ u(\a)<0\ {\rm and}\ \sca{\Lambda}{\a}>0\}.$$

$(\imath)$ Let $\a$ be a positive root with $\sca{\Lambda}{\a^\vee}=0$, we
need to prove that $wx(\a)$ is positive. Because $x\in W^P$, we have
$x(\a)>0$. Assume first that $x(\a)$ is a root of $H_x$, then we have
$\sca{x(\Lambda)}{x(\a)^\vee}=0$ thus, because $w\in W_x^{Q_x}$, we
have $w(x(\a))>0$. If $x(\a)$ is not a root of $H_x$, then it has a
positive coefficient on a simple root not in the root system of
$H_x$. But as $w\in W_x^{Q_x}$, the root $w(x(\a))$ has the same
coefficient on that root and $wx(\a)>0$.

($\imath\imath$) We have the inequality $l(wx)\leq l(w)+l(x)$. To prove the
converse inequality, we prove the following inclusion (and thus
equality) on the set of inversions:
$${\rm Inv}(x)\cup x({\rm Inv}(w))\subset {\rm Inv}(wx).$$
We will also prove that the first two sets are disjoint proving the
result. 

Let $\a$ a positive root with $\sca{\Lambda}{\a^\vee}>0$ and
$x(\a)<0$. Assume that $x(\a)$ is in the root system of $H_x$. We may
write $x(\a)$ has a linear combination of positive roots in $H_x$ with
non positive coefficients. Thus by definition of $H_x$, we get
$\sca{x(\Lambda)}{x(\a)^\vee}\leq0$. But we have the equality
$\sca{x(\Lambda)}{x(\a)^\vee}=\sca{\Lambda}{\a^\vee}>0$ a
contradiction. This implies, by the same argument as in the end of
($\imath$) that $wx(\a)<0$. Thus ${\rm Inv}(x)\subset {\rm Inv}(wx)$.

Let $\b$ a positive root of $H_x$ with $w(\b)<0$ and
$\sca{x(\Lambda)}{\b^\vee}>0$. We have
$\sca{\Lambda}{x^{-1}(\b)^\vee}>0$ thus $x^{-1}(\b)>0$ and
$x^{-1}(\b)\in {\rm Inv}(wx)$. The second inclusion follows. The sets 
${\rm Inv}(x)$
and $x^{-1}({\rm Inv}(w))$ are disjoint since by our proof $x({\rm Inv}(x))$ is
disjoint from the root system of $H_x$ while $x(x^{-1}({\rm Inv}(w)))=
{\rm Inv}(w)$ is
contained in that root system.
\end{proo}

Let $B$ be a Borel subgroup of $G$ and $U^-$ an opposite unipotent subgroup
(see \cite[Page 215]{kumar} for more details). Given $w \in W^P$ we denote
with $X_w$ resp. $X^w$ the closure of the $B$-orbit resp. $U^-$-orbit in $G/P$
through the point $wP/P$ in $G/P$. For $u \in W_x^{Q_x}$ we define
similarly the subvarieties $Y_u$ and $Y^u$ of $H_x/Q_x$. We also
denote with $i:H_x/Q_x \to G/P$ the natural injection.

\begin{lemm}
\label{lemm-intersection}
Let $x$ be $\Lambda$-(co)minuscule and let $u,w \in {(W_x)}^{Q_x}$.
We have $X^{ux} \cap X_{wx} = i ( Y^v \cap Y_u )$, as subvarieties
of $G/P$.
\end{lemm}
\begin{proo}
For $v \in W^P$ let $[v] \in G/P$ denote the corresponding
$T$-fixed point, and define similarly $[u] \in H_x/Q_x$ for
$u \in W_x^{Q_x}$. Let $U(x) \subset B$ resp.
$U(w) \subset B_x$ denote the unipotent subgroups corresponding
to $x$ resp. $w$.
We have $X_x = \overline {U(x) \cdot [e]}$ thus
$x \in \overline {U(x) \cdot [e]}$,
from which it follows that
$U(w) \cdot x \subset \overline {U(w)U(x) \cdot [e]} = X_{wx}$. Since
$i([e])=[x]$ and $i$ is $H_x$-equivariant, it follows that
$i(Y_w) \subset X_{wx}$. Similarly we have
$i(Y^u) \subset X^{ux}$. Thus we have an injection
$i:Y^u \cap Y_w \to X^{ux} \cap X_{wx}$.

By \cite[Lemma 7.3.10]{kumar}, both intersections are transverse
and irreducible, so that,
by Lemma \ref{nicolas}, the intersections
$Y^u \cap Y_w$ and $X^{ux} \cap_{wx}$ have the same dimension, namely
$l(w)-l(u)$, and thus the lemma is proved.
\end{proo}

For $u\in {(W_x)}^{Q_x}$, let us denote with $\tau_u$ resp. $\tau^u$ the
Schubert class in the homology group $H_*(H_x/Q_x,\Z)$ resp. its dual in
$H^*(H_x/Q_x,\Z)$. 

\begin{lemm}
\label{lemm-inter}
Let $x$ be $\Lambda$-(co)minuscule and let $u,w \in {(W_x)}^{Q_x}$.
We have
$\sigma^{ux} \cap \sigma_{wx} = i_* ( \tau^u \cap \tau_w )$,
in $H_*(G/P)$.
\end{lemm}
\begin{proo}
We still denote with $\sigma^{ux}$ the restriction of the cohomology
class $\sigma^{ux}$ to $X_{wx}$. We choose a reduced expression $\ww$ for $wx$
and denote with 
$q : \tilde X_{wx} \to X_{wx}$ the Bott-Samelson resolution associated to this
expression (see for example \cite[Chapter 7]{kumar}).
Recall that, since the expression is reduced, the morphism $q$ is birational.
We denote with $p$ its inverse which is a rational morphism. Observe that
$p$ is defined at $[wx]$.

Since $\tilde X_{wx}$ is smooth, homology and cohomology are identified
via Poincar{\'e} duality and moreover the cup product identifies with
the intersection product in the Chow ring.
We assume that $u \leq w$, since otherwise the terms of the lemma
both equal 0. In this case
$[wx] \in X^{ux} \cap X_{wx}$ and we
define $\tilde X^{ux} = \overline {p(X^{ux} \cap X_{wx})}$.
We claim that
$[\tilde X^{ux}] = q^* \sigma^{ux} \in H^*(\tilde X_{wx})$.
Note that $q^* \s^{ux}$ is caracterised by the equality
$\scal{q^* \s^{ux} , \g } = \scal { \s^{ux} , q_* \g } $.
To prove our claim, we use the fact that
$H_{2l(ux)}(\tilde X_{wx})$ has a basis consisting of the
classes $[\tilde X_\wv]$ where $\tilde X_\wv$ is the Bott-Samelson subvariety
of $\tilde X_{wx}$ defined by the subword $\wv$ of $\ww$ and the length of
$\wv$ is $l(ux)$. The claim is now implied by the fact that the intersection
$\tilde X^{ux} \cap \tilde X_\wv$ is a reduced point
if $q(X_\wv) = X_{ux}$ and is empty otherwise. Indeed, first remark that
$q(\tilde X_\wv)$ is a Schubert variety. We may thus use Lemma 7.1.22 and
Lemma 7.3.10 in \cite{kumar}. If $\dim q(X_\wv) < l(u)+l(x)$ then
$q(X_\wv)$ will not meet $X^{ux}$ and we are done. If $\dim q(X_\wv) =
l(u)+l(x)$, then $q(X_\wv)$ can meet $X^{ux}$ only if $q(X_\wv) = X_{ux}$, in
which case they meet transversely at $[ux]$. Moreover, since $p$ is defined
at $[ux]$, it follows that $\scal { \tilde X^{ux} ,\tilde X_\wv}=1$
in this case.

Remark that because $q$ is birational, we have the equality $q_*[\tilde
X_{wx}]=\s_{wx}$. Since furthermore $p$ is defined at $[wx]$, we have the
equality $q_*[\tilde X^{ux}\cap \tilde X_{wx}]=q_*[\tilde X^{ux}] =
[X^{ux}\cap X_{wx}]$. Applying projection formula we get: 
$$\s^{ux} \cap \s_{wx}=q_*(q^*
\s^{ux} \cap [\tilde X_{wx}])=q_*( [ \tilde X^{ux} \cap \tilde X_{wx}] ) =
[X^{ux} \cap X_{wx}].$$
The same argument gives $\tau^{u}\cap\tau_w=[Y^u\cap Y_w]$ and the lemma
follows from Lemma \ref{lemm-intersection}.
\end{proo}

\begin{theo}
\label{theo-bruhat-recursion}
Let $x$ be $\Lambda$-(co)minuscule, let $u,w \in {(W_x)}^{Q_x}$
and let $v \in W^P$.
Then we have
$$
\displaystyle{c_{ux,v}^{wx} (G/P) = \sum_{s \in [e,w]}
  c_{u,s}^{w} (H_x/Q_x)  \cdot c_{x,v}^{sx} (G/P).}
$$
\end{theo}
\noindent
In other words, $x$ is a Bruhat recursion.
\begin{proo}
The proof goes as in \cite{TY}.
Let $x,u,v,w \in W$ be as in the hypothesis of the theorem.

The left hand side of the equality in Lemma \ref{lemm-inter} is equal to
$\sum_v c_{ux,v}^{wx}(G/P) \sigma_v$, and the right hand side is 
equal to 
$i_* \sum_s c_{u,s}^w(H_x/Q_x) \tau_s$. By Lemma \ref{lemm-inter}
again, we have the equalities 
$i_* \tau_s = \sigma^x \cap \sigma_{sx} = \sum_v c_{x,v}^{sx}(G/P) \sigma_v$,
so the right
hand side is
$\sum_{v,s} c_{u,s}^w(H_x/Q_x) \cdot c_{x,v}^{sx}(G/P) \sigma_v$. Equating 
the coefficient of $\sigma_v$ we get the theorem.
\end{proo}

\subsection{System of posets associated with a dominant weight}
\label{subsection-posets}

Contrary to the situation of \cite{TY}, to compute the intersection numbers 
in a general homogeneous space, it will not be possible to use only one 
poset. Therefore it is necessary to show that the notion of ideals, of skew  
ideals, of tableaux, and of rectification make sense for a system of posets.

Let $J$ be a poset. A $J$-system $\pos$ of posets
is the data of a poset $P_i$ for each $i$ in $J$ and an injective morphism 
of posets $f_{i,j}:P_i \to P_j$ for all pairs $(i,j)$ with $i \leq j$, such 
that $f_{i,j}(P_i)$ is an
order ideal in $P_j$ and $f_{j,k} \circ f_{i,j} = f_{i,k}$ if $i \leq j \leq k$.
We assume that $J$ and each $P_i$'s are bounded below. Thus if
$\lambda \subset P_i$ is an order ideal and $i \leq j$
then $f_{i,j}(\lambda) \subset P_j$ is also an order ideal in $P_j$,
and we consider the order in the set
$S := \{ ( \lambda , P_i ) : \lambda \mbox{ is an  order ideal in } P_i \}$
generated by the relations
$(\lambda , P_i) \leq (f_{i,j}(\lambda) , P_j)$. The set of order
ideals of the system $\pos$ is by definition the direct limit of $S$. A skew 
ideal is a pair $(\nu,\lambda)$ of order ideals of $\pos$ such that
$\lambda \subset \nu$; it will be denoted with $\nu / \l$. A tableau
$T$ in $\pos$  of skew shape $\nu / \l$,
where $\nu / \l$ is a skew ideal, is a list of compatible tableaux in each of 
the $P_i$ where
$\nu$ is defined, of skew shape $\nu_i / \l_i$.

We say that $\pos$ has the jeu de taquin property if each $P_i$ has
this property. Let $T_i$ be a tableau of skew shape $\lambda / \nu$ in
$P_i$, let $i \leq j$, and denote with $T_j := f_{i,j}(T_i)$. If $R_i$
(resp. $R_j$) denotes the rectification of $T_i$ (resp. $T_j$) in
$P_i$ (resp. $P_j$), then note that $R_j = f_{i,j}(R_i)$ (informally,
the rectification of a tableau does not depend on what is above this
tableau). Therefore the rectification of a tableau in the system of
posets $\pos$ is well-defined as a tableau in $\pos$. Moreover an
analogue of Proposition~\ref{prop_independant_t} holds in this
context, thus defining the integer $t_{\lambda,\mu}^\nu$ for three
order ideals in $\pos$.

\vskip .4cm

Recall that $\Lambda$ is a dominant weight in a root system with Weyl
group $W$. We now show that $\Lambda$ defines a system of posets with
the jeu de taquin property. Let $J$ be the set of
$\Lambda$-(co)minuscule elements in $W$, equipped with the  weak
Bruhat order 
(which coincides with the strong Bruhat order). If $v,w \in J$ and $v
\leq w$, then we may write $w = s_{i_1} \cdots s_{i_k} \cdot v$, thus
the heap $H(v)$ of $v$ embeds naturally in $H(w)$ as an order ideal of
$H(w)$. This gives a map $f_{v,w}$ and defines the system
$\pos_\Lambda$ associated with $\Lambda$. Note that the set of order
ideals of $\pos_\Lambda$ is the set of heaps of
$\Lambda$-(co)minuscule elements in $W$. We refer to the pictures
(\ref{equa-systeme-d7}) in Subsection \ref{section-decomp} for
pictures of such posets.

\subsection{Algebra associated with a system of posets having the jeu
  de taquin property} 

\label{subsection-algebra-taquin}

Using the jeu de taquin, we now define a $\Z$-algebra $H(\pos)$
attached to any system of posets $\pos$ having the jeu de taquin property.
As a $\Z$-module, $H(\pos)$ is just
a free $\Z$-module with basis $\{x_\lambda\}$ indexed by all
order ideals $\lambda$ of $\pos$.
We then define a product on $H(\pos)$ by
$$ x_\lambda *_\pos 
x_\mu := \sum_\nu t_{\lambda,\mu}^\nu x_\nu,$$
where $t_{\lambda,\mu}^\nu$ is the integer defined in Proposition
\ref{prop_independant_t}.
If $T'$ is a tableau of skew shape $\nu / \l$, we denote with
$x_{T'} := x_\nu$ and say that $T'$ is {\it relative to $\l$}. We also write
$T' \leadsto T$ when the rectification of $T'$ is a standard tableau $T$.
Our definition
of the algebra $H(\pos)$
may thus be rewritten as
$ x_\lambda *_\pos 
x_\mu := \sum_{T' \leadsto T} x_{T'}$, where
the sum runs over all $T'$ relative to $\l$ and where $T$ is a fixed standard
tableau of shape $\mu$.

\begin{prop}
\label{prop_h_associatif}
Let $\pos$ be a system of posets having the jeu de taquin property. Then
the algebra $H(\pos)$ with the product $*_\pos$ is commutative and
associative.
\end{prop}
\begin{proo}
The commutativity of $H(\pos)$ amounts to the fact that
$t_{\l,\mu}^\nu = t_{\mu,\l}^\nu$, which is proved in Proposition
\ref{prop_independant_t}. Let us prove that $H(\pos)$ is associative.

So let $\lambda,\mu,\nu$ be order ideals. We choose standard 
tableaux $U$ and $V$, of shapes $\mu$ and $\nu$, and labelled respectively
with the indices
$\{1,\ldots,|\mu|\}$ and $\{|\mu|+1,\ldots,|\mu|+|\nu|\}$.
If $\gamma$ is a skew ideal, let $sh(\gamma)$ denote its shape.
By definition, we have
\begin{equation}
\label{equa_asso_1}
(x_\lambda *_\pos x_\mu) *_\pos x_\nu = \sum_{U' \leadsto U,V'' 
\leadsto V} x_{V''}
\end{equation}
where $U'$ is relative to $\lambda$ and $V''$ to $\lambda \cup sh(U')$.
Since by definition
we have $x_\mu *_\pos x_\nu = \sum_{V' \leadsto V} x_{V'}$, where $V'$
is relative to $\mu$, and since
for each such $V'$, $U\cup V'$ is a standard tableau, we also have
by definition
\begin{equation}
\label{equa_asso_2}
x_\lambda *_\pos (x_\mu *_\pos x_\nu) = \sum_{V' \leadsto V,W' \leadsto U 
\cup V'} x_{W'}
\end{equation}
where $V'$ is relative to $\mu$ and $W'$ is relative to $\lambda$.

We finish the proof of the proposition exhibiting
a bijection between the set of pairs $(U',W'')$
in (\ref{equa_asso_1}) and the set of pairs $(V',W')$
in (\ref{equa_asso_2}). We hope that the following scheme will help
following the argument (the order ideals $\l,\mu,\nu$ correspond to
the shapes: circle, rectangle, triangle). 

\vskip 0.1 cm

\centerline{
\psfrag{U}{\small{$U$}}
\psfrag{U'}{\small{$U'$}}
\psfrag{V}{\small{$V$}}
\psfrag{V'}{\small{$V'$}}
\psfrag{V''}{\small{$V''$}}
\psfrag{W}{\small{$W$}}
\psfrag{W'}{\small{$W'$}}
    \epsfxsize=2in
    \epsfysize=2in
    \epsfbox{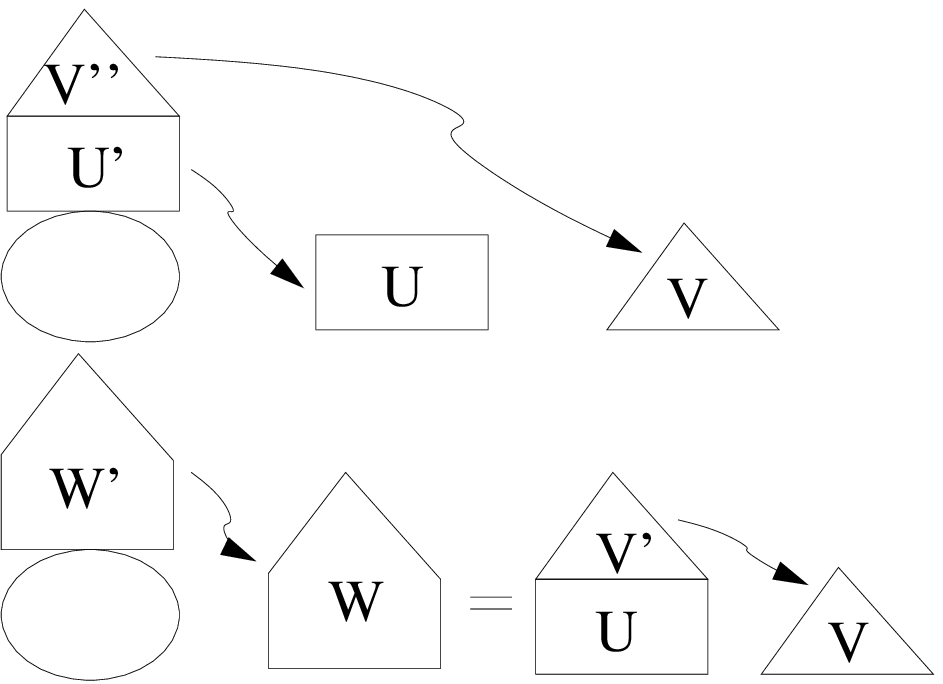}
    } 

\vskip 0.1 cm

Given a pair $(U',V'')$ as in (\ref{equa_asso_1}), we may consider the
standard
tableau $W'=U' \cup V''$. While performing the rectification of $W'$, we get
at each step a union of two tableaux which are obtained from $U'$ and $V''$
applying suitable jeu de taquin slides. At the end,
the rectification $W$ of $W'$
is a standard tableau $W = U_1 \cup T_1$, with $U_1$ (resp. $T_1$)
obtained by jeu de taquin slides from $U'$ (resp. $V''$). Therefore,
$U_1=U$, and $V_1$ rectifies to $V$. Therefore, if we set $V'=V_1$, we get
a pair $(V',W')$ in (\ref{equa_asso_2}). The inverse of this bijection is given
by setting $U'$ (resp. $V''$) to be the tableau made of all elements of $W$ 
with labels
less or equal to $|\mu|$ (resp. bigger than $|\mu|$). We thus have
proved that $(x_\lambda *_\pos x_\mu) *_\pos x_\nu = 
x_\lambda *_\pos (x_\mu *_\pos x_\nu)$.
\end{proo}

In the situation of a system of posets $\pos$ associated to a dominant 
weight $\Lambda$ as defined in Section~\ref{subsection-posets}, we define a 
perturbation  of this product by the numbers
$m_{\lt,\mu}^\nu$ as follows:
$$ x_\lambda \pp 
x_\mu := \sum_\nu t_{\lambda,\mu}^\nu m_{\lt,\mu}^\nu x_\nu.$$

Using Equation (\ref{eq:equa_m}) of the proof of Proposition \ref{recu},
we obtain:

\begin{coro}
Let $\pos$ be a system of posets having the jeu de taquin property. Then
the algebra $H(\pos)$ with the product $\pp$ is commutative and
associative.
\end{coro}

\vskip .5cm

We therefore have a purely combinatorially-defined algebra $H(\pos)$. On
the cohomology side there is also a natural algebra
with basis indexed by the $\Lambda$-minuscule (resp. 
$\Lambda$-cominuscule) elements of $W$,
because of the 
following fact (here we denote with $W_{mi}$ resp. $W_{co}$ the set of 
\lmin resp. $\Lambda$-cominuscule elements).

\begin{fact}
\label{fait-ideal}
The $\Z$-modules $\bigoplus_{w\not \in W_{mi} } \Z \cdot \sigma^w$ and 
$\bigoplus_{w\not \in W_{co} } \Z \cdot \sigma^w$ are ideals in $H^*(G/P)$.
\end{fact}
\begin{proo}
Let $v \in W$ be non $\Lambda$-(co)minuscule and let $x \in H^*(G/P)$.
We want to
show that $\sigma^v \cup x$ is a linear combinaison of some
$\sigma^w$'s with $w$ non $\Lambda$-(co)minuscule.
To this end we may assume that $x$
is a Schubert cohomology class of degree $d$; thus
$x \leq h^d$ ($h$ denotes the degree 1 Schubert cohomology class).

Write $\sigma^v \cup x = \sum c_w \sigma^w$,
and let $w$ be such that $c_w >0$.
Thus the coefficient of $\sigma^w$ in $\sigma^v \cdot h^d$ is positive.
Thus in the strong Bruhat order we have $v \leq w$. By
\cite[Proposition 2.1]{stembridge}, $w$ cannot be $\Lambda$-(co)minuscule.
So the fact is
proved.
\end{proo} 

\begin{fact}
\label{fait-troncation}
Let $w_1,\ldots,w_s \in W$. Then
the $\Z$-module $\bigoplus_{\forall i, w\not \leq w_i} \Z \cdot \sigma^w$
is an ideal in $H^*(G/P)$. We denote with $H^*_{(w_i)}(X)$ the
corresponding quotient algebra.
\end{fact}
\begin{proo}
For all $i\in[1,s]$, if $v \geq u$ and $u \not \leq w_i$, then
$v \not \leq w_i$. Thus the argument is the same as for the
previous fact.
\end{proo}

%%% Local Variables: 
%%% mode: latex
%%% TeX-master: t
%%% End: 

\section{Main result and strategy for the proof}

\label{section-general-argument}

\subsection{Statement of the main result}

Let $X=G/P$ be a homogeneous space and let $W$ resp. $\Lambda$ denote
the Weyl group of $G$ resp. the dominant weight associated to $P$. Denote 
with $D$ the Dynkin diagram of $G$.
Let $w \in W$ be $\Lambda$-(co)minuscule.
As in Definition \ref{defi-heap} we associate
to $w$ a heap $H(w)$. By \cite[Proposition A]{proctor-classification} 
(see also the end of Subsection \ref{subsection-taquin-poset}),
we may decompose
$H(w)$ into a so-called slant product of irreducible heaps that we
denote with 
$(H_i)_{0 \leq i \leq k}$.
We also denote with $D(H_i)=c(H_i)\subset D$ the Dynkin diagram 
corresponding to $H_i$.

\begin{defi}
\label{defi-slant}
Let $w \in W$ be $\Lambda$-(co)minuscule. We say that $w$ is
slant-finite-dimensional if all the Dynkin diagrams $D(H_i)$ are
Dynkin diagrams of finite-dimensional algebraic groups, in other words
$D(H_i)$ belongs to $\{ A_n,B_n,C_n,D_n,E_6,E_7,E_8,F_4,G_2 \}$ for
all $i$.
\end{defi}

\noindent
Our main result is the following.

\begin{theo}
\label{main-theo}
Let $G/P$ be a Kac-Moody homogeneous space where $P$ corresponds
to the dominant weight $\Lambda$. Let $u,v,w \in W$ be
$\Lambda$-(co)minuscule. Assume that $w$ is slant-finite-dimensional.
Then we have $c_{u,v}^w = m_{u,v}^wt_{u,v}^w$.
\end{theo}

\subsection{Definition of some systems of posets}
\label{section-decomp}

In order to prove Theorem \ref{main-theo}
we may assume, thanks to Corollary \ref{coro-sous-groupe},
that $P$ is a maximal parabolic subgroup of $G$.
The proof of Theorem \ref{main-theo} will be done by induction
on the rank of $G$, considering
the different possible cases for the irreducible component $H_0(w)$
of $H(w)$ containing the minimal element of $H(w)$.
In this subsection we give general lemmas to enable this.

\vskip .3cm

For the basic definitions concerning posets, we refer the reader to
Subsection \ref{subsection-taquin}.
We fix a marked Dynkin diagram $(D_0,\Lambda)$ which has no cycle,
and we consider a system of $\Lambda$-(co)minuscule $D_0$-colored
posets that we denote with $\pos_0$. We denote with $I_0$ the poset
indexing this system, so that for all $i \in I_0$ we are given a
$\Lambda$-(co)minuscule $D_0$-colored
poset $\pos_0(i)$. The choice of $\Lambda$ equips $D_0$ with the structure
of a poset, because we set $d_1 \leq d_2$ in $D_0$ if
$d_1$ and $\Lambda$ belong to the same
connected component of $D_0 - \{ d_2 \}$.
We assume that any $\alpha \in D_0$ is the color 
of at least one element in $\pos_0(i)$ for each $i$ in $I_0$,
thus the rooted tree of
$\pos_0(i)$ is equivalent, as a poset, with $D_0$.

We denote with $S_0$ the set
of maximal elements in $D_0$. For each
$\alpha \in S_0$ we suppose we are given a marked Dynkin
diagram $(D_\alpha,\Lambda_\alpha)$ and a $\Lambda_\alpha$-(co)minuscule
$D_\alpha$-colored poset $P_\a$, and
we now define a system of posets $\pos$ which contains all the possible
ways of adjoining the posets $P_\alpha$ to the posets $\pos_0(i)$.
Let $D$ be the Dynkin diagram obtained
from the disjoint union of $D_0$ and the $D_\alpha$'s for
$\alpha \in S_0$, where we connect $\alpha \in S_0$ with
$\Lambda_\alpha \in D_\a$ with an arbitrary number of edges.  
The colors of $\pos$ will be the elements of $D$.

The system $\pos$ is indexed
by the set of triples $(i,S_1,S_2)$ where $i \in I_0$ and
$(S_1,S_2)$ are subsets of $S_0$
with $S_1 \supset S_2$. This index set is itself a poset if we set
$(i,S_1,S_2) \leq (j,T_1,T_2)$ if $i \leq j$,
$S_1 \subset T_1$ and $S_2 \subset T_2$.

To any subset $S_1 \subset S_0$ and $i \in I_0$ we associate
the subposet $\pos_0(i,S_1)$ of $\pos_0(i)$
which is the maximal subposet such that all the
colors $\alpha$ in $S_0-S_1$ occur only once in $\pos_0(i,S_1)$ (in other
words $\pos_0(i,S_1)$ contains all the elements in
$\pos_0(i)$ which are not
bigger or equal to some element $(\alpha,2)$ with $\alpha \in S_0-S_1$).
Thus if $(i,S_1) \leq (j,T_1)$ then
$\pos_0(i,S_1) \subset \pos_0(j,T_1)$, and
$\pos_0(i,\emptyset) = \pos_0(i)$. 
We define $\pos(i,S_1,S_2)$ to be the slant product
of $\pos_0(i,S_1)$ and the posets $P_\alpha$ for $\alpha \in S_2$, where
the poset $P_\alpha$ is attached to $\pos_0(i,S_1)$ on the unique node
colored by $\alpha$ in $\pos_0(i,S_1)$. By
\cite[Proposition A]{proctor-classification}, $\pos(i,S_1,S_2)$
is $\Lambda$-minuscule (resp. $\Lambda$-cominuscule) if $\pos_0(i,S_1)$ 
is $\Lambda$-minuscule ($\Lambda$-cominuscule) and $\Lambda_\a$ is not 
shorter (resp. longer) that $\a$. Moreover for $(i,S_1,S_2) 
\leq (j,T_1,T_2)$ we obviously have an injection $\pos(i,S_1,S_2) 
\subset \pos(j,T_1,T_2)$, so that $\pos$ is indeed a system of 
$\Lambda$-(co)minuscule $D$-colored posets.

\begin{nota}
\label{nota-system}
We denote with $\pos_{\pos_0,(P_\alpha)}$ the system of posets
constructed above.
\end{nota}

\begin{exam}
  In the following array we give explicitly the obtained system of
posets when $\pos_0$ contains only one element which is the heap
of the maximal Schubert cell in $D_7/P_6$.
Note that in this case $S_0 = \{1,7\}$.
Since $I_0$ has only one element we abbreviate $\pos(i,S_1,S_2)$ into
$\pos(S_1,S_2)$. In the drawings we represent the rooted tree with
solid dots and solid diamonds (for the maximal elements), we represent
the elements which must belong to an ideal in order for this ideal to be
slant-irreducible with $\otimes$, and the other elements are depicted
with hollow dots. The posets $P_\a$ for $\a \in S_0$ are represented by
angular sectors.

\begin{equation}
\label{equa-systeme-d7}
\begin{array}{ccc}
%auto-ignore

\psset{unit=6mm}

\begin{pspicture*}(0.5,0.5)(-5.5,-10.5)

\psline(-0.12,-9.88)(-0.88,-9.12)
\psline(-1,-8.84)(-1,-8.16)
\psline(-1.12,-8.88)(-1.88,-8.12)
\psline(-1,-7.84)(-1,-7.16)
\psline(-1.88,-7.88)(-1.12,-7.12)
\psline(-2.12,-7.88)(-2.88,-7.12)
\psline(-0.88,-6.88)(-0.12,-6.12)
\psline(-1.12,-6.88)(-1.88,-6.12)
\psline(-2.88,-6.88)(-2.12,-6.12)
\psline(-3.12,-6.88)(-3.88,-6.12)
\psline(-0.12,-5.88)(-0.88,-5.12)
\psline(-1.88,-5.88)(-1.12,-5.12)
\psline(-2.12,-5.88)(-2.88,-5.12)
\psline(-3.88,-5.88)(-3.12,-5.12)
\psline(-4.12,-5.88)(-4.88,-5.12)
\psline(-1,-4.84)(-1,-4.16)
\psline(-1.12,-4.88)(-1.88,-4.12)
\psline(-2.88,-4.88)(-2.12,-4.12)
\psline(-3.12,-4.88)(-3.88,-4.12)
\psline(-4.88,-4.88)(-4.12,-4.12)
\psline(-1,-3.84)(-1,-3.16)
\psline(-1.88,-3.88)(-1.12,-3.12)
\psline(-2.12,-3.88)(-2.88,-3.12)
\psline(-3.88,-3.88)(-3.12,-3.12)
\psline(-0.88,-2.88)(-0.12,-2.12)
\psline(-1.12,-2.88)(-1.88,-2.12)
\psline(-2.88,-2.88)(-2.12,-2.12)
\psline(-0.12,-1.88)(-0.88,-1.12)
\psline(-1.88,-1.88)(-1.12,-1.12)
\psline(-1,-0.84)(-1,-0.16)

\pscircle*(-0,-10){0.16}
\pscircle(-0,-6){0.16}  \psline(0.12,-5.88)(-0.12,-6.12)  \psline(0.12,-6.12)(-0.12,-5.88)
\pscircle(-0,-2){0.16}
\pscircle*(-1,-9){0.16}
\psdiamond[fillstyle=solid,fillcolor=black](-1,-8)(0.256,0.256)
\pscircle(-1,-7){0.16}  \psline(-0.88,-6.88)(-1.12,-7.12)  \psline(-0.88,-7.12)(-1.12,-6.88)
\pscircle(-1,-5){0.16}
\pscircle(-1,-4){0.16}
\pscircle(-1,-3){0.16}
\pscircle(-1,-1){0.16}
\pscircle(-1,-0){0.16}
\pscircle*(-2,-8){0.16}
\pscircle(-2,-6){0.16}  \psline(-1.88,-5.88)(-2.12,-6.12)  \psline(-1.88,-6.12)(-2.12,-5.88)
\pscircle(-2,-4){0.16}
\pscircle(-2,-2){0.16}
\pscircle*(-3,-7){0.16}
\pscircle(-3,-5){0.16}  \psline(-2.88,-4.88)(-3.12,-5.12)  \psline(-2.88,-5.12)(-3.12,-4.88)
\pscircle(-3,-3){0.16}
\pscircle*(-4,-6){0.16}
\pscircle(-4,-4){0.16}  \psline(-3.88,-3.88)(-4.12,-4.12)  \psline(-3.88,-4.12)(-4.12,-3.88)
\psdiamond[fillstyle=solid,fillcolor=black](-5,-5)(0.256,0.256)

\end{pspicture*} & %auto-ignore

\psset{unit=6mm}

\begin{pspicture*}(0.5,0.5)(-6.7,-10.5)

\psline(-0.12,-9.88)(-0.88,-9.12)
\psline(-1,-8.84)(-1,-8.16)
\psline(-1.12,-8.88)(-1.88,-8.12)
\psline(-1,-7.84)(-1,-7.16)
\psline(-1.88,-7.88)(-1.12,-7.12)
\psline(-2.12,-7.88)(-2.88,-7.12)
\psline(-0.88,-6.88)(-0.12,-6.12)
\psline(-1.12,-6.88)(-1.88,-6.12)
\psline(-2.88,-6.88)(-2.12,-6.12)
\psline(-3.12,-6.88)(-3.88,-6.12)
\psline(-0.12,-5.88)(-0.88,-5.12)
\psline(-1.88,-5.88)(-1.12,-5.12)
\psline(-2.12,-5.88)(-2.88,-5.12)
\psline(-3.88,-5.88)(-3.12,-5.12)
\psline(-4.12,-5.88)(-4.88,-5.12)
\psline(-1,-4.84)(-1,-4.16)
\psline(-1.12,-4.88)(-1.88,-4.12)
\psline(-2.88,-4.88)(-2.12,-4.12)
\psline(-3.12,-4.88)(-3.88,-4.12)
\psline(-4.88,-4.88)(-4.12,-4.12)
\psline(-1,-3.84)(-1,-3.16)
\psline(-1.88,-3.88)(-1.12,-3.12)
\psline(-2.12,-3.88)(-2.88,-3.12)
\psline(-3.88,-3.88)(-3.12,-3.12)
\psline(-0.88,-2.88)(-0.12,-2.12)
\psline(-1.12,-2.88)(-1.88,-2.12)
\psline(-2.88,-2.88)(-2.12,-2.12)
\psline(-0.12,-1.88)(-0.88,-1.12)
\psline(-1.88,-1.88)(-1.12,-1.12)
\psline(-1,-0.84)(-1,-0.16)

\pscircle*(-0,-10){0.16}
\pscircle(-0,-6){0.16}  \psline(0.12,-5.88)(-0.12,-6.12)  \psline(0.12,-6.12)(-0.12,-5.88)
\pscircle(-0,-2){0.16}
\pscircle*(-1,-9){0.16}
\psdiamond[fillstyle=solid,fillcolor=black](-1,-8)(0.256,0.256)
\pscircle(-1,-7){0.16}  \psline(-0.88,-6.88)(-1.12,-7.12)  \psline(-0.88,-7.12)(-1.12,-6.88)
\pscircle(-1,-5){0.16}
\pscircle(-1,-4){0.16}
\pscircle(-1,-3){0.16}
\pscircle(-1,-1){0.16}
\pscircle(-1,-0){0.16}
\pscircle*(-2,-8){0.16}
\pscircle(-2,-6){0.16}  \psline(-1.88,-5.88)(-2.12,-6.12)  \psline(-1.88,-6.12)(-2.12,-5.88)
\pscircle(-2,-4){0.16}
\pscircle(-2,-2){0.16}
\pscircle*(-3,-7){0.16}
\pscircle(-3,-5){0.16}  \psline(-2.88,-4.88)(-3.12,-5.12)  \psline(-2.88,-5.12)(-3.12,-4.88)
\pscircle(-3,-3){0.16}
\pscircle*(-4,-6){0.16}
\pscircle(-4,-4){0.16}  \psline(-3.88,-3.88)(-4.12,-4.12)  \psline(-3.88,-4.12)(-4.12,-3.88)
\psdiamond[fillstyle=solid,fillcolor=black](-5,-5)(0.256,0.256)\pswedge[fillstyle=solid,fillcolor=blue,linecolor=blue](-5,-5){1.6}{130}{200}

\end{pspicture*} & %auto-ignore

\psset{unit=6mm}

\begin{pspicture*}(0.5,0.5)(-5.5,-7.5)

\psline(-0.12,-6.88)(-0.88,-6.12)
\psline(-1,-5.84)(-1,-5.16)
\psline(-1.12,-5.88)(-1.88,-5.12)
\psline(-1,-4.84)(-1,-4.16)
\psline(-1.88,-4.88)(-1.12,-4.12)
\psline(-2.12,-4.88)(-2.88,-4.12)
\psline(-0.88,-3.88)(-0.12,-3.12)
\psline(-1.12,-3.88)(-1.88,-3.12)
\psline(-2.88,-3.88)(-2.12,-3.12)
\psline(-3.12,-3.88)(-3.88,-3.12)
\psline(-0.12,-2.88)(-0.88,-2.12)
\psline(-1.88,-2.88)(-1.12,-2.12)
\psline(-2.12,-2.88)(-2.88,-2.12)
\psline(-3.88,-2.88)(-3.12,-2.12)
\psline(-4.12,-2.88)(-4.88,-2.12)
\psline(-1.12,-1.88)(-1.88,-1.12)
\psline(-2.88,-1.88)(-2.12,-1.12)
\psline(-3.12,-1.88)(-3.88,-1.12)
\psline(-4.88,-1.88)(-4.12,-1.12)
\psline(-2.12,-0.88)(-2.88,-0.12)
\psline(-3.88,-0.88)(-3.12,-0.12)

\pscircle*(-0,-7){0.16}
\pscircle(-0,-3){0.16}  \psline(0.12,-2.88)(-0.12,-3.12)  \psline(0.12,-3.12)(-0.12,-2.88)
\pscircle*(-1,-6){0.16}
\psdiamond[fillstyle=solid,fillcolor=black](-1,-5)(0.256,0.256)
\pscircle(-1,-4){0.16}  \psline(-0.88,-3.88)(-1.12,-4.12)  \psline(-0.88,-4.12)(-1.12,-3.88)
\pscircle(-1,-2){0.16}
\pscircle*(-2,-5){0.16}
\pscircle(-2,-3){0.16}  \psline(-1.88,-2.88)(-2.12,-3.12)  \psline(-1.88,-3.12)(-2.12,-2.88)
\pscircle(-2,-1){0.16}
\pscircle*(-3,-4){0.16}
\pscircle(-3,-2){0.16}  \psline(-2.88,-1.88)(-3.12,-2.12)  \psline(-2.88,-2.12)(-3.12,-1.88)
\pscircle(-3,-0){0.16}
\pscircle*(-4,-3){0.16}
\pscircle(-4,-1){0.16}  \psline(-3.88,-0.88)(-4.12,-1.12)  \psline(-3.88,-1.12)(-4.12,-0.88)
\psdiamond[fillstyle=solid,fillcolor=black](-5,-2)(0.256,0.256)

\end{pspicture*} \\
\pos(\emptyset,\emptyset) = \pos(\{1\},\emptyset) & \pos(\{1\},\{1\}) &
\pos(\{7\},\emptyset) = \pos(\{1,7\},\emptyset) \\
\end{array}
\end{equation}
$$\begin{array}{ccc}
%auto-ignore

\psset{unit=6mm}

\begin{pspicture*}(0.5,0.5)(-6.7,-7.5)

\psline(-0.12,-6.88)(-0.88,-6.12)
\psline(-1,-5.84)(-1,-5.16)
\psline(-1.12,-5.88)(-1.88,-5.12)
\psline(-1,-4.84)(-1,-4.16)
\psline(-1.88,-4.88)(-1.12,-4.12)
\psline(-2.12,-4.88)(-2.88,-4.12)
\psline(-0.88,-3.88)(-0.12,-3.12)
\psline(-1.12,-3.88)(-1.88,-3.12)
\psline(-2.88,-3.88)(-2.12,-3.12)
\psline(-3.12,-3.88)(-3.88,-3.12)
\psline(-0.12,-2.88)(-0.88,-2.12)
\psline(-1.88,-2.88)(-1.12,-2.12)
\psline(-2.12,-2.88)(-2.88,-2.12)
\psline(-3.88,-2.88)(-3.12,-2.12)
\psline(-4.12,-2.88)(-4.88,-2.12)
\psline(-1.12,-1.88)(-1.88,-1.12)
\psline(-2.88,-1.88)(-2.12,-1.12)
\psline(-3.12,-1.88)(-3.88,-1.12)
\psline(-4.88,-1.88)(-4.12,-1.12)
\psline(-2.12,-0.88)(-2.88,-0.12)
\psline(-3.88,-0.88)(-3.12,-0.12)

\pscircle*(-0,-7){0.16}
\pscircle(-0,-3){0.16}  \psline(0.12,-2.88)(-0.12,-3.12)  \psline(0.12,-3.12)(-0.12,-2.88)
\pscircle*(-1,-6){0.16}
\psdiamond[fillstyle=solid,fillcolor=black](-1,-5)(0.256,0.256)
\pscircle(-1,-4){0.16}  \psline(-0.88,-3.88)(-1.12,-4.12)  \psline(-0.88,-4.12)(-1.12,-3.88)
\pscircle(-1,-2){0.16}
\pscircle*(-2,-5){0.16}
\pscircle(-2,-3){0.16}  \psline(-1.88,-2.88)(-2.12,-3.12)  \psline(-1.88,-3.12)(-2.12,-2.88)
\pscircle(-2,-1){0.16}
\pscircle*(-3,-4){0.16}
\pscircle(-3,-2){0.16}  \psline(-2.88,-1.88)(-3.12,-2.12)  \psline(-2.88,-2.12)(-3.12,-1.88)
\pscircle(-3,-0){0.16}
\pscircle*(-4,-3){0.16}
\pscircle(-4,-1){0.16}  \psline(-3.88,-0.88)(-4.12,-1.12)  \psline(-3.88,-1.12)(-4.12,-0.88)
\psdiamond[fillstyle=solid,fillcolor=black](-5,-2)(0.256,0.256)\pswedge[fillstyle=solid,fillcolor=blue,linecolor=blue](-5,-2){1.6}{130}{200}

\end{pspicture*} & %auto-ignore

\psset{unit=6mm}

\begin{pspicture*}(0.7,0.5)(-5.5,-7.5)

\psline(-0.12,-6.88)(-0.88,-6.12)
\psline(-1,-5.84)(-1,-5.16)
\psline(-1.12,-5.88)(-1.88,-5.12)
\psline(-1,-4.84)(-1,-4.16)
\psline(-1.88,-4.88)(-1.12,-4.12)
\psline(-2.12,-4.88)(-2.88,-4.12)
\psline(-0.88,-3.88)(-0.12,-3.12)
\psline(-1.12,-3.88)(-1.88,-3.12)
\psline(-2.88,-3.88)(-2.12,-3.12)
\psline(-3.12,-3.88)(-3.88,-3.12)
\psline(-0.12,-2.88)(-0.88,-2.12)
\psline(-1.88,-2.88)(-1.12,-2.12)
\psline(-2.12,-2.88)(-2.88,-2.12)
\psline(-3.88,-2.88)(-3.12,-2.12)
\psline(-4.12,-2.88)(-4.88,-2.12)
\psline(-1.12,-1.88)(-1.88,-1.12)
\psline(-2.88,-1.88)(-2.12,-1.12)
\psline(-3.12,-1.88)(-3.88,-1.12)
\psline(-4.88,-1.88)(-4.12,-1.12)
\psline(-2.12,-0.88)(-2.88,-0.12)
\psline(-3.88,-0.88)(-3.12,-0.12)

\pscircle*(-0,-7){0.16}
\pscircle(-0,-3){0.16}  \psline(0.12,-2.88)(-0.12,-3.12)  \psline(0.12,-3.12)(-0.12,-2.88)
\pscircle*(-1,-6){0.16}
\psdiamond[fillstyle=solid,fillcolor=black](-1,-5)(0.256,0.256)\pswedge[fillstyle=solid,fillcolor=blue,linecolor=blue](-1,-5){1.6}{-20}{50}
\pscircle(-1,-4){0.16}  \psline(-0.88,-3.88)(-1.12,-4.12)  \psline(-0.88,-4.12)(-1.12,-3.88)
\pscircle(-1,-2){0.16}
\pscircle*(-2,-5){0.16}
\pscircle(-2,-3){0.16}  \psline(-1.88,-2.88)(-2.12,-3.12)  \psline(-1.88,-3.12)(-2.12,-2.88)
\pscircle(-2,-1){0.16}
\pscircle*(-3,-4){0.16}
\pscircle(-3,-2){0.16}  \psline(-2.88,-1.88)(-3.12,-2.12)  \psline(-2.88,-2.12)(-3.12,-1.88)
\pscircle(-3,-0){0.16}
\pscircle*(-4,-3){0.16}
\pscircle(-4,-1){0.16}  \psline(-3.88,-0.88)(-4.12,-1.12)  \psline(-3.88,-1.12)(-4.12,-0.88)
\psdiamond[fillstyle=solid,fillcolor=black](-5,-2)(0.256,0.256)

\end{pspicture*} & %auto-ignore

\psset{unit=6mm}

\begin{pspicture*}(0.7,0.5)(-6.7,-7.5)

\psline(-0.12,-6.88)(-0.88,-6.12)
\psline(-1,-5.84)(-1,-5.16)
\psline(-1.12,-5.88)(-1.88,-5.12)
\psline(-1,-4.84)(-1,-4.16)
\psline(-1.88,-4.88)(-1.12,-4.12)
\psline(-2.12,-4.88)(-2.88,-4.12)
\psline(-0.88,-3.88)(-0.12,-3.12)
\psline(-1.12,-3.88)(-1.88,-3.12)
\psline(-2.88,-3.88)(-2.12,-3.12)
\psline(-3.12,-3.88)(-3.88,-3.12)
\psline(-0.12,-2.88)(-0.88,-2.12)
\psline(-1.88,-2.88)(-1.12,-2.12)
\psline(-2.12,-2.88)(-2.88,-2.12)
\psline(-3.88,-2.88)(-3.12,-2.12)
\psline(-4.12,-2.88)(-4.88,-2.12)
\psline(-1.12,-1.88)(-1.88,-1.12)
\psline(-2.88,-1.88)(-2.12,-1.12)
\psline(-3.12,-1.88)(-3.88,-1.12)
\psline(-4.88,-1.88)(-4.12,-1.12)
\psline(-2.12,-0.88)(-2.88,-0.12)
\psline(-3.88,-0.88)(-3.12,-0.12)

\pscircle*(-0,-7){0.16}
\pscircle(-0,-3){0.16}  \psline(0.12,-2.88)(-0.12,-3.12)  \psline(0.12,-3.12)(-0.12,-2.88)
\pscircle*(-1,-6){0.16}
\psdiamond[fillstyle=solid,fillcolor=black](-1,-5)(0.256,0.256)\pswedge[fillstyle=solid,fillcolor=blue,linecolor=blue](-1,-5){1.6}{-20}{50}
\pscircle(-1,-4){0.16}  \psline(-0.88,-3.88)(-1.12,-4.12)  \psline(-0.88,-4.12)(-1.12,-3.88)
\pscircle(-1,-2){0.16}
\pscircle*(-2,-5){0.16}
\pscircle(-2,-3){0.16}  \psline(-1.88,-2.88)(-2.12,-3.12)  \psline(-1.88,-3.12)(-2.12,-2.88)
\pscircle(-2,-1){0.16}
\pscircle*(-3,-4){0.16}
\pscircle(-3,-2){0.16}  \psline(-2.88,-1.88)(-3.12,-2.12)  \psline(-2.88,-2.12)(-3.12,-1.88)
\pscircle(-3,-0){0.16}
\pscircle*(-4,-3){0.16}
\pscircle(-4,-1){0.16}  \psline(-3.88,-0.88)(-4.12,-1.12)  \psline(-3.88,-1.12)(-4.12,-0.88)
\psdiamond[fillstyle=solid,fillcolor=black](-5,-2)(0.256,0.256)\pswedge[fillstyle=solid,fillcolor=blue,linecolor=blue](-5,-2){1.6}{130}{200}

\end{pspicture*} \\
\pos(\{1,7\},\{1\}) & \pos(\{7\},\{7\}) = \pos(\{1,7\},\{7\}) &
\pos(\{1,7\},\{1,7\})
\end{array}$$
\end{exam}

We also consider the Kac-Moody
homogeneous space defined by the marked
Dynkin diagram $(D_0,\Lambda)$ resp. $(D,\Lambda)$, that we denote
with $X_0$
resp. $X$.
Let $W$ be the Weyl group corresponding to $D$, and for each triple
$(i,S_1,S_2)$ let $w_{i,S_1,S_2} \in W$ be the $\Lambda$-(co)minuscule
element corresponding to
the $\Lambda$-(co)minuscule poset $\pos(i,S_1,S_2)$. We denote with 
$H^*_t(X)$ the
truncation $H^*_{\{w_{i,S_1,S_2}\}}(X)$ of $H^*(X)$ obtained with the elements
$w_{i,S_1,S_2}$ (see Fact \ref{fait-troncation}).

\subsection{Reduction to indecomposable posets}

In the rest of this subsection $\pos_0,S_0,(P_\alpha),\pos,X,H^*_t(X)$
are as above, and we assume that Conjecture \ref{main_conj} holds
for any marked Dynkin diagram $(D',d')$ and any
$D'$-colored $d'$-minuscule poset as soon as $D' \varsubsetneq D$.

We will give some lemmas which help
comparing
$H^*(\pos)$ with $H^*_t(X)$. Note that these two $\Z$-modules have a basis
indexed by the same set, namely the set of ideals of $\pos$. Thus, in
order to simplify notation, we will identify these $\Z$-modules
and denote with $x \cdot y$ resp. $x \pp y$
the product in $H^*_t(X)$ resp. $H^*(\pos)$.

\vskip .8cm

\begin{nota}
For $\alpha \in S_0$ we denote with 
$\lt_\alpha = \scal { (\alpha,1) }$ the ideal in $\pos_0$,
and we define the cohomology class $\sigma^\alpha =
\sigma^{\lt_\alpha} \in H^*(\pos_0)$.
\end{nota}

We now make use of Theorem \ref{theo-bruhat-recursion}.

\begin{lemm}
\label{lemm-plus-grand}
Let $\s \in H^*(\pos)$.

\begin{enumerate}
\item
Let $\a \in D$ and let $i$ be an integer. Assume that
$\ok{\s}{\s^\lt}$ for $\lt = \scal{(\a,i)}$. Then
$\oknu{\s}{\mu}{\nu}$ for $\mu,\nu \in I(\pos)$ such that
$(\a,i) \in \mu$ and $(\a,i+1) \not \in \nu$.

\item
In particular, assume that $\a$ and $i$ are such that
for each poset $P$ in the system $\pos$
the number of elements of $P$ colored by $\a$ is not bigger than $i$,
and that $\ok{\s}{\s^\lt}$ for $\lt = \scal{(\a,i)}$.
Then $\ok{\s}{\s^\mu}$ if $(\a,i) \in \mu$.

\item
In particular, if $\a \in S_0$ and $\ok{\s}{\s^\a}$, then
$\ok{\s}{\s^\l}$ for 
$\l$ containing $(\a,1)$.
\end{enumerate}
\end{lemm}
\begin{proo}
Let $\s \in H^*(\pos)$ and let $\a,i$ as in the first point, and
let $\mu,\nu \in I(\pos)$ such that
$\mu \supset \scal{(\a,i)}$ and $(\a,i+1) \not \in \nu$.

Let $x$ resp $u,w$ be the elements in $W$ corresponding to the ideals
$\scal{(\a,i)}$ resp. $\mu,\nu$.
Since $\scal{(\a,i)}$ has only one peak namely $(\a,i)$,
by Corollary \ref{recu-carquois} and the assumption on $\l,\mu$
we have $u,w \in W_x \cdot x$. By assumption, Conjecture \ref{main_conj}
holds for posets colored by $D-\{\a\}$. Thus for $s \in W_x$ we have
$\oknu{ux^{-1}}{s}{wx^{-1}}$. Moreover the hypothesis that
$\ok{\s}{\s^x}$ says that
$\oknu{x}{\s}{sx}$ for $s \in W_x$. Thus by Theorem
\ref{theo-bruhat-recursion}
it follows that $\oknu{u}{\s}{w}$. This proves the first point.

The second point follows because under the hypothesis
for any ideal $\mu$ we have
$(\a,i+1) \not \in \mu$. The third point is a special case, for $i=1$,
of the second one.
\end{proo}

\begin{lemm}
\label{lemm-egalite}
Let $\g\in H^*(\pos)$ such that $\ok{\g}{\s}$ for $\s\in H^*(\pos_0)$, 
then $\ok{\g}{\s}$ for $\s\in H^*(\pos)$.
\end{lemm}

\begin{proo}
Assume $\s=\s^\lt$. For $\lt\in I(\pos)-I(\pos_0)$, 
there exists $\a\in S_0$ such that $\lt\supset\lt^\a$. But the class 
$\s^\a$ is in $H^*(\pos_0)$ thus the 
equality $\ok{\g}{\s^\a}$ holds and by Lemma \ref{lemm-plus-grand}
we have $\ok{\g}{\s}$.
\end{proo}

In the following lemma $(\g^i)$ is a list of elements
in $H^*(\pos_0)$ and we denote with $\scal{(\g^i)}$ 
the subalgebra they generate in $H^*(\pos_0)$.
For $d$ an integer we denote with $\scal{(\g^i)}_d$ the classes of
$\scal{(\g^i)}$ of degree at most $d$.
Moreover we denote with 
$\pi:H^*(\pos) \to H^*(\pos_0)$ the algebra morphism obtained by moding
out by the ideal of $H^*(\pos)$ linearly
generated by the Schubert classes $\s^\l$
with $\l \in I(\pos) - I(\pos_0)$. Finally let $H_d \subset H^*(\pos)$
denote the space of linear combinaisons of $\s^\l$ for
$\l \in I(\pos) - I(\pos_0)$ such that there exists 
$\a \in S_0$ with $\deg(\s^\a) \leq d$ and $\l \supset \lt_\a$.

\begin{lemm}
\label{lemm-sous-algebre-d}
Let $(\g^i)_{i\in[1,k]}$ be elements in $H^*(\pos_0)$ and $d$ an integer.
Assume that
\begin{itemize}
\item
For all $i$
and all $\s \in H^*(\pos_0)$ we have
$\ok{\s}{\g^i}$.
\item
For each $\a$ in $S_0$ with $\deg(\s^\a) \leq d$,
we have $\s^\a \in \scal{(\g^i)}$.
\end{itemize}
Then for all $\s$ in $H^*(\pos)$ and for all
$\tau$ in $H^*(\pos)$ such that
$\pi(\tau) \in \scal{(\g^i)}_d$ and $\tau - \pi(\tau) \in H_d$
we have $\ok{\s}{\tau}$.
\end{lemm}

\begin{proo}
By Lemma \ref{lemm-egalite} we have the equality
$\ok{\s}{\g^i}$ for general $\s \in H^*(\pos)$.
In particular a polynomial expression in the $\g^i$'s
is the same whether it is computed with the product $\cdot$ or
$\pp$.
If $P$ is a polynomial and $\s \in H^*(\pos)$
we moreover have $\ok{\s}{P(\g^i)}$.

We then prove by induction on $d' \leq d$ that if $\a \in S_0$
with $\deg(\s^\a) \leq d'$
and $\s \in H^*(\pos)$, then
$$
\ok{\s}{\s^\l} \mbox{ if } \l \supset \lt_\a.
$$
Let $d' \leq d$ be an integer and let
$\a$ such that $\deg(\s^\a) = d'$.
Let $P$ be a polynomial such that
$\s^\a = \pi(P(\g^1,\ldots,\g^k))$ (such a $P$ exists
because of the hypothesis that $\s^\a \in \scal{(\g^i)}$). In
$H^*(\pos)$ we therefore have
$P(\g^1,\ldots,\g^k) = \s^\a + \sum_{m \in M} x_m \s^{\l_m}$ with
$\l_m$ some elements in $I(\pos) - I(\pos_0)$. For each $m$ in $M$,
since $\l_m \not \in I(\pos_0)$, $\l_m$ must contain some element
$\lt_\b$ with $\b \in S_0$ and $\deg(\s^\b) < d'$ and by induction
hypothesis $\ok{\s}{\s^{\l_m}}$. Thus from
$\ok{\s}{P(\g^i)}$ we get $\ok{\s}{\s^\a}$. By recursion
with respect to $\lt_\a$ (Lemma \ref{lemm-plus-grand} point 3)
it follows that $\ok{\s}{\s^\l}$ if
$\l \supset \lt_\a$ and we are done.

We thus have proved that if $\s \in H^*(\pos)$ and
$\tau' \in H_d$ then $\ok{\s}{\tau'}$.
Let finally $\tau \in H^*(\pos)$ such that
$\pi(\tau) \in \scal{(\g^i)}_d$ and $\tau - \pi(\tau) \in H_d$,
and let $\s \in H^*(\pos)$
be arbitrary. Let $P$ as before such that
$P(\g^i) = \tau + \tau'$ with $\tau' \in H_d$.
Since $\ok{\s}{P(\g^i)}$ and we already know that
$\ok{\s}{\tau'}$, we deduce $\ok{\s}{\tau}$.
\end{proo}

We now specialise this lemma.

\begin{lemm}
\label{lemm-sous-algebre}
Let $(\g^i)_{i\in[1,k]}$ be elements in $H^*(\pos_0)$.
Assume that
\begin{itemize}
\item
For all $i$
and all $\s \in H^*(\pos_0)$ we have
$\ok{\s}{\g^i}$.
\item
For each $\a$ in $S_0$,
we have $\s^\a \in \scal{(\g^i)}$.
\end{itemize}
Then for all $\s$ in $H^*(\pos)$ and for all
$\tau$ in $H^*(\pos)$ such that
$\pi(\tau) \in \scal{(\g^i)}$,
we have $\ok{\s}{\tau}$.
\end{lemm}

\begin{lemm}
\label{oreille}
Let $(\g^i)_{i\in[1,k]}$ be elements in $H^*(\pos_0)$
such that:
\begin{itemize}
\item
For all $i \in \{1,\ldots,k\}$ and for all $\s \in H^*(\pos_0)$
we have $\ok{\g^i}{\s}$.
\item
$\g^1,\ldots,\g^k$ generate $H^*(\pos_0)$.
\end{itemize}
Then for all $\s,\tau \in H^*(\pos)$ we have $\ok{\s}{\tau}$.
\end{lemm}

For $\sigma \in H^*(\pos)$ let us denote with $\sigma^{\cdot n}$ resp.
$\sigma^{\pp n}$ the $n$-th power of $\s$ computed with the product
$\cdot$ resp. $\pp$.

\begin{lemm}
\label{lemm-petite-puissance} 
Let $\s,\g^1,\ldots,\g^k \in H^*(\pos_0)$ and $d$ an integer
such that we have:
\begin{itemize}
\item
$\forall \tau \in H^*(\pos_0)\, , \, \ok{\g^i}{\tau}$. 
\item
$\forall n\leq d\, , \, \s^{\cdot n} = \s^{\pp n}$.
\end{itemize}
Then for any polynomial $P(X,X_1,\ldots,X_n)$ of degree at most $d-1$ in
$X$ we have the relation
$$\ok{\s}{P(\s,\g^1,\ldots,\g^k)}.$$
\end{lemm}
In particular $P(\s,\g^1,\ldots,\g^k)$ itself does not depend on the product.
\begin{proo}
In fact we may assume that $P=X^n \, Q(X_1,\ldots,X_k)$ with $n\leq d-1$.
We compute
$$
\begin{array}{lllll}
 \s \cdot P(\s,\g^1,\ldots,\g^k) & = & \s \cdot \s^{\cdot n} \cdot
 Q(\g^1,\ldots,\g^k) 
 & = & \s^{\cdot (n+1)} \cdot Q(\g^1,\ldots,\g^k) \\
 & = & \s^{\pp n+1} \pp Q(\g^1,\ldots,\g^k) & = & \s \pp ( \s^{\pp n} \pp
 Q(\g^1,\ldots,\g^k) \\
& = & \s \pp P(\s,\g^1,\ldots,\g^k). \\
\end{array}
$$
\end{proo}

\begin{lemm}
\label{exterieur}
Let $\lambda,\mu \in I(\pos_0)$, and assume the following:

(\i)
$\forall \nu \in I(\pos_0)$ we have
$c_{\lambda,\mu}^\nu = m_{\lt,\mu}^\nu t_{\lambda,\mu}^\nu$.

(\i\i)
For all $\alpha$ in $S_0$, we have either
\vskip -.4cm
$$
\begin{array}{cl}
\mu \supset \lt_\alpha \mbox{ and }
\sigma^{\lambda} \cdot \sigma^{\alpha}
= \sigma^{\lambda} \odot \sigma^{\alpha} & \mbox{ or }\\

\lambda \supset \lt_\alpha \mbox{ and }
\sigma^{\mu} \cdot \sigma^{\alpha} =
\sigma^{\mu} \odot \sigma^{\alpha}.
\end{array}
$$
\vskip -.2cm

\noindent
Then $\sigma^\lambda \cdot \sigma^\mu = \sigma^\lambda \odot \sigma^\mu$.
\end{lemm}
\begin{proo}
The lemma amounts to the fact that $\forall \nu \in I(\pos)$ we have
$c_{\lambda,\mu}^\nu = m_{\lt,\mu}^\nu t_{\lambda,\mu}^\nu$. This
holds by assumption if $\nu \in I(\pos_0)$. Otherwise there exists a
simple root $\alpha$ in $S_0$ such that $\mu \supset \lt_\alpha$.
By $(\imath\imath)$ we may assume that
$\mu \supset \lt_\alpha \mbox{ and }
\ok{\s^\l}{\s^\alpha}$.
The result follows by the third part of Lemma \ref{lemm-plus-grand}.
\end{proo}

Recall from Proposition \ref{prop-chevalley}, the equality of
combinatorial and cohomogical Chevalley formula.

\begin{nota}
Let $\lambda,\nu \in I(\pos)$ and let $d$ be an integer. We define
\begin{itemize}
\item
$\lambda \cap \pos_0$ the ideal in
$\pos$ defined
by $(\lambda \cap \pos_0)(i,S_1,S_2) = \lambda(i,S_1,S_2) \cap \pos_0(i)$.
\item
$A_{\lambda,d} = \{ \mu \in I(\pos_0) : \deg(\mu) = d , 
\sigma^\lambda \cdot \sigma^\mu \not = \sigma^\lambda \pp \sigma^\mu \}$.
\item
$A_{\lambda,d}^\nu = \{ \mu \in I(\pos_0) : \deg(\mu) = d , 
c_{\lambda,\mu}^\nu \not = t_{\lambda,\mu}^\nu \cdot m_{\lambda,\mu}^\nu \}$.
\end{itemize}
\end{nota}

\begin{lemm}
\label{tous-sauf-1}
Let $\lambda \in I(\pos)$ be a fixed ideal and $d$ be an integer.

(\i) Assume that for all $\mu \in I(\pos)$ such that
$\deg(\mu \cap \pos_0) < d$ we have $\ok{\s^\l}{\s^\mu}$ and
that $\# A_{\lambda,d} \leq 1$. Then $A_{\lambda,d} = \emptyset$.

(\i\i) More specifically, let $\nu \in I(\pos)$ be another ideal and
assume that for all $\mu \in I(\pos)$ such that
$\deg(\mu \cap \pos_0) < d$ we have $\oknu{\l}{\mu}{\nu}$ and that
$\# A_{\lambda,d}^\nu \leq 1$. Then $A_{\lambda,d}^\nu = \emptyset$.
\end{lemm}

\begin{proo}
Let us prove $(\imath)$. Let $\lambda,d$ be as in the lemma.
By Proposition \ref{prop-chevalley} the $d$-th powers of $h$ computed
in $H^*_t(X)$ and $H^*(\pos)$ are equal. Let $\mu \in I(\pos)$
with $\deg(\mu)=d$. We have the following properties:
\begin{itemize}
\item
By Chevalley formula the coefficient of $\sigma^\mu$ in $h^d$
is positive.
\item
If $\sigma \not \in I(\pos_0)$ then $\deg(\mu \cap \pos_0) < d$.
\item
$\ok{h^d}{\s^\l}$.
\end{itemize}
Thus it follows from the hypothesis that $\ok{\s^\l}{\s^\mu}$.
The proof of $(\imath\imath)$ is similar.
\end{proo}

Recall that $X_0$ is the homogeneous space associated
to the marked Dynkin diagram $(D_0,\Lambda)$.

\begin{lemm}
  \label{lefs}
Let $\s$ be a fixed Schubert class and $d$ be an integer.

(\i) Assume that $\dim H^d(X_0) \geq \dim H^{d+1}(X_0)$ and assume that 
$\ok{\s}{\s^\mu}$
for any $\mu \in I(\pos_0)$ such that $\deg(\mu) \leq d$.
Assume moreover that
$X_0$ is finite dimensional.
Then for any $\mu \in I(\pos)$ such that $\deg(\mu \cap \pos_0) \leq d+1$
we have $\ok{\s}{\s^\mu}$.

(\i\i) Assume there exists a subset $C$ of $I(\pos_0)$ such that
for all $\mu \in C$ we have $\ok{\s}{\s^\mu}$. Assume
furthermore that the natural map given by multiplication by $h$:
$${\bigoplus_{\stackrel{\stackrel{\mu\in I(P_0)_d,}{\mu \not \in
        C}}{}}\Z} \cdot \sigma^\mu \to \bigoplus_{\stackrel{\stackrel{\mu'\in 
  I(P_0)_{d+1},}{\mu'\not\in C}}{}}\Z \cdot \sigma^{\mu'}$$
is surjective and that $\ok{\s}{\s^\mu}$ for $\mu \in I(\pos_0)$ such that
$\deg(\mu) \leq d$. Then for any
$\mu$ in $I(\pos)$ such that $\deg(\mu \cap \pos_0) \leq d+1$
we have $\ok{\s}{\s^\mu}$.
\end{lemm}

\begin{proo}
$(\imath)$
Let $h_d : H^d(X_0,\Q) \to H^{d+1}(X_0,\Q)$ and
$\kappa_d : H^d(\pos_0,\Q) \to H^{d+1}(\pos_0,\Q)$
be the maps
induced by multiplication by the class of degree 1.
If $2d \geq \dim(X_0)$, then
by Lefschetz
Theorem (see for example \cite[Theorem 3.1.39]{lazarsfeld}), $h_d$ is
surjective. If $2d < \dim(X_0)$ again by Lefschetz Theorem $h_d$
is injective and hence, under hypothesis $(\imath)$, surjective.
It follows that the induced quotient map $H^d_t(X_0) \to H^{d+1}_t(X_0)$
is also surjective. Since this map identifies with $\kappa_d$, $\kappa_d$
is surjective.

We first prove that for any $\mu \in I(\pos)$ such that
$\deg(\mu \cap \pos_0) \leq d$ we have $\ok{\s}{\s^\mu}$.
In fact, let $\mu$ be such an ideal and let $\a \in S_0$.
If $\mu \supset \lt_\a$ then $\deg(\lt_\a) \leq d$ and thus by assumption
$\ok{\s}{\s^\a}$. By Lemma \ref{lemm-plus-grand}(3) we deduce
that $\ok{\s}{\s^\mu}$. If $\mu \not \supset \lt_\a$ for all elements
$\a \in S_0$ then $\mu \in I(\pos_0)$ and this equality is true by
assumption.

Now we consider $\mu \in I(\pos)$ such that $\deg(\mu)=d+1$. If
$\mu \not \in I(\pos_0)$ then we have already proved the result, so
assume $\mu \in I(\pos_0)$. Since $\kappa_d$ is surjective, there exists
$\rho \in H^*(\pos_0)$ such that $h \cdot \rho = \s^\mu + \tau$,
where $\tau$ is a linear combinaison of some $\s^\mu$ with
$\mu \in I(\pos) - I(\pos_0)$ and therefore $\ok{\s}{\tau}$.
Since $\ok{\s}{(h \cdot \rho)}$, we get $\ok{\s}{\s^\mu}$.

Finally we prove that if $\deg(\mu \cap \pos_0) \leq d+1$ then
$\ok{\s}{\s^\mu}$. This is similar to what we have done in the first case.

$(\imath\imath)$ In this case the proof is as for $(\imath)$.
\end{proo}

We end this subsection with a lemma specific to the finite dimension and even 
specific to the minuscule and cominuscule case. This lemma corresponds to 
Lemma 5.8.(iii) in \cite{TY}. In the following lemma we assume that the 
longest element $w^P$ in $W^P$ is $\Lambda$-(co)minuscule. This is 
equivalent to saying that $\Lambda$ itself is (co)minuscule. We define 
$\pos_0$ 
as the heap of $w^P$.

\begin{lemm}
  \label{poincare}
Let $\lt$ and $\mu$ be two ideals in $\pos_0$ and assume that for all 
ideals $\nu$ in $\pos_0$ except one we have $\oknu{\lt}{\mu}{\nu}$, then 
we have $\oknu{\lt}{\mu}{\nu}$ for all $\nu$.
\end{lemm}

\begin{proo}
We denote with $\s^\lt$, $\s^\mu$ and $\s^\nu$ the classes corresponding to 
the ideals  $\lt$, $\mu$ and $\nu$.

This lemma amounts to the fact that Poincar{\'e} duality is compatible with 
jeu de taquin in this situation. In other words, if $\lt$ is an ideal in 
$\pos_0$, then there exists a unique ideal $\lt^c$ in $\pos_0$ of degree 
$\deg(\pos_0)-\deg(\lt)$ such that for any $\mu$ with 
$\deg(\mu)=\deg(\pos_0)-\deg(\lt)$, we have
$$\s^\lt\cdot\s^\mu=\delta_{\mu,\lt^c}\cdot [{\rm pt}]=\s^\lt\pp\s^\mu$$
where $[{\rm pt}]\in H^{\dim X}(X)$ is the Poincar{\'e} dual of the class 
of a point. This result was proved in \cite[Corollary 4.7]{TY}.

Let us prove the lemma. Let $m=\deg(\pos_0)-(\deg(\lt)+\deg(\mu))$ 
and $h$ the hyperplane class. We have 
$$\s^\lt\cdot\s^\mu=\sum_{\stackrel{\stackrel{\nu\subset
\pos_0}{\deg(\nu)=m}}{}}c_{\lt,\mu}^\nu\s^\nu\ {\rm 
and}\  
\s^\lt\pp\s^\mu=\sum_{\stackrel{\stackrel{\nu\subset
\pos_0}{\deg(\nu)=m}}{}}t_{\lt,\mu}^\nu m_{\lt,\mu}^\nu
\s^\nu.$$
By the discussion above, we have $\ok{\s}{\tau}$ for any classes $\s$ 
and $\tau$ such that $\deg(\s)+\deg(\tau)=\deg(\pos_0)$. Because 
$\deg(h^m\cdot\s^\lt)+\deg(\mu)=\deg(\pos_0)$, we have 
$$(h^m\cdot\s^\lt)\cdot\s^\mu=(h^m\cdot\s^\lt)\pp \s^\mu=
(h^m\pp\s^\lt)\pp\s^\mu.$$
But this is also equal to
$$h^m\cdot(\s^\lt\cdot\s^\mu)= \sum_{\stackrel{\stackrel{\nu\subset
\pos_0}{\deg(\nu)=m}}{}}c_{\lt,\mu}^\nu 
(h^m\cdot\s^\nu)
\ {\rm and}\ 
h^m\pp(\s^\lt\pp\s^\mu)=\sum_{\stackrel{\stackrel{\nu\subset
\pos_0}{\deg(\nu)=m}}{}}t_{\lt,\mu}^\nu m_{\lt,\mu}^\nu
(h^m\pp\s^\nu).$$
As for all $\nu$ of degree $m$ the class 
$h^m\cdot\s^\nu=h^m\pp\s^\nu$ is non zero, and because 
$\oknu{\lt}{\mu}{\nu}$ for all 
$\nu$ but one we get the result.
\end{proo}

\subsection{Strategy for the proof of the main Theorem}
\label{subsection-decomposable}

We now reduce the proof of Theorem \ref{main-theo} to some tractable
cases. 
So let $D$ be a Dynkin
diagram with Weyl group $W$, $\Lambda$ a dominant weight,
$X$ the corresponding Kac-Moody homogeneous
space, and let $u,v,w$ be $\Lambda$-minuscule elements in $W$.
By Corollary \ref{coro-picard-superieur} we may assume that
$\Lambda$ is a fundamental weight; let $d \in D$ be the corresponding
node.

Let us first
introduce some notation. If $(D,d)$ is a marked Dynkin diagram and
$w$ is a $\Lambda_d$-(co)minuscule element, we denote with $P_0(w)$
the slant-irreducible component of $H(w)$ containing the minimal element
$(d,1)$ of $H(w)$. We also denote with $D_0(w) \subset D$ the set of colors
of $P_0(H(w))$.

The heap of $w$ is a slant product of $P_0(w)$ and some
$P_\alpha$'s.
We first prove Theorem \ref{main-theo} in the case $D_0(w)$ is simply laced. 
Arguing by induction, we may assume that Theorem \ref{main-theo}
holds for $P_\alpha$ and for any $u',v',w'$ with
$D_0(w') \varsubsetneq D_0(w)$. Note moreover that $P_0(w)$,
being slant-irreducible, must fall in one of the cases of
\cite{proctor-classification} and its correspoding Dynkin
diagram must correspond to a finite-dimensional
Kac-Moody group (by our assumption). In the following array, we
indicate, depending on $P_0(w)$, which lemma allows to finish
the proof.

$$
\begin{array}{ccc}
P_0(w) \mbox{ as in Proctor's case} & D_0(w) & \mbox{Lemma} \\
1 & A_n & \ref{lemm-an} \\
2 & D_n & \ref{lemm-dnpn} \\
3\ (f=1\, ; \, g=2\, ; \, 2\leq h\leq 4) & E_{h+4} &
\ref{lemm-e6p2} , \ref{lemm-e7p2} , \ref{lemm-e8p2} \\
4\ (f \geq 2\, ;\, h=1) & D_n & \ref{lemm-dnp1} \\
4,5,6,7\ (f = 2\, ;\, 2 \leq h \leq 4) & E_{4+h} &
\ref{lemm-e6p1} , \ref{lemm-e7p1} , \ref{lemm-e8p1} \\
4,5,6,7,8,9,10,11,12,13,14\ (3\leq f \leq 4\, ;\, h=2)
& E_{f+4} & \ref{lemm-e7p7} , \ref{lemm-e8p8}\\
15 & E_7 & \ref{lemm-e7p7}
\end{array}
$$
Once Theorem \ref{main-theo} is proved in case $D_0(w)$ is simply laced, 
we prove it in general thanks to Lemmas \ref{lemm-bnpn}, \ref{lemm-f4p4}, 
\ref{lemm-cnpn} and \ref{lemm-f4p1}.
We end this section with the following notation we shall use in the
sequel:

\begin{nota}
  A generator $\gamma$ of the algebra $H^*(\pos)$ will be called a
  good generator if $\gamma\cdot\s=\gamma\pp\s$ for all classes $\s$
  in $H^*(\pos)$.
\end{nota}
%%% Local Variables: 
%%% mode: latex
%%% TeX-master: t
%%% End: 

\section{Simply laced case}
\label{section-calc}

\subsection{Generators for the cohomology}

For the convenience of the reader, we reproduce here arguments of \cite{CMP} on
well known fact concerning the cohomology of a rational finite dimensional
homogeneous space $G/P$. As we have seen we may assume that
$P$ is maximal. The cohomology with coefficients in a ring $k$ will be denoted
with $H^*(X,k)$.

First, we recall the {\it Borel presentation} of the cohomology ring
with rational coefficients. Let $W$ (resp. $W_P$) be the Weyl group of
$G$ (resp. of $P$). Let $\cP$ denote the weight lattice of $G$. The
Weyl group $W$ acts on $\cP$. We have
$$H^*(G/P,\QQ)\simeq \QQ [\cP]^{W_P}/\QQ [\cP]^W_+, $$
where $\QQ [\cP]^{W_P}$ denotes the ring of 
$W_P$-invariants polynomials on the 
weight lattice, and $\QQ [\cP]^W_+$ is the 
ideal of $\QQ [\cP]^{W_P}$ generated 
by $W$-invariants without constant term 
(see \cite[Proposition 27.3]{borel} or \cite[Theorem 5.5]{BGG}).

Recall that the full invariant algebra $\QQ [\cP]^W$ is a polynomial
algebra $\QQ[F_{e_1+1},\ldots ,F_{e_{max}+1}]$, where $e_1,\ldots
,e_{max}$ is the set $E(G)$ of exponents of $G$.  
If $d_1,\ldots ,d_{max}$ denote the
exponents of a Levi subgroup $L(P)$of $P$, we get that $\QQ
[\cP]^{W_P}=\QQ[I_1,I_{d_1+1},\ldots ,I_{d_{max}+1}]$, where $I_1$
represents the fundamental weight $\varpi_P$ defining
$P$. Geometrically, it corresponds to the 
hyperplane class. 

Each $W$-invariant $F_{e_i+1}$ must be
interpreted as a polynomial relation between the $W_P$-invariants
$I_1,I_{d_1+1},\ldots ,I_{d_{max}+1}$. In particular, if $e_i$ is also 
an exponent of the semi-simple $L(P)$ part of $P$, 
this relation allows to eliminate $I_{e_i+1}$. 
We thus get the presentation, by generators and relations,
$$H^*(G/P,{\QQ})\simeq \QQ[I_1,I_{p_1+1},\ldots ,I_{p_n+1}]/
(R_{q_1+1},\ldots ,R_{q_r+1}),$$
where $\{p_1,\ldots ,p_n\}=E(L(P))-E(G)$ and $\{q_1,\ldots
,q_r\}=E(G)-E(L(P))$. 

\subsection{Quadrics}

Let us start with the case of quadrics. Thus we consider  the system
of $\varpi_1$-minuscule $D_n$-colored posets $\pos_0$ given by the
following maximal element:

\centerline{\begin{pspicture*}(0,0)(4.2000003,6.0)

\psellipse(0.6,5.4)(0.096,0.096)
\psline(0.532,5.332)(0.668,5.468)
\psline(0.532,5.468)(0.668,5.332)
\psellipse[fillstyle=solid,fillcolor=black](0.6,0.6)(0.096,0.096)
\psellipse(1.2,4.8)(0.096,0.096)
\psline(1.132,4.732)(1.268,4.868)
\psline(1.132,4.868)(1.268,4.732)
\psellipse[fillstyle=solid,fillcolor=black](1.2,1.2)(0.096,0.096)
\psellipse(1.8,4.2)(0.096,0.096)
\psline(1.732,4.132)(1.868,4.268)
\psline(1.732,4.268)(1.868,4.132)
\psellipse[fillstyle=solid,fillcolor=black](1.8,1.8)(0.096,0.096)
\psellipse(2.4,3.6)(0.096,0.096)
\psline(2.332,3.532)(2.468,3.668)
\psline(2.332,3.668)(2.468,3.532)
\pspolygon[fillstyle=solid,fillcolor=black](2.52,3.0)(2.4,3.12)(2.28,3.0)(2.4,2.88)
\psellipse[fillstyle=solid,fillcolor=black](2.4,2.4)(0.096,0.096)
\pspolygon[fillstyle=solid,fillcolor=black](3.12,3.0)(3.0,3.12)(2.88,3.0)(3.0,2.88)
\psline(0.668,5.332)(1.132,4.868)
\psline(1.868,4.132)(2.332,3.668)
\psline(2.468,3.532)(2.94,3.06)
\psline(2.4,3.504)(2.4,3.12)
\psline(2.4,2.88)(2.4,2.496)
\psline(2.94,2.94)(2.468,2.468)
\psline(2.332,2.332)(1.868,1.868)
\psline(1.132,1.132)(0.668,0.668)
\psline[linestyle=dashed](1.268,4.732)(1.732,4.268)
\psline[linestyle=dashed](1.732,1.732)(1.268,1.268)
\uput[l](2.4,3.0){$\s^{n-1}$}
\uput[r](3.0,3.0){$\tau^{n-1}$}

\end{pspicture*}}

We have $S_0=\{n-1,n\}$. For $i\in \{n-1,n\}$, let $(D_i,d_i)$ be a
marked Dynkin diagram and $P_i$ be any  $d_i$-minuscule $D_i$-colored
poset. Set $\pos=\pos_{\pos_0,\{P_{n-1},P_n\}}$. 

\begin{lemm}
\label{lemm-dnp1}
With the above notation, assume that Conjecture \ref{main_conj} holds
for $P_{n-1}$, $P_n$ and any $\lt$ in $I(\pos)$ with
$D_0(\lt)\varsubsetneq D_n$. Then Conjecture \ref{main_conj} holds for
$\pos$. 
\end{lemm}

\begin{proo}
Let us define the degree $n-1$ ideals $\lt_{n-1}=\scal{(\a_{n-1},1)}$
and $\mu_{n-1}=\scal{(\a_n,1)}$ in $\pos_0$. The corresponding Schubert
classes are denoted with $\s^{n-1}$ and with $\tau^{n-1}$. Let
$\{\gamma^1,\gamma^{n-1}\}$ be a set of generators of the cohomology ring of
the quadric, with $\deg(\gamma^i)=i$. The variety $D_n/P_1$ has dimension
$2(n-1)$, the dimensions of $H^d(D_n/P_1)$ are

\vskip 0.2 cm

\centerline{\begin{tabular}{c|cc}
\hline
$d$&$d\neq n-1$&$n-1$\\
\hline
$\dim H^d(D_n/P_1)$&1&2\\
\hline
\end{tabular}}

\vskip 0.2 cm

Since by assumption the conjecture holds for any $\lt\in I(\pos)$ with
$D(\lt)\varsubsetneq D_n$, we have $c_{\lt,\mu}^\nu=t_{\lt,\mu}^\nu$
as soon as $\nu\not\supset\pos_0$.

By Proposition \ref{prop-chevalley}, $\gamma^1$ is a good
generator. For $\gamma^{n-1}$, by Lemma \ref{lefs}, we have
the equality $\gamma^{n-1}\cdot\s^\lt=\gamma^{n-1}\pp\s^\lt$ for any
class $\s^\lt$  with $\deg(\lt\cap\pos_0)\leq n-2$. Furthermore, we
have the equality 
$c_{\gamma^{n-1},\s^\lt}^{\s^\nu}=t_{\gamma^{n-1},\s^\lt}^{s^\nu}$ for
$\nu\not\supset \pos_0$. For $\deg(\s^\lt)=n-1$, we are left with the
equality $c_{\gamma^{n-1},\s^\lt}^{\s^\nu}=t_{\gamma^{n-1},\s^\lt}^{s^\nu}$ for
$\nu= \pos_0$. But in this case we are reduced to the same computation in
the quadric and the result follows, for example by Poincar{\'e} duality.
For higher degree, we use Lemma \ref{lefs}.
\end{proo}

\subsection{Type $A_n$}

In this case, we consider  the system
of $\varpi_p$-minuscule $A_n$-colored posets $\pos_0$ given by the
poset of a Grassmannian $\G(p,n+1)$:

\centerline{\begin{pspicture*}(0,0)(8.5,7.8)

\pspolygon[fillstyle=solid,fillcolor=black](1.32,3.6)(1.2,3.72)(1.08,3.6)(1.2,3.48)
\psellipse(1.8,4.2)(0.096,0.096)
\psline(1.732,4.132)(1.868,4.268)
\psline(1.732,4.268)(1.868,4.132)
\psellipse[fillstyle=solid,fillcolor=black](1.8,3.0)(0.096,0.096)
\psellipse(2.4,4.8)(0.096,0.096)
\psellipse(2.4,3.6)(0.096,0.096)
\psline(2.332,3.532)(2.468,3.668)
\psline(2.332,3.668)(2.468,3.532)
\psellipse(3.0,4.2)(0.096,0.096)
\psellipse[fillstyle=solid,fillcolor=black](3.0,1.8)(0.096,0.096)
\psellipse(3.6,6.0)(0.096,0.096)
\psellipse(3.6,2.4)(0.096,0.096)
\psline(3.532,2.332)(3.668,2.468)
\psline(3.532,2.468)(3.668,2.332)
\psellipse[fillstyle=solid,fillcolor=black](3.6,1.2)(0.096,0.096)
\psellipse(4.2,6.6)(0.096,0.096)
\psellipse(4.2,5.4)(0.096,0.096)
\psellipse(4.2,3.0)(0.096,0.096)
\psline(4.132,2.932)(4.268,3.068)
\psline(4.132,3.068)(4.268,2.932)
\psellipse[fillstyle=solid,fillcolor=black](4.2,1.8)(0.096,0.096)
\psellipse(4.8,7.2)(0.096,0.096)
\psellipse(4.8,6.0)(0.096,0.096)
\psellipse[fillstyle=solid,fillcolor=black](4.8,2.4)(0.096,0.096)
\psellipse(5.4,6.6)(0.096,0.096)
\psellipse(5.4,4.2)(0.096,0.096)
\psline(5.332,4.132)(5.468,4.268)
\psline(5.332,4.268)(5.468,4.132)
\psellipse(6.0,4.8)(0.096,0.096)
\psline(5.932,4.732)(6.068,4.868)
\psline(5.932,4.868)(6.068,4.732)
\psellipse[fillstyle=solid,fillcolor=black](6.0,3.6)(0.096,0.096)
\psellipse(6.6,5.4)(0.096,0.096)
\psline(6.532,5.332)(6.668,5.468)
\psline(6.532,5.468)(6.668,5.332)
\psellipse[fillstyle=solid,fillcolor=black](6.6,4.2)(0.096,0.096)
\pspolygon[fillstyle=solid,fillcolor=black](7.32,4.8)(7.2,4.92)(7.08,4.8)(7.2,4.68)
\psline(3.668,1.268)(4.132,1.732)
\psline(4.268,1.868)(4.732,2.332)
\psline(4.732,2.468)(4.268,2.932)
\psline(4.132,1.868)(3.668,2.332)
\psline(3.532,1.268)(3.068,1.732)
\psline(3.068,1.868)(3.532,2.332)
\psline(3.668,2.468)(4.132,2.932)
\psline(6.068,3.668)(6.532,4.132)
\psline(6.668,4.268)(7.14,4.74)
\psline(5.932,3.668)(5.468,4.132)
\psline(6.532,4.268)(6.068,4.732)
\psline(7.14,4.86)(6.668,5.332)
\psline(5.468,4.268)(5.932,4.732)
\psline(6.068,4.868)(6.532,5.332)
\psline(1.868,3.068)(2.332,3.532)
\psline(2.468,3.668)(2.932,4.132)
\psline(1.732,3.068)(1.26,3.54)
\psline(2.332,3.668)(1.868,4.132)
\psline(2.932,4.268)(2.468,4.732)
\psline(1.26,3.66)(1.732,4.132)
\psline(1.868,4.268)(2.332,4.732)
\psline(4.268,5.468)(4.732,5.932)
\psline(4.868,6.068)(5.332,6.532)
\psline(4.104,5.4)(4.104,5.4)
\psline(4.132,5.468)(3.668,5.932)
\psline(4.732,6.068)(4.268,6.532)
\psline(5.332,6.668)(4.868,7.132)
\psline(3.668,6.068)(4.132,6.532)
\psline(4.268,6.668)(4.732,7.132)
\psline[linestyle=dashed](6.532,5.468)(5.468,6.532)
\psline[linestyle=dashed](5.932,4.868)(4.868,5.932)
\psline[linestyle=dashed](5.332,4.268)(4.268,5.332)
\psline[linestyle=dashed](4.868,2.468)(5.932,3.532)
\psline[linestyle=dashed](4.268,3.068)(5.332,4.132)
\psline[linestyle=dashed](4.132,3.068)(3.068,4.132)
\psline[linestyle=dashed](3.068,4.268)(4.132,5.332)
\psline[linestyle=dashed](2.468,4.868)(3.532,5.932)
\psline[linestyle=dashed](3.532,2.468)(2.468,3.532)
\psline[linestyle=dashed](2.932,1.868)(1.868,2.932)
\uput[r](7.2,4.8){$\tau^{n-p+1}$}
\uput[r](4.8,2.4){$\tau^{3}$}
\uput[r](4.2,1.8){$\tau^{2}$}
\uput[l](3.0,1.8){$\s^{2}$}
\uput[l](1.2,3.6){$\s^{p}$}
\uput[d](3.6,1.2){$\s^{1}$}

\end{pspicture*}}

We have $S_0=\{1,n\}$. For $i\in \{1,n\}$, let $(D_i,d_i)$ be a
marked Dynkin diagram and $P_i$ be any  $d_i$-minuscule $D_i$-colored
poset. Set $\pos=\pos_{\pos_0,\{P_{1},P_n\}}$. 

\begin{lemm}
\label{lemm-an}
With the above notation, assume that Conjecture \ref{main_conj} holds
for $P_{1}$, $P_n$, and any $\lt$ in $I(\pos)$ with
$D_0(\lt)\varsubsetneq A_n$. Then Conjecture \ref{main_conj} holds for
$\pos$. 
\end{lemm}

\begin{proo}
By Proposition \ref{prop-chevalley}, we may assume $n\geq2$, by
irreducibility of $\pos_0$ we may then assume $n\geq3$ and by Proposition
\ref{lemm-dnpn}, we may assume that $n\geq4$.

Let us define the degree $i$ ideals $\lt_{i}=\scal{(\a_{p+1-i},1)}$
for $i\in[1,p]$ and $\mu_{i}=\scal{(\a_{p+i-1},1)}$ for $i\in[1,n+1-p]$
  in $\pos_0$. The corresponding Schubert cells are denoted with
  $\s^{i}$ and with $\tau^{i}$. Take
  $(\gamma^i)_{i\in[1,p]}$ 
a set of generators of the cohomology ring of the Grassmannian, with
$\deg(\gamma^i)=i$.

Since by assumption the conjecture holds for any $\lt\in I(\pos)$ with
$D_0(\lt)\varsubsetneq A_n$, we have $c_{\lt,\mu}^\nu=t_{\lt,\mu}^\nu$
as soon as $\deg(\nu\cap\pos_0)\leq2n-3$.
In particular because for $n\geq4$ we have $i+j\leq 2n-3$ for $i\leq
p$ and $j\leq n+1-p$, the equality $\g^i\cdot\s^j=\g^i\pp\s^j$ holds
for all $i\leq p$ and $j\leq p$ and the equality
$\g^i\cdot\tau^j=\g^i\pp\tau^j$ holds for all $i\leq p$ and $j\leq
n+1-p$.

Now let $\lt\in I(\pos_0)$. If $\lt\supset\lt_p$ or
$\lt\supset\mu_{n+1-p}$, then by recursion with respect to $\lt_p$ or
$\mu_{n+1-p}$ we have $\g^i\cdot\s^\lt=\g^i\pp\s^\lt$.

If $\lt\not\supset\lt_p$ and $\lt\not\supset\mu_{n+1-p}$, then we
first consider the case where $\lt$ is an ideal of the form
$\scal{(\a_k,l)}$ for some simple root $\a_k$ and some integer $l$. We
prove the equality $\g^i\cdot\s^\lt=\g^i\pp\s^\lt$ by induction on
$\deg(\lt)$ in that case. We may of course assume that $\lt$ is
distinct from all the $\lt_i$ and the $\mu_j$. Consider the two
subideals $\lt'$ and $\lt''$ in $\lt$ described by
$\lt'=\scal{(\a_{k-1},l')}$ and $\lt''=\scal{(\a_{k+1},l'')}$ where
$l'=\max\{a\ /\ (\a_{k-1},a)\in \lt\}$ and $l''=\max\{a\ /\
(\a_{k+1},a)\in \lt\}$. By recursion with respect to $\lt'$ or
$\lt''$, we have $c_{\g^i,\s^\lt}^{\s^\nu}=t_{\g^i,\s^\lt}^{\s^\nu}$
for any $\nu$ not containing $(\a_{k-1},l'+1)$ or
$(\a_{k+1},l''+1)$. By induction 
on $\pos_0$ it is also true if $\nu$ does not contain $(\a_1,1)$ or
$(\a_n,1)$. For an ideal $\nu$ in $\pos$ containing all these elements
of $\pos_0$, we have $\deg(\nu\cap \pos_0)\geq \deg(\lt)+n-1$. For such a
$\nu$ we have $c_{\g^i,\s^\lt}^{\s^\nu}=0=t_{\g^i,\s^\lt}^{\s^\nu}$
for degree reasons. 

We finish by dealing with $\lt\in I(\pos_0)$ not of the previous
form. Let us consider the set $M(\lt)$ of maximal elements in
$\lt$. For $(\a_k,l)\in M(\lt)$, define the ideal
$\lt(\a_k,l)=\scal{(\a_k,l)}$. By what we have just done, we have
$\g^i\cdot \s^{\lt(\a_k,l)}=\g^i\pp\s^{\lt(\a_k,l)}$. In particular we
can use recursion with respect to $\lt(\a_k,l)$ and we deduce that
$c_{\g^i,\s^\lt}^{\s^\nu}=t_{\g^i,\s^\lt}^{\s^\nu}$ for any $\nu$ not
containing $(\a_{k},l+1)$. By induction on $\pos_0$ it is also true if
$\nu$ does not contain $(\a_1,1)$ or $(\a_n,1)$. For an ideal $\nu$ in
$\pos$ containing all the elements $(\a_k,l+1)$ for $(\a_k,l)\in
M(\lt)$ as well as $(\a_1,1)$ and $(\a_n,1)$, we have $\deg(\nu\cap
\pos_0)\geq \deg(\lt)+n$. For such a $\nu$ we have
$c_{\g^i,\s^\lt}^{\s^\nu}=0=t_{\g^i,\s^\lt}^{\s^\nu}$ for degree reasons.
\end{proo}

\subsection{Type $D_n$}

In this case, we consider  the system
of $\varpi_{n-1}$-minuscule $D_n$-colored posets $\pos_0$ given by the
posets of an orthogonal Grassmannian $\G_Q(n,2n)$.
We have $S_0=\{1,n\}$. For $i\in \{1,n\}$, let $(D_i,d_i)$ be a
marked Dynkin diagram and $P_i$ be any  $d_i$-minuscule $D_i$-colored
poset. Set $\pos=\pos_{\pos_0,\{P_{1},P_n\}}$. 
The quiver $\pos$ for $D_7$ was described in (\ref{equa-systeme-d7}).

\begin{lemm}
\label{lemm-dnpn}
With the above notation, assume that Conjecture \ref{main_conj} holds
for $P_{1}$, $P_n$, and any $\lt$ in $I(\pos)$ with
$D_0(\lt)\varsubsetneq D_n$. Then Conjecture \ref{main_conj} holds for
$\pos$. 
\end{lemm}

\begin{proo}
By irreducibility of $\pos_0$ we may assume $n\geq4$ and by Proposition
\ref{lemm-dnp1}, we may assume that $n\geq5$.

Let us define the degree $i$ ideals $\lt_{i}=\scal{(\a_{n-i},1)}$
for $i\in[1,n-1]$ and the degree 3 ideal $\mu_3=\scal{(\a_{n},1)}$. The
corresponding Schubert cells are denoted with $\s^{i}$ and with
$\tau^{3}$. Take $(\g^i)_{i\in[1,n-1]}$
a set of generators of the cohomology ring of the isotropic Grassmannian, with
$\deg(\gamma^i)=i$.

Since by assumption the conjecture holds for any $\lt\in I(\pos)$ with
$D_0(\lt)\varsubsetneq D_n$, we have $c_{\lt,\mu}^\nu=t_{\lt,\mu}^\nu$
as soon as $\deg(\nu\cap\pos_0)\leq2n-3$.
In particular because for $n\geq5$ we have $i+3\leq 2n-3$ for $i\leq
n-1$, the equality $\g^i\cdot\tau^3=\g^i\pp\tau^3$ holds
for all $i\leq n-1$.

For any ideal $\lt$ in $\pos$
containing $\mu_3$, we obtain by recursion with respect to $\tau^3$
that $c_{\g^i,\s^\lt}^{\s^\nu}=t_{\g^i,\s^\lt}^{\s^\nu}$ for $\nu$ with
$\deg(\nu\cap\pos_0)\leq 2n-1$. In particular if $\deg(\lt)\leq n-1$
we have $\deg(\lt)+i\leq 2n-1$ and
$\g^i\cdot\s^\lt=\g^i\pp\s^\lt$. But there is a unique class in
$H^*(\pos)$ of degree $j\in[1,n-1]$ not bigger than
$\tau^3$: the class $\s^j$, thus by Lemma \ref{tous-sauf-1} we obtain
$\g^i\cdot\s^j=\g^i\pp\s^j$ for all $i$ and $j$ in $[1,n-1]$. 

If $\lt\supset\lt_{n-1}$, then by recursion with respect to $\lt_{n-1}$
we have $\g^i\cdot\s^\lt=\g^i\pp\s^\lt$ for $i\in[1,n-1]$. 

If $\lt\not\supset\lt_{n-1}$, then we first consider the case where
$\lt$ is an ideal of the form $\scal{(\a_k,l)}$ for some simple root
$\a_k$ and some integer $l$. We prove the equality
$\g^i\cdot\s^\lt=\g^i\pp\s^\lt$ by induction on $\deg(\lt)$ 
in that case. We may of course assume that $\lt$ is distinct from all
the $\lt_i$ and from $\mu_3$. We have to discuss two cases. If
$k\not\in\{n-2,n-1,n\}$, then consider the three subideals $\lt'$,
$\lt''$ and $\lt'''$ in $\lt$ described by $\scal{(\a_{k-1},l')}$,
$\scal{(\a_{k+1},l'')}$ and $\scal{(\a_{k'},l''')}$
where $l'=\max\{a\ /\ (\a_{k-1},a)\in \lt\}$, $l''=\max\{a\ /\
(\a_{k+1},a)\in \lt\}$ and $(\a_{k'},l''')$ is the largest
element in $\lt$ with $k'\in\{n-1,n\}$. If $k\in\{n-2,n-1,n\}$, then
consider the subideal $\lt'$ in $\lt$ described by
$\scal{(\a_{k'},l')}$ where 
$(\a_{k'},l')$ is the largest element in $\lt$ with $\{k,k'\}=\{n-1,n\}$. 
By recursion with respect to
$\lt'$, $\lt''$ or $\lt'''$, we have
$c_{\g^i,\s^\lt}^{\s^\nu}=t_{\g^i,\s^\lt}^{\s^\nu}$
for any $\nu$ not containing
$(\a_{k-1},l'+1)$, $(\a_{k+1},l''+1)$ and $(\a_{k'},l'''+1)$ in the
first case and 
$(\a_{k'},l'+1)$ in the second one. By induction
on $\pos_0$ it is also true if $\nu$ does not contain $(\a_1,1)$. For
an ideal $\nu$ in $\pos$ containing all these elements of $\pos_0$, we have
$\deg(\nu\cap \pos_0)\geq \deg(\lt)+n-1$. For such a $\nu$ we have
$c_{\g^i,\s^\lt}^{\s^\nu}=0=t_{\g^i,\s^\lt}^{\s^\nu}$
for degree reasons. Remark that here this
method does not work for $i=n-1$ however, we proved the equality
$c_{\g^{n-1},\s^\lt}^{\s^\nu}=t_{\g^{n-1},\s^\lt}^{\s^\nu}$
for all $\lt$ and $\nu$ except for some $\nu=\scal{(\a_k,l+1)}$.
We obtain $c_{\g^{n-1},\s^\lt}^{\s^\nu}=t_{\g^{n-1},\s^\lt}^{\s^\nu}$
for this $\nu$ by Lemma \ref{poincare}.

We finish by dealing with $\lt$ not of the previous form. Let us
consider the set $M(\lt)$ of maximal elements in $\lt$. For $(\a_k,l)\in
M(\lt)$, define the ideal $\lt(\a_k,l)=\scal{(\a_k,l)}$. We have
$\g^i\cdot \s^{\lt(\a_k,l)}=\g^i\pp\s^{\lt(\a_k,l)}$. In particular we can
use recursion with respect to $\lt(\a_k,l)$ and we deduce that
$c_{\g^i,\s^\lt}^{\s^\nu}=t_{\g^i,\s^\lt}^{\s^\nu}$
for any $\nu$ not containing
$(\a_{k},l+1)$. By induction on $\pos_0$ it is also true if $\nu$ does
not contain $(\a_1,1)$. For an ideal $\nu$ in $\pos$ containing all
the elements $(\a_k,l+1)$ for $(\a_k,l)\in M(\lt)$ as well as
$(\a_1,1)$, we have $\deg(\nu\cap \pos_0)\geq \deg(\lt)+n-1$. For such
a $\nu$ we have
$c_{\g^i,\s^\lt}^{\s^\nu}=t_{\g^i,\s^\lt}^{\s^\nu}$
for degree reasons. Once more, for $i=n-1$, we proved the equality 
$c_{\g^{n-1},\s^\lt}^{\s^\nu}=t_{\g^{n-1},\s^\lt}^{\s^\nu}$
for all $\nu$ except for $\nu=\scal{(\a_1,1),(\a_k,l)\in M(\lt)}$. We
conclude by Lemma \ref{poincare}.
\end{proo}

\subsection{Type $E_6$ case}
\label{section-e_6}

Let us start with the case of $E_6/P_1$.
Thus we consider the system of $\varpi_1$-minuscule
$E_6$-colored posets $\pos_0$ given by the
following maximal element:

\centerline{%auto-ignore

\psset{unit=6mm}

\begin{pspicture*}(1.5,0.5)(-6.5,-10.5)

\psline(-3.88,-9.88)(-3.12,-9.12)
\psline(-2.88,-8.88)(-2.12,-8.12)
\psline(-1.88,-7.88)(-1.12,-7.12)
\psline(-2,-7.84)(-2,-7.16)
\psline(-0.88,-6.88)(-0.12,-6.12)
\psline(-1.12,-6.88)(-1.88,-6.12)
\psline(-2,-6.84)(-2,-6.16)
\psline(-0.12,-5.88)(-0.88,-5.12)
\psline(-1.88,-5.88)(-1.12,-5.12)
\psline(-2.12,-5.88)(-2.88,-5.12)
\psline(-1.12,-4.88)(-1.88,-4.12)
\psline(-2.88,-4.88)(-2.12,-4.12)
\psline(-3.12,-4.88)(-3.88,-4.12)
\psline(-2,-3.84)(-2,-3.16)
\psline(-2.12,-3.88)(-2.88,-3.12)
\psline(-3.88,-3.88)(-3.12,-3.12)
\psline(-2,-2.84)(-2,-2.16)
\psline(-2.88,-2.88)(-2.12,-2.12)
\psline(-1.88,-1.88)(-1.12,-1.12)
\psline(-0.88,-0.88)(-0.12,-0.12)

\psdiamond[fillstyle=solid,fillcolor=black](-0,-6)(0.256,0.256)
\put(0.374,-6.125){$\tau^5$} 
\pscircle(-0,-0){0.16}
\pscircle*(-1,-7){0.16}
\pscircle(-1,-5){0.16}  \psline(-0.88,-4.88)(-1.12,-5.12)
\psline(-0.88,-5.12)(-1.12,-4.88) 
\pscircle(-1,-1){0.16}
\pscircle*(-2,-8){0.16}
\psdiamond[fillstyle=solid,fillcolor=black](-2,-7)(0.256,0.256)
\put(-2.958,-7.125){$\sigma^4$} 
\pscircle(-2,-6){0.16}  \psline(-1.88,-5.88)(-2.12,-6.12)
\psline(-1.88,-6.12)(-2.12,-5.88) 
\pscircle(-2,-4){0.16}
\pscircle(-2,-3){0.16}
\pscircle(-2,-2){0.16}
\pscircle*(-3,-9){0.16}
\pscircle(-3,-5){0.16}  \psline(-2.88,-4.88)(-3.12,-5.12)
\psline(-2.88,-5.12)(-3.12,-4.88) 
\pscircle(-3,-3){0.16}
\pscircle*(-4,-10){0.16}
\pscircle(-4,-4){0.16}  \psline(-3.88,-3.88)(-4.12,-4.12)
\psline(-3.88,-4.12)(-4.12,-3.88)   \put(-4.958,-4.125){$\sigma^8$} 

\end{pspicture*}}

We have $S_0 = \{2,6\}$. For $i \in \{2,6\}$
let $(D_i,d_i)$ be a marked Dynkin diagram and $P_i$ be any $d_i$-minuscule
$D_i$-colored poset. Set $\pos = \pos_{\pos_0,(P_2,P_6)}$
with notation \ref{nota-system}.

\begin{lemm}
\label{lemm-e6p1}
With the above notation, assume that Conjecture \ref{main_conj}
holds for $P_2$,$P_6$, and any $\lambda$ in $I(\pos)$ with
$D_0(\lambda) \varsubsetneq E_6$. Then Conjecture \ref{main_conj}
holds for $\pos$.
\end{lemm}

\begin{proo}
We consider the ideals
$\l_4 = \scal{(\alpha_2,1)}$ resp. $\mu_5=\scal{(\a_6,1)}$ in $\pos_0$, 
of degree 4 resp. 5. The corresponding Schubert cells are denoted with
$\s^4$ resp. $\tau^5$. Let $\{ \g^1,\g^4 \}$ be a set of generators of the
cohomology ring of $E_6/P_1$,
with $\deg(\g^i)=i$. The variety $E_6/P_1$ has dimension 16 and the
dimensions of $H^d(E_6/P_1)$ are

\vskip 0.2 cm

\centerline{\begin{tabular}{c|llllllllllllllllllllllll}
\hline
$d$&0&1&2&3&4&5&6&7&8\\
\hline
$\dim H^d(E_6/P_1)$&1&1&1&1&2&2&2&2&3\\
\hline
\end{tabular}}

\vskip 0.2 cm

Since by assumption the conjecture holds for any
$\lambda \in I(\pos)$ with $D_0(\lambda) \varsubsetneq E_6$,
we have $\oknu{\l}{\mu}{\nu}$ as soon as
$\deg(\nu \cap \pos_0) \leq 9$.

By Proposition \ref{prop-chevalley}, $\g^1$ is a good generator.
By the above argument, we have $\gamma^4\cdot\s=\gamma^4\pp\s$ for all 
$\s$ of degree at most 5. Furthermore, for any ideal $\lt$ in $I(\pos)$ 
such that $\deg(\lt\cap \pos_0)=5$, we have $\lt\supset \lt_4$, $\lt\supset 
\mu_5$ or $\lt\subset \pos_0$. In any case we have 
$\gamma^4\cdot\s^\lt=\gamma^4\pp\s^\lt$
either by recursion with respect to $\s^4$, to $\tau^5$ or by the previous 
argument. By Proposition \ref{lefs} we get the same equality for $\s^\lt$ 
with $\deg(\lt\cap\pos_0)\leq7$.

Let $\s^\lt$ be a degree 8 class associated to an ideal $\lt$ in $\pos$. If
$\lt$ is not contained in $\pos_0$, then $\deg(\lt\cap \pos_0)\leq7$ 
and 
we have $\gamma^4\cdot\s^\lt=\gamma^4\pp\s^\lt$. Moreover, if
$\lt\supset \mu_5$, then by recursion with respect to $\tau^5$ we have
$\gamma^4\cdot\s^\lt=\gamma^4\pp\s^\lt$. Finally, there is a unique ideal 
$\lt$ in
$\pos$ satisfying $\lt\subset \pos_0$ and $\lt\not\supset\mu_5$. For this
class we conclude by Lemma \ref{tous-sauf-1}.

Let $\s^\lt$ be a class associated to an ideal $\lt$ in $\pos$ such that
$\deg(\lt\cap \pos_0)=8$. If $\lt\not\subset \pos_0$, then $\lt\supset\lt_4$
or $\lt\supset\mu_5$ and we have $\gamma^4\cdot\s^\lt=\gamma^4\pp\s^\lt$ by 
recursion
with respect to $\s^4$ or $\tau^5$. If $\lt\subset \pos_0$, then we
already proved the equality $\gamma^4\cdot\s^\lt=\gamma^4\pp\s^\lt$. By 
Lemma \ref{lefs}, we get equality $\gamma^4\cdot\s^\lt=\gamma^4\pp\s^\lt$ 
for higher degree classes.
\end{proo}

We now consider the case of $E_6/P_2$. Thus we consider the system of 
$\varpi_2$-minuscule $E_6$-colored posets $\pos_0$ given by:

\centerline{%auto-ignore

\psset{unit=6mm}

\begin{pspicture*}(1.5,0.5)(-5.5,-4.5)

\psline(-2,-3.84)(-2,-3.16)
\psline(-1.88,-2.88)(-1.12,-2.12)
\psline(-2.12,-2.88)(-2.88,-2.12)
\psline(-0.88,-1.88)(-0.12,-1.12)
\psline(-1.12,-1.88)(-1.88,-1.12)
\psline(-2.88,-1.88)(-2.12,-1.12)
\psline(-3.12,-1.88)(-3.88,-1.12)
\psline(-0.12,-0.88)(-0.88,-0.12)
\psline(-1.88,-0.88)(-1.12,-0.12)
\psline(-2,-0.84)(-2,-0.16)
\psline(-2.12,-0.88)(-2.88,-0.12)
\psline(-3.88,-0.88)(-3.12,-0.12)

\psdiamond[fillstyle=solid,fillcolor=black](-0,-1)(0.256,0.256)   \put(0.374,-1.125){$\tau^4$}
\pscircle*(-1,-2){0.16}
\pscircle(-1,-0){0.16}  \psline(-0.88,0.12)(-1.12,-0.12)  \psline(-0.88,-0.12)(-1.12,0.12)
\pscircle*(-2,-4){0.16}
\pscircle*(-2,-3){0.16}
\pscircle(-2,-1){0.16}  \psline(-1.88,-0.88)(-2.12,-1.12)  \psline(-1.88,-1.12)(-2.12,-0.88)
\pscircle(-2,-0){0.16}  \psline(-1.88,0.12)(-2.12,-0.12)  \psline(-1.88,-0.12)(-2.12,0.12)
\pscircle*(-3,-2){0.16}
\pscircle(-3,-0){0.16}  \psline(-2.88,0.12)(-3.12,-0.12)  \psline(-2.88,-0.12)(-3.12,0.12)
\psdiamond[fillstyle=solid,fillcolor=black](-4,-1)(0.256,0.256)   \put(-4.958,-1.125){$\sigma^4$}

\end{pspicture*}}

We have $S_0 = \{1,6\}$. For $i \in \{1,6\}$
let $(D_i,d_i)$ be a marked Dynkin diagram and $P_i$ be any $d_i$-minuscule
$D_i$-colored poset. Set $\pos = \pos_{\pos_0,(P_1,P_6)}$
with notation \ref{nota-system}.

\begin{lemm}
\label{lemm-e6p2}
With the above notation, assume that Conjecture \ref{main_conj}
holds for $P_1$,$P_6$, and any $\lambda$ in $I(\pos)$ with
$D_0(\lambda) \varsubsetneq E_6$. Then Conjecture \ref{main_conj}
holds for $\pos$.
\end{lemm}

\begin{proo}
We consider the ideals $\l_4 = \scal{(\alpha_1,1)}$ resp. 
$\mu_4=\scal{(\a_6,1)}$ in $\pos_0$ both are of degree 4. The corresponding 
Schubert cells are denoted with $\s^4$ resp. $\tau^4$. Let $\{
\g^1,\g^3,\g^4 \}$  
be a set of generators of the cohomology ring of $E_6/P_2$, with
$\deg(\g^i)=i$.

Since by assumption the conjecture holds for any
$\lambda \in I(\pos)$ with $D_0(\lambda) \varsubsetneq E_6$,
we have $\oknu{\l}{\mu}{\nu}$ as soon as
$\deg(\nu \cap \pos_0) \leq 9$.

By Proposition \ref{prop-chevalley}, $\g^1$ is a good generator.
By the previous argument, we have $\g^3\cdot\s=\g^3\pp\s$ for all $\s$ in
$\pos_0$ of degree at most 6. In particular this holds for $\s=\s^4$ or 
$\s=\tau^4$. Let $\lt\subset \pos_0$ with
$\deg(\lt)\geq7$. We have $\lt\supset\lt_4$ or $\lt\supset\mu_4$ and
by recursion with respect to $\s^4$ or $\tau^4$ we get the equality
$\g^3\cdot\s=\g^3\pp\s$ for $\s=\s^\lt$.

By induction on $\pos_0$, we have $\g^4\cdot\s=\g^4\pp\s$ for all $\s$ in
$\pos_0$ of degree at most 5. Let $\lt\in I(\pos)$ with
$\deg(\lt\cap \pos_0)\leq6$. If $\deg(\lt\cap \pos_0)<6$, then by
induction on $\pos_0$ we have $\g^4\cdot\s^\lt=\g^4\pp\s^\lt$. If
$\deg(\lt\cap \pos_0)=6$, then we have $\lt\supset\lt_4$ or
$\lt\supset\mu_4$ for all $\lt$ except one and by recursion with respect
to $\s^4$ or $\tau^4$ we get equality $\g^4\cdot\s^\lt=\g^4\pp\s^\lt$ for 
those $\lt$. Equation $\g^4\cdot\s^\lt=\g^4\pp\s^\lt$ holds for the last 
ideal $\lt$ by Lemma \ref{tous-sauf-1}. Finally, let $\lt\subset \pos_0$ 
with $\deg(\lt)\geq7$. We have $\lt\supset\lt_4$ or $\lt\supset\mu_4$ and
by recursion with respect to $\s^4$ or $\tau^4$ we get equality
$\g^4\cdot\s^\lt=\g^4\pp\s^\lt$ for $\s^\lt$.
\end{proo}

\subsection{Type $E_7$ case}

Let us start with the case of $E_7/P_1$.
Thus we consider the system of $\varpi_1$-minuscule
$E_7$-colored posets $\pos_0$ given by the
following maximal elements:

$$
\begin{array}{cccc}
%auto-ignore

\psset{unit=6mm}

\begin{pspicture*}(-1.5,0.5)(6,-8.5)

\psline(0.12,-7.88)(0.88,-7.12)
\psline(1.12,-6.88)(1.88,-6.12)
\psline(2,-5.84)(2,-5.16)
\psline(2.12,-5.88)(2.88,-5.12)
\psline(2,-4.84)(2,-4.16)
\psline(2.88,-4.88)(2.12,-4.12)
\psline(3.12,-4.88)(3.88,-4.12)
\psline(1.88,-3.88)(1.12,-3.12)
\psline(2.12,-3.88)(2.88,-3.12)
\psline(3.88,-3.88)(3.12,-3.12)
\psline(4.12,-3.88)(4.88,-3.12)
\psline(0.88,-2.88)(0.12,-2.12)
\psline(1.12,-2.88)(1.88,-2.12)
\psline(2.88,-2.88)(2.12,-2.12)
\psline(3.12,-2.88)(3.88,-2.12)
\psline(4.88,-2.88)(4.12,-2.12)
\psline(0.12,-1.88)(0.88,-1.12)
\psline(1.88,-1.88)(1.12,-1.12)
\psline(2.12,-1.88)(2.88,-1.12)
\psline(3.88,-1.88)(3.12,-1.12)
\psline(1.12,-0.88)(1.88,-0.12)
\psline(2.88,-0.88)(2.12,-0.12)

\pscircle*(0,-8){0.16}
\pscircle(0,-2){0.16}  \psline(-0.12,-1.88)(0.12,-2.12)  \psline(-0.12,-2.12)(0.12,-1.88)   \put(-0.958,-2.125){$\sigma^8$}
\pscircle*(1,-7){0.16}
\pscircle(1,-3){0.16}  \psline(0.88,-2.88)(1.12,-3.12)  \psline(0.88,-3.12)(1.12,-2.88)
\pscircle(1,-1){0.16}
\pscircle*(2,-6){0.16}
\psdiamond[fillstyle=solid,fillcolor=black](2,-5)(0.256,0.256)   \put(1.042,-5.125){$\sigma^4$}
\pscircle(2,-4){0.16}  \psline(1.88,-3.88)(2.12,-4.12)  \psline(1.88,-4.12)(2.12,-3.88)
\pscircle(2,-2){0.16}
\pscircle(2,0){0.16}
\pscircle*(3,-5){0.16}
\pscircle(3,-3){0.16}  \psline(2.88,-2.88)(3.12,-3.12)  \psline(2.88,-3.12)(3.12,-2.88)
\pscircle(3,-1){0.16}
\pscircle*(4,-4){0.16}
\pscircle(4,-2){0.16}  \psline(3.88,-1.88)(4.12,-2.12)  \psline(3.88,-2.12)(4.12,-1.88)
\psdiamond[fillstyle=solid,fillcolor=black](5,-3)(0.256,0.256)   \put(5.374,-3.125){$\tau^6$}

\end{pspicture*} &
%auto-ignore

\psset{unit=6mm}

\begin{pspicture*}(-0.3,0.5)(5.3,-7.5)

\psline(0.12,-6.88)(0.88,-6.12)
\psline(1.12,-5.88)(1.88,-5.12)
\psline(2,-4.84)(2,-4.16)
\psline(2.12,-4.88)(2.88,-4.12)
\psline(2,-3.84)(2,-3.16)
\psline(2.88,-3.88)(2.12,-3.12)
\psline(3.12,-3.88)(3.88,-3.12)
\psline(1.88,-2.88)(1.12,-2.12)
\psline(2.12,-2.88)(2.88,-2.12)
\psline(3.88,-2.88)(3.12,-2.12)
\psline(4.12,-2.88)(4.88,-2.12)
\psline(0.88,-1.88)(0.12,-1.12)
\psline(1.12,-1.88)(1.88,-1.12)
\psline(2.88,-1.88)(2.12,-1.12)
\psline(3.12,-1.88)(3.88,-1.12)
\psline(4.88,-1.88)(4.12,-1.12)
\psline(0.12,-0.88)(0.88,-0.12)
\psline(1.88,-0.88)(1.12,-0.12)
\psline(2,-0.84)(2,-0.16)
\psline(2.12,-0.88)(2.88,-0.12)
\psline(3.88,-0.88)(3.12,-0.12)

\pscircle*(0,-7){0.16}
\pscircle(0,-1){0.16}  \psline(-0.12,-0.88)(0.12,-1.12)
\psline(-0.12,-1.12)(0.12,-0.88) 
\pscircle*(1,-6){0.16}
\pscircle(1,-2){0.16}  
\psline(0.88,-1.88)(1.12,-2.12)
\psline(0.88,-2.12)(1.12,-1.88) 
\pscircle(1,0){0.16}
\pscircle*(2,-5){0.16}
\psdiamond[fillstyle=solid,fillcolor=black](2,-4)(0.256,0.256)
\pscircle(2,-3){0.16}  \psline(1.88,-2.88)(2.12,-3.12)
\psline(1.88,-3.12)(2.12,-2.88) 
\pscircle(2,-1){0.16}
\pscircle(2,0){0.16}
\pscircle*(3,-4){0.16}
\pscircle(3,-2){0.16}  \psline(2.88,-1.88)(3.12,-2.12)
\psline(2.88,-2.12)(3.12,-1.88) 
\pscircle(3,0){0.16}  
%\psline(2.88,0.12)(3.12,-0.12)
%\psline(2.88,-0.12)(3.12,0.12) 
\pscircle*(4,-3){0.16}
\pscircle(4,-1){0.16}  
\psline(3.88,-0.88)(4.12,-1.12)
\psline(3.88,-1.12)(4.12,-0.88) 
\psdiamond[fillstyle=solid,fillcolor=black](5,-2)(0.256,0.256)

\end{pspicture*} &
%auto-ignore

\psset{unit=6mm}

\begin{pspicture*}(-0.5,0.5)(5.5,-8.5)

\psline(0.12,-7.88)(0.88,-7.12)
\psline(1.12,-6.88)(1.88,-6.12)
\psline(2,-5.84)(2,-5.16)
\psline(2.12,-5.88)(2.88,-5.12)
\psline(2,-4.84)(2,-4.16)
\psline(2.88,-4.88)(2.12,-4.12)
\psline(3.12,-4.88)(3.88,-4.12)
\psline(1.88,-3.88)(1.12,-3.12)
\psline(2.12,-3.88)(2.88,-3.12)
\psline(3.88,-3.88)(3.12,-3.12)
\psline(4.12,-3.88)(4.88,-3.12)
\psline(0.88,-2.88)(0.12,-2.12)
\psline(1.12,-2.88)(1.88,-2.12)
\psline(2.88,-2.88)(2.12,-2.12)
\psline(3.12,-2.88)(3.88,-2.12)
\psline(4.88,-2.88)(4.12,-2.12)
\psline(2,-1.84)(2,-1.16)
\psline(2.12,-1.88)(2.88,-1.12)
\psline(3.88,-1.88)(3.12,-1.12)
\psline(2,-0.84)(2,-0.16)
\psline(2.88,-0.88)(2.12,-0.12)

\pscircle*(0,-8){0.16}
\pscircle(0,-2){0.16}  \psline(-0.12,-1.88)(0.12,-2.12)  \psline(-0.12,-2.12)(0.12,-1.88)
\pscircle*(1,-7){0.16}
\pscircle(1,-3){0.16}  \psline(0.88,-2.88)(1.12,-3.12)  \psline(0.88,-3.12)(1.12,-2.88)
\pscircle*(2,-6){0.16}
\psdiamond[fillstyle=solid,fillcolor=black](2,-5)(0.256,0.256)
\pscircle(2,-4){0.16}  \psline(1.88,-3.88)(2.12,-4.12)  \psline(1.88,-4.12)(2.12,-3.88)
\pscircle(2,-2){0.16}
\pscircle(2,-1){0.16}
\pscircle(2,0){0.16}
\pscircle*(3,-5){0.16}
\pscircle(3,-3){0.16}  \psline(2.88,-2.88)(3.12,-3.12)  \psline(2.88,-3.12)(3.12,-2.88)
\pscircle(3,-1){0.16}
\pscircle*(4,-4){0.16}
\pscircle(4,-2){0.16}  \psline(3.88,-1.88)(4.12,-2.12)  \psline(3.88,-2.12)(4.12,-1.88)
\psdiamond[fillstyle=solid,fillcolor=black](5,-3)(0.256,0.256)

\end{pspicture*} & 
%auto-ignore

\psset{unit=6mm}

\begin{pspicture*}(-0.5,0.5)(5.5,-8.5)

\psline(0.12,-7.88)(0.88,-7.12)
\psline(1.12,-6.88)(1.88,-6.12)
\psline(2,-5.84)(2,-5.16)
\psline(2.12,-5.88)(2.88,-5.12)
\psline(2,-4.84)(2,-4.16)
\psline(2.88,-4.88)(2.12,-4.12)
\psline(3.12,-4.88)(3.88,-4.12)
\psline(1.88,-3.88)(1.12,-3.12)
\psline(2.12,-3.88)(2.88,-3.12)
\psline(3.88,-3.88)(3.12,-3.12)
\psline(4.12,-3.88)(4.88,-3.12)
\psline(0.88,-2.88)(0.12,-2.12)
\psline(1.12,-2.88)(1.88,-2.12)
\psline(2.88,-2.88)(2.12,-2.12)
\psline(3.12,-2.88)(3.88,-2.12)
\psline(4.88,-2.88)(4.12,-2.12)
\psline(0.12,-1.88)(0.88,-1.12)
\psline(1.88,-1.88)(1.12,-1.12)
\psline(2,-1.84)(2,-1.16)
\psline(1.12,-0.88)(1.88,-0.12)
\psline(2,-0.84)(2,-0.16)

\pscircle*(0,-8){0.16}
\pscircle(0,-2){0.16}  \psline(-0.12,-1.88)(0.12,-2.12)  \psline(-0.12,-2.12)(0.12,-1.88)
\pscircle*(1,-7){0.16}
\pscircle(1,-3){0.16}  \psline(0.88,-2.88)(1.12,-3.12)  \psline(0.88,-3.12)(1.12,-2.88)
\pscircle(1,-1){0.16}
\pscircle*(2,-6){0.16}
\psdiamond[fillstyle=solid,fillcolor=black](2,-5)(0.256,0.256)
\pscircle(2,-4){0.16}  \psline(1.88,-3.88)(2.12,-4.12)  \psline(1.88,-4.12)(2.12,-3.88)
\pscircle(2,-2){0.16}
\pscircle(2,-1){0.16}
\pscircle(2,0){0.16}
\pscircle*(3,-5){0.16}
\pscircle(3,-3){0.16}  \psline(2.88,-2.88)(3.12,-3.12)  \psline(2.88,-3.12)(3.12,-2.88)
\pscircle*(4,-4){0.16}
\pscircle(4,-2){0.16}  \psline(3.88,-1.88)(4.12,-2.12)  \psline(3.88,-2.12)(4.12,-1.88)
\psdiamond[fillstyle=solid,fillcolor=black](5,-3)(0.256,0.256)

\end{pspicture*} \\
\end{array}
$$

We have $S_0 = \{2,7\}$. For $i \in \{2,7\}$
let $(D_i,d_i)$ be a marked Dynkin diagram and $P_i$ be any $d_i$-minuscule
$D_i$-colored poset. Set $\pos = \pos_{\pos_0,(P_2,P_7)}$
with notation \ref{nota-system}.

\begin{lemm}
\label{lemm-e7p1}
With the above notation, assume that Conjecture \ref{main_conj}
holds for $P_2$,$P_7$, and any $\lambda$ in $I(\pos)$ with
$D_0(\lambda) \varsubsetneq E_7$. Then Conjecture \ref{main_conj}
holds for $\pos$.
\end{lemm}

\begin{proo}
We consider the ideals
$\l_4 = \scal{(\alpha_2,1)}$ resp. $\mu_6=\scal{(\a_7,1)}$ and 
$\lt_8=\scal{(\a_1,2)}$ in $\pos_0$, of degree 4 resp. 6 and 8. 
The corresponding 
Schubert cells are denoted with $\s^4$  resp. $\tau^6$ and $\s^8$. Let 
$\{ \g^1,\g^4,\g^6 \}$ be a set of generators of the cohomology ring of
$E_7/P_1$, with $\deg(\g^i)=i$.
The variety $E_7/P_1$ has dimension 33 and the
dimensions of $H^d(E_7/P_1)$ are

\vskip 0.2 cm

\centerline{\begin{tabular}{c|llllllllllllllllllllllllllllllll}
\hline
$d$&0&1&2&3&4&5&6&7&8&9&10&11&12&13&14&15&16\\
\hline
$\dim H^d(E_7/P_1)$&1&1&1&1&2&2&3&3&4&4&5&5&6&6&6&6&7\\
\hline
\end{tabular}}

\vskip 0.2 cm

Since by assumption the conjecture holds for any
$\lambda \in I(\pos)$ with $D_0(\lambda) \varsubsetneq E_7$,
we have $\oknu{\l}{\mu}{\nu}$ as soon as
$\deg(\nu \cap \pos_0) \leq 11$.

By Proposition \ref{prop-chevalley}, $\g^1$ is a good generator.
By the above argument, 
we have $\g^4\cdot\s=\g^4\pp\s$ for all $\s$ in
$\pos_0$ of degree at most 7. 
In particular this is valid for $\s=\s^4$, for $\tau^6$, for any class 
$\s^\lt$ with $\lt\not\supset\lt_4$ and for any class 
$\s^\lt$ with $\lt\in I(\pos)-I(\pos_0)$. By recursion with respect 
to $\s^4$, 
an ideal $\nu$ such that $c_{\lt_4,\lt}^\nu\neq t_{\lt_4,\lt}^\nu$ 
has to contain $(\a_2,2)$ and in particular $\deg(\nu\cap \pos_0)\geq14$ (we 
also use induction with respect to $\pos_0$). We thus have 
$\g^4\cdot\s=\g^4\pp\s$ for any class $\s$ with $\deg(\s)\leq 9$. By 
recursion with respect to $\tau^6$ or $\s^8$ we have 
$\g^4\cdot\s^\lt=\g^4\pp\s^\lt$ for any ideal $\lt\supset\mu_6$ 
or $\lt\supset\lt_8$. As there is only one ideal $\lt\in \pos_0$ 
with $\lt\not\supset\mu_6$ and $\lt\not\supset\lt_8$ in degree 10 and 11 and 
none in higher degree, we have $\g^4\cdot\s=\g^4\pp\s$ for any $\s\in 
H^*(\pos)$ by Lemma \ref{tous-sauf-1}.

By induction on $\pos_0$, we have $\g^6\cdot\s=\g^6\pp\s$ for all $\s$ in
$\pos_0$ of degree at most 5. Let $\s^\lt\in H^*(\pos)$ of degree 6 and 
$\s^\lt\neq\tau^6$, then $\lt\supset\lt_4$ and by recursion with respect to 
$\s^4$ we have $c_{\mu_6,\lt}^\nu=t_{\mu_6,\lt}^\nu$ for $\nu\not\ni(\a_2,2)$ 
(or $\nu\not\ni(\a_6,2)$ or $\nu\not\ni(\a_1,2)$ by induction on $\pos_0$). 
But for degree reasons we have $\deg(\nu\cap \pos_0)\leq12$ thus 
$\g^6\cdot\s^\lt=\g^6\pp\s^\lt$. By Lemma \ref{tous-sauf-1} we 
obtain $\g^6\cdot\tau^6=\g^6\pp\tau^6$. In particular 
$\g^6\cdot\s=\g^6\pp\s$ for $\s$ with associated ideal $\lt$ in 
$I(\pos)-I(\pos_0)$. By Lemma 
\ref{lefs}, we obtain $\g^6\cdot\s=\g^6\pp\s$ for $\s$ of degree 7. 
As any degree 8 class in $H^*(\pos_0)$ is a linear combination, in 
$H^*(\pos_0)$, of $(\g^4)^{\pp2}$ and multiples of $\g^1$ we conclude by Lemma 
\ref{lemm-sous-algebre}  for degree 8 classes. For degree 9 classes we 
conclude by Lemma \ref{lefs} and for higher degree classes we conclude 
as for $\g^4$.
\end{proo}

We now deal with the case $E_7/P_2$. Thus we consider the system of 
$\varpi_2$-minuscule $E_7$-colored posets $\pos_0$ given by:

\centerline{%auto-ignore

\psset{unit=6mm}

\begin{pspicture*}(1.5,0.5)(-6.5,-5.5)

\psline(-3,-4.84)(-3,-4.16)
\psline(-2.88,-3.88)(-2.12,-3.12)
\psline(-3.12,-3.88)(-3.88,-3.12)
\psline(-1.88,-2.88)(-1.12,-2.12)
\psline(-2.12,-2.88)(-2.88,-2.12)
\psline(-3.88,-2.88)(-3.12,-2.12)
\psline(-4.12,-2.88)(-4.88,-2.12)
\psline(-0.88,-1.88)(-0.12,-1.12)
\psline(-1.12,-1.88)(-1.88,-1.12)
\psline(-2.88,-1.88)(-2.12,-1.12)
\psline(-3,-1.84)(-3,-1.16)
\psline(-3.12,-1.88)(-3.88,-1.12)
\psline(-4.88,-1.88)(-4.12,-1.12)
\psline(-0.12,-0.88)(-0.88,-0.12)
\psline(-1.88,-0.88)(-1.12,-0.12)

\psdiamond[fillstyle=solid,fillcolor=black](-0,-1)(0.256,0.256)   \put(0.374,-1.125){$\tau^5$}
\pscircle*(-1,-2){0.16}
\pscircle(-1,-0){0.16}  \psline(-0.88,0.12)(-1.12,-0.12)  \psline(-0.88,-0.12)(-1.12,0.12)
\pscircle*(-2,-3){0.16}
\pscircle(-2,-1){0.16}  \psline(-1.88,-0.88)(-2.12,-1.12)  \psline(-1.88,-1.12)(-2.12,-0.88)   \put(-2.208,-0.526){$\sigma^7$}
\pscircle*(-3,-5){0.16}
\pscircle*(-3,-4){0.16}
\pscircle(-3,-2){0.16}  \psline(-2.88,-1.88)(-3.12,-2.12)  \psline(-2.88,-2.12)(-3.12,-1.88)
\pscircle(-3,-1){0.16}  \psline(-2.88,-0.88)(-3.12,-1.12)  \psline(-2.88,-1.12)(-3.12,-0.88)   \put(-3.208,-0.526){$\sigma^6$}
\pscircle*(-4,-3){0.16}   \put(-4.958,-3.125){$\sigma^3$}
\pscircle(-4,-1){0.16}  \psline(-3.88,-0.88)(-4.12,-1.12)  \psline(-3.88,-1.12)(-4.12,-0.88)
\psdiamond[fillstyle=solid,fillcolor=black](-5,-2)(0.256,0.256)   \put(-5.958,-2.125){$\sigma^4$}

\end{pspicture*}}

We have $S_0 = \{1,7\}$. For $i \in \{1,7\}$
let $(D_i,d_i)$ be a marked Dynkin diagram and $P_i$ be any $d_i$-minuscule
$D_i$-colored poset. Set $\pos = \pos_{\pos_0,(P_1,P_7)}$
with notation \ref{nota-system}.

\begin{lemm}
\label{lemm-e7p2}
With the above notation, assume that Conjecture \ref{main_conj}
holds for $P_1$,$P_7$, and any $\lambda$ in $I(\pos)$ with
$D_0(\lambda) \varsubsetneq E_7$. Then Conjecture \ref{main_conj}
holds for $\pos$.
\end{lemm}

\begin{proo}
We consider the ideals
$\l_4 = \scal{(\alpha_1,1)}$ resp. $\mu_5=\scal{(\a_7,1)}$ and 
$\lt_6=\scal{(\a_2,2)}$ in $\pos_0$, of degree 4 resp. 5 and 6. 
The corresponding 
Schubert cells are denoted with $\s^4$  resp. $\tau^5$ and $\s^6$. Let 
$\{ \g^1,\g^3,\g^4,\g^5,\g^7 \}$ be a set of generators of the cohomology ring
of $E_7/P_2$, with
$\deg(\g^i)=i$.

\vskip 0.2 cm

Since by assumption the conjecture holds for any
$\lambda \in I(\pos)$ with $D_0(\lambda) \varsubsetneq E_7$,
we have $\oknu{\l}{\mu}{\nu}$ as soon as
$\deg(\nu \cap \pos_0) \leq 11$.

By Proposition \ref{prop-chevalley}, $\g^1$ is a good generator.
By the above argument, 
we have $\g^3\cdot\s=\g^3\pp\s$ for all $\s$ in
$\pos_0$ of degree at most 8. In particular this equation is valid for 
$\s=\s^4$ and $\tau^5$. As any class $\s^\lt$ of degree at least 9 in 
$H^*(\pos_0)$ satisfies $\lt\supset\lt_4$ or $\lt\supset\mu_5$ we conclude 
by recursion with respect to $\s^4$ or $\tau^5$.

We know that
$\g^4\cdot\s^\lt=\g^4\pp\s^\lt$ for all  
$\lt$ in $I(\pos)$ with $\deg(\lt \cap \pos_0)\leq7$. In particular this 
equation  is valid for  $\lt=\lt_4$ and $\lt=\mu_5$ and thus for any class 
$\s^\lt$ with $\lt\in I(\pos)-I(\pos_0)$. As any class $\s^\lt$ of degree 
at least 8 in 
$H^*(\pos_0)$ satisfies $\lt\supset\lt_4$ or $\lt\supset\mu_5$ except one in 
degree 8, we conclude  by recursion with respect to $\s^4$ or $\tau^5$ and 
Lemma \ref{tous-sauf-1}.

We know that 
$\tau^5\cdot\s^\lt=\tau^5\pp\s^\lt$ for 
all $\lt$ in $I(\pos)$ with $\deg(\lt \cap \pos_0)\leq6$. In particular this 
equation  is valid for $\lt=\lt_4$ and $\lt=\mu_5$ and thus for any class 
$\s^\lt$ with $\lt\in I(\pos)-I(\pos_0)$. As any class $\s^\lt$ of degree 
at least 7 in 
$H^*(\pos_0)$ satisfies $\lt\supset\lt_4$ or $\lt\supset\mu_5$ except one in 
degree 7 and one in degree 8, we conclude  by recursion with respect 
to $\s^4$ or $\tau^5$ and Lemma \ref{tous-sauf-1}.

By what we already did and Lemma \ref{lemm-sous-algebre}, we have 
$\g^7\cdot\s=\g^7\pp\s$ for all $\s$ in $H^*(\pos_0)$ of degree at most 6.
In particular this equation 
is valid for  $\s=\s^4$, $\s=\tau^5$ and $\s=\s^6$. As a consequence, it is 
also valid for any class $\s^\lt$ with $\lt\in I(\pos)-I(\pos_0)$.
As any class $\s^\lt$ of degree at least 7 in 
$H^*(\pos_0)$ satisfies $\lt\supset\lt_4$, $\lt\supset\mu_5$ or 
$\lt\supset\lt_8$, we conclude  by recursion with respect 
to $\s^4$, $\tau^5$ or $\s^6$.
\end{proo}

We now deal with the case of $E_7/P_7$. Thus $\pos_0 $ 
contains only one poset which is the following:

\centerline{\begin{pspicture*}(0,0)(4.8,10.8)

\psellipse(1.2,7.2)(0.096,0.096)
\pspolygon[fillstyle=solid,fillcolor=black](1.32,3.6)(1.2,3.72)(1.08,3.6)(1.2,3.48)
\psellipse(1.8,7.8)(0.096,0.096)
\psellipse(1.8,6.6)(0.096,0.096)
\psellipse(1.8,4.2)(0.096,0.096)
\psline(1.732,4.132)(1.868,4.268)
\psline(1.732,4.268)(1.868,4.132)
\psellipse[fillstyle=solid,fillcolor=black](1.8,3.0)(0.096,0.096)
\psellipse(2.4,8.4)(0.096,0.096)
\psellipse(2.4,7.8)(0.096,0.096)
\psellipse(2.4,7.2)(0.096,0.096)
\psellipse(2.4,6.0)(0.096,0.096)
\psellipse(2.4,5.4)(0.096,0.096)
\psellipse(2.4,4.8)(0.096,0.096)
\psellipse(2.4,3.6)(0.096,0.096)
\psline(2.332,3.532)(2.468,3.668)
\psline(2.332,3.668)(2.468,3.532)
\pspolygon[fillstyle=solid,fillcolor=black](2.52,3.0)(2.4,3.12)(2.28,3.0)(2.4,2.88)
\psellipse[fillstyle=solid,fillcolor=black](2.4,2.4)(0.096,0.096)
\psellipse(3.0,9.0)(0.096,0.096)
\psellipse(3.0,6.6)(0.096,0.096)
\psellipse(3.0,5.4)(0.096,0.096)
\psellipse(3.0,4.2)(0.096,0.096)
\psline(2.932,4.132)(3.068,4.268)
\psline(2.932,4.268)(3.068,4.132)
\psellipse[fillstyle=solid,fillcolor=black](3.0,1.8)(0.096,0.096)
\psellipse(3.6,9.6)(0.096,0.096)
\psellipse(3.6,6.0)(0.096,0.096)
\psellipse(3.6,4.8)(0.096,0.096)
\psline(3.532,4.732)(3.668,4.868)
\psline(3.532,4.868)(3.668,4.732)
\psellipse[fillstyle=solid,fillcolor=black](3.6,1.2)(0.096,0.096)
\psellipse(4.2,10.2)(0.096,0.096)
\psellipse(4.2,5.4)(0.096,0.096)
\psline(4.132,5.332)(4.268,5.468)
\psline(4.132,5.468)(4.268,5.332)
\psellipse[fillstyle=solid,fillcolor=black](4.2,0.6)(0.096,0.096)
\psline(1.268,7.132)(1.732,6.668)
\psline(1.26,3.54)(1.732,3.068)
\psline(1.868,7.732)(2.332,7.268)
\psline(1.732,7.732)(1.268,7.268)
\psline(1.868,6.532)(2.332,6.068)
\psline(1.868,4.132)(2.332,3.668)
\psline(1.732,4.132)(1.26,3.66)
\psline(1.868,2.932)(2.332,2.468)
\psline(2.4,8.304)(2.4,7.896)
\psline(2.332,8.332)(1.868,7.868)
\psline(2.4,7.704)(2.4,7.296)
\psline(2.468,7.132)(2.932,6.668)
\psline(2.332,7.132)(1.868,6.668)
\psline(2.468,5.932)(2.932,5.468)
\psline(2.4,5.904)(2.4,5.496)
\psline(2.4,5.304)(2.4,4.896)
\psline(2.468,4.732)(2.932,4.268)
\psline(2.332,4.732)(1.868,4.268)
\psline(2.4,3.504)(2.4,3.12)
\psline(2.332,3.532)(1.868,3.068)
\psline(2.4,2.88)(2.4,2.496)
\psline(2.468,2.332)(2.932,1.868)
\psline(2.932,8.932)(2.468,8.468)
\psline(3.068,6.532)(3.532,6.068)
\psline(2.932,6.532)(2.468,6.068)
\psline(3.068,5.332)(3.532,4.868)
\psline(2.932,5.332)(2.468,4.868)
\psline(2.932,4.132)(2.468,3.668)
\psline(3.068,1.732)(3.532,1.268)
\psline(3.532,9.532)(3.068,9.068)
\psline(3.668,5.932)(4.132,5.468)
\psline(3.532,5.932)(3.068,5.468)
\psline(3.532,4.732)(3.068,4.268)
\psline(3.668,1.132)(4.132,0.668)
\psline(4.132,10.132)(3.668,9.668)
\psline(4.132,5.332)(3.668,4.868)
\uput[r](2.4,3.0){$\s^5$}
\uput[l](1.2,3.6){$\tau^6$}

\end{pspicture*}}

We have $S_0 = \{1,2\}$. For $i \in \{1,2\}$
let $(D_i,d_i)$ be a marked Dynkin diagram and $P_i$ be any $d_i$-minuscule
$D_i$-colored poset. Set $\pos = \pos_{\pos_0,(P_1,P2)}$
with notation \ref{nota-system}.

\begin{lemm}
\label{lemm-e7p7}
With the above notation, assume that Conjecture \ref{main_conj}
holds for $P_1$,$P_2$, and any $\lambda$ in $I(\pos)$ with
$D_0(\lambda) \varsubsetneq E_7$. Then Conjecture \ref{main_conj}
holds for $\pos$.
\end{lemm}
\begin{proo}
Let $\s^5$ resp. $\tau^6$ 
be the Schubert classes corresponding
to the ideals generated by $(\a_2,1)$ resp. $(\a_7,1)$. 
They are of degree $5$ resp. $6$. 
Let $\{ \g^1,\g^5,\g^9 \}$ be a set of generators of $H^*(E_7/P_7)$,
where $\g^i$ has degree $i$.
The variety $E_7/P_7$ has dimension 27
and the dimensions of $H^d(E_7/P_7)$ are

\vskip 0.2 cm

\centerline{\begin{tabular}{c|lllllllllllllllllllllllllllll}
\hline
$d$&0&1&2&3&4&5&6&7&8&9&10&11&12&13\\
\hline
$\dim H^d(E_7/P_7)$&1&1&1&1&1&2&2&2&2&3&3&3&3&3\\
\hline
\end{tabular}}

\vskip 0.2 cm

By Proposition \ref{prop-chevalley}, $\g^1$ is a good generator.
If $\notoknu{\l}{\mu}{\nu}$, then $\nu \cap \pos_0$ must have degree at
least 12. Thus $\ok{\g^5}{\s}$ if $\deg(\s) \leq 6$. By recursion with
respect to $\s^5$ and $\tau^6$ we have $\ok{\g^5}{\s^\l}$
if $\l \in I(\pos) - I(\pos_0)$. Thus by Lemma \ref{lefs} we have
$\ok{\g^5}{\s}$ if $\deg(\s) \leq 8$.
Let $\mu,\nu \in I(\pos_0)$ such that $\notoknu{\g^5}{\mu}{\nu}$.
Assume $\deg(\mu) = 9$. If $(\a_1,1) \in \mu$
then by recursion with respect to $\tau^6$ we have $(\a_1,2) \in \nu$,
thus $\deg(\nu) \geq 18$, a contradiction. Thus $\mu$ cannot
contain $(\a_1,1)$. Since there is a unique ideal in $\pos_0$
of degree $9$ not containing
$(\a_1,1)$ (namely $\scal{(\a_6,2)}$), we conclude by Lemma
\ref{tous-sauf-1} that $\ok{\g^5}{\s}$ if $\deg(\s)=9$. By Lemma \ref{lefs},
$\g^5$ is a good generator.

Since we know that $\ok{\g^9}{\g^5}$, by Lemma \ref{lemm-sous-algebre}
we deduce that $\ok{\g^9}{\s}$ if $\s$ is in the subalgebra generated by
$\g^1$ and $\g^5$ in $H^*(E_7/P_7)$ and in particular for
$\deg(\s)\leq 9$. By recursion with respect to $\s^5$ and $\tau^6$ we
have $\ok{\g^9}{\s^\l}$ if $\l \in I(\pos) - I(\pos_0)$. Let $\mu,\nu
\in I(\pos_0)$ such that $\notoknu{\g^9}{\mu}{\nu}$. Assume $\deg(\mu)
= 9$. 
If $(\a_1,1) \in \mu$ then by recursion with respect to $\tau^6$ we
have $(\a_1,2) \in \nu$, thus $\nu=\scal{(\a_1,2),(\a_7,2)}$. Since
there is a unique ideal in $\pos_0$ of degree $9$ not containing
$(\a_1,1)$ (namely $\scal{(\a_6,2)}$), we conclude by Lemma
\ref{tous-sauf-1} that $\oknu{\g^5}{\s}{\nu}$ except for
$\nu=\scal{(\a_1,2),(\a_7,2)}$. For $\nu=\scal{(\a_1,2),(\a_7,2)}$, we
only need to compute in $H^*(\pos_0)$ and because $\pos_0$ is a
complete $d$-poset we conclude by Lemma \ref{poincare}.
Then we conclude
that $\g^9$ is a good generator by Lemma \ref{lefs}. 
\end{proo}

\subsection{Type $E_8$ case}

Let us start with the case of $E_8/P_1$.
Thus we consider the system of $\varpi_1$-minuscule
$E_8$-colored posets $\pos_0$ given by the three
following maximal elements and their obvious intersections:

$$
\begin{array}{ccc}
%auto-ignore

\psset{unit=6mm}

\begin{pspicture*}(-1.5,0.5)(7.5,-9.5)

\psline(0.12,-8.88)(0.88,-8.12)
\psline(1.12,-7.88)(1.88,-7.12)
\psline(2,-6.84)(2,-6.16)
\psline(2.12,-6.88)(2.88,-6.12)
\psline(2,-5.84)(2,-5.16)
\psline(2.88,-5.88)(2.12,-5.12)
\psline(3.12,-5.88)(3.88,-5.12)
\psline(1.88,-4.88)(1.12,-4.12)
\psline(2.12,-4.88)(2.88,-4.12)
\psline(3.88,-4.88)(3.12,-4.12)
\psline(4.12,-4.88)(4.88,-4.12)
\psline(0.88,-3.88)(0.12,-3.12)
\psline(1.12,-3.88)(1.88,-3.12)
\psline(2.88,-3.88)(2.12,-3.12)
\psline(3.12,-3.88)(3.88,-3.12)
\psline(4.88,-3.88)(4.12,-3.12)
\psline(5.12,-3.88)(5.88,-3.12)
\psline(0.12,-2.88)(0.88,-2.12)
\psline(1.88,-2.88)(1.12,-2.12)
\psline(2.12,-2.88)(2.88,-2.12)
\psline(3.88,-2.88)(3.12,-2.12)
\psline(4.12,-2.88)(4.88,-2.12)
\psline(5.88,-2.88)(5.12,-2.12)
\psline(1.12,-1.88)(1.88,-1.12)
\psline(2.88,-1.88)(2.12,-1.12)
\psline(3.12,-1.88)(3.88,-1.12)
\psline(4.88,-1.88)(4.12,-1.12)
\psline(2.12,-0.88)(2.88,-0.12)
\psline(3.88,-0.88)(3.12,-0.12)

\pscircle*(0,-9){0.16}
\pscircle(0,-3){0.16}  \psline(-0.12,-2.88)(0.12,-3.12)  \psline(-0.12,-3.12)(0.12,-2.88)   \put(-0.958,-3.125){$\sigma^8$}
\pscircle*(1,-8){0.16}
\pscircle(1,-4){0.16}  \psline(0.88,-3.88)(1.12,-4.12)
\psline(0.88,-4.12)(1.12,-3.88)   %\put(0.042,-4.125){$\sigma_7$} 
\pscircle(1,-2){0.16}
\pscircle*(2,-7){0.16}
\psdiamond[fillstyle=solid,fillcolor=black](2,-6)(0.256,0.256)   \put(1.042,-6.125){$\sigma^4$}
\pscircle(2,-5){0.16}  \psline(1.88,-4.88)(2.12,-5.12)  \psline(1.88,-5.12)(2.12,-4.88)
\pscircle(2,-3){0.16}
\pscircle(2,-1){0.16}
\pscircle*(3,-6){0.16}
\pscircle(3,-4){0.16}  \psline(2.88,-3.88)(3.12,-4.12)  \psline(2.88,-4.12)(3.12,-3.88)
\pscircle(3,-2){0.16}
\pscircle(3,0){0.16}
\pscircle*(4,-5){0.16}
\pscircle(4,-3){0.16}  \psline(3.88,-2.88)(4.12,-3.12)  \psline(3.88,-3.12)(4.12,-2.88)
\pscircle(4,-1){0.16}
\pscircle*(5,-4){0.16}
\pscircle(5,-2){0.16}  \psline(4.88,-1.88)(5.12,-2.12)  \psline(4.88,-2.12)(5.12,-1.88)
\psdiamond[fillstyle=solid,fillcolor=black](6,-3)(0.256,0.256)   \put(6.374,-3.125){$\tau^7$}

\end{pspicture*} &
%auto-ignore

\psset{unit=6mm}

\begin{pspicture*}(-0.5,0.5)(6.5,-9.5)

\psline(0.12,-8.88)(0.88,-8.12)
\psline(1.12,-7.88)(1.88,-7.12)
\psline(2,-6.84)(2,-6.16)
\psline(2.12,-6.88)(2.88,-6.12)
\psline(2,-5.84)(2,-5.16)
\psline(2.88,-5.88)(2.12,-5.12)
\psline(3.12,-5.88)(3.88,-5.12)
\psline(1.88,-4.88)(1.12,-4.12)
\psline(2.12,-4.88)(2.88,-4.12)
\psline(3.88,-4.88)(3.12,-4.12)
\psline(4.12,-4.88)(4.88,-4.12)
\psline(0.88,-3.88)(0.12,-3.12)
\psline(1.12,-3.88)(1.88,-3.12)
\psline(2.88,-3.88)(2.12,-3.12)
\psline(3.12,-3.88)(3.88,-3.12)
\psline(4.88,-3.88)(4.12,-3.12)
\psline(5.12,-3.88)(5.88,-3.12)
\psline(0.12,-2.88)(0.88,-2.12)
\psline(1.88,-2.88)(1.12,-2.12)
\psline(2,-2.84)(2,-2.16)
\psline(2.12,-2.88)(2.88,-2.12)
\psline(3.88,-2.88)(3.12,-2.12)
\psline(4.12,-2.88)(4.88,-2.12)
\psline(5.88,-2.88)(5.12,-2.12)
\psline(3.12,-1.88)(3.88,-1.12)
\psline(4.88,-1.88)(4.12,-1.12)

\pscircle*(0,-9){0.16}
\pscircle(0,-3){0.16}  \psline(-0.12,-2.88)(0.12,-3.12)  \psline(-0.12,-3.12)(0.12,-2.88)
\pscircle*(1,-8){0.16}
\pscircle(1,-4){0.16}  \psline(0.88,-3.88)(1.12,-4.12)  \psline(0.88,-4.12)(1.12,-3.88)
\pscircle(1,-2){0.16}
\pscircle*(2,-7){0.16}
\psdiamond[fillstyle=solid,fillcolor=black](2,-6)(0.256,0.256)
\pscircle(2,-5){0.16}  \psline(1.88,-4.88)(2.12,-5.12)  \psline(1.88,-5.12)(2.12,-4.88)
\pscircle(2,-3){0.16}
\pscircle(2,-2){0.16}   \put(1.626,-1.609){$\sigma^{11}$}
\pscircle*(3,-6){0.16}
\pscircle(3,-4){0.16}  \psline(2.88,-3.88)(3.12,-4.12)  \psline(2.88,-4.12)(3.12,-3.88)
\pscircle(3,-2){0.16}
\pscircle*(4,-5){0.16}
\pscircle(4,-3){0.16}  \psline(3.88,-2.88)(4.12,-3.12)  \psline(3.88,-3.12)(4.12,-2.88)
\pscircle(4,-1){0.16}
\pscircle*(5,-4){0.16}
\pscircle(5,-2){0.16}  \psline(4.88,-1.88)(5.12,-2.12)  \psline(4.88,-2.12)(5.12,-1.88)
\psdiamond[fillstyle=solid,fillcolor=black](6,-3)(0.256,0.256)

\end{pspicture*} &
%auto-ignore

\psset{unit=6mm}

\begin{pspicture*}(-0.5,0.5)(6.5,-9.5)

\psline(0.12,-8.88)(0.88,-8.12)
\psline(1.12,-7.88)(1.88,-7.12)
\psline(2,-6.84)(2,-6.16)
\psline(2.12,-6.88)(2.88,-6.12)
\psline(2,-5.84)(2,-5.16)
\psline(2.88,-5.88)(2.12,-5.12)
\psline(3.12,-5.88)(3.88,-5.12)
\psline(1.88,-4.88)(1.12,-4.12)
\psline(2.12,-4.88)(2.88,-4.12)
\psline(3.88,-4.88)(3.12,-4.12)
\psline(4.12,-4.88)(4.88,-4.12)
\psline(0.88,-3.88)(0.12,-3.12)
\psline(1.12,-3.88)(1.88,-3.12)
\psline(2.88,-3.88)(2.12,-3.12)
\psline(3.12,-3.88)(3.88,-3.12)
\psline(4.88,-3.88)(4.12,-3.12)
\psline(5.12,-3.88)(5.88,-3.12)
\psline(2,-2.84)(2,-2.16)
\psline(2.12,-2.88)(2.88,-2.12)
\psline(3.88,-2.88)(3.12,-2.12)
\psline(4.12,-2.88)(4.88,-2.12)
\psline(5.88,-2.88)(5.12,-2.12)
\psline(2,-1.84)(2,-1.16)
\psline(2.88,-1.88)(2.12,-1.12)
\psline(3.12,-1.88)(3.88,-1.12)
\psline(4.88,-1.88)(4.12,-1.12)
\psline(2.12,-0.88)(2.88,-0.12)
\psline(3.88,-0.88)(3.12,-0.12)

\pscircle*(0,-9){0.16}
\pscircle(0,-3){0.16}  \psline(-0.12,-2.88)(0.12,-3.12)  \psline(-0.12,-3.12)(0.12,-2.88)
\pscircle*(1,-8){0.16}
\pscircle(1,-4){0.16}  \psline(0.88,-3.88)(1.12,-4.12)  \psline(0.88,-4.12)(1.12,-3.88)
\pscircle*(2,-7){0.16}
\psdiamond[fillstyle=solid,fillcolor=black](2,-6)(0.256,0.256)
\pscircle(2,-5){0.16}  \psline(1.88,-4.88)(2.12,-5.12)  \psline(1.88,-5.12)(2.12,-4.88)
\pscircle(2,-3){0.16}
\pscircle(2,-2){0.16}
\pscircle(2,-1){0.16}
\pscircle*(3,-6){0.16}
\pscircle(3,-4){0.16}  \psline(2.88,-3.88)(3.12,-4.12)  \psline(2.88,-4.12)(3.12,-3.88)
\pscircle(3,-2){0.16}
\pscircle(3,0){0.16}
\pscircle*(4,-5){0.16}
\pscircle(4,-3){0.16}  \psline(3.88,-2.88)(4.12,-3.12)  \psline(3.88,-3.12)(4.12,-2.88)
\pscircle(4,-1){0.16}
\pscircle*(5,-4){0.16}
\pscircle(5,-2){0.16}  \psline(4.88,-1.88)(5.12,-2.12)  \psline(4.88,-2.12)(5.12,-1.88)
\psdiamond[fillstyle=solid,fillcolor=black](6,-3)(0.256,0.256)

\end{pspicture*}
\end{array}
$$

We have $S_0 = \{2,8\}$. For $i \in \{2,8\}$
let $(D_i,d_i)$ be a marked Dynkin diagram and $P_i$ be any $d_i$-minuscule
$D_i$-colored poset. Set $\pos = \pos_{\pos_0,(P_2,P_8)}$
with notation \ref{nota-system}.

\begin{lemm}
\label{lemm-e8p1}
With the above notation, assume that Conjecture \ref{main_conj}
holds for $P_2$,$P_8$, and any $\lambda$ in $I(\pos)$ with
$D_0(\lambda) \varsubsetneq E_8$. Then Conjecture \ref{main_conj}
holds for $\pos$.
\end{lemm}
\begin{proo}
For $i=4$ resp. $7,8,11$ we consider the ideals
$\l_i = \scal{(\alpha_2,1)}$ resp. $\scal{(\a_8,1)},\scal{(\a_1,2)}$,
$\scal{(\a_2,2)}$ in $\pos_0$, of degree $i$. The corresponding Schubert
cells are denoted with $\s^i$. We denote $\tau^7$ the Schubert cell corresponding
to the ideal $\scal{(\a_8,1)}$, and $\s^{11}$ the cell corresponding to
the ideal $\scal{(\a_2,2)}$ in the two last posets.
Let $\{ \g^1,\g^4,\g^6,\g^7,\g^{10} \}$ be a set of generators of the
cohomology ring of $E_8/P_1$, with
$\deg(\g^i)=i$.

Since by assumption the conjecture holds for any
$\lambda \in I(\pos)$ with $D_0(\lambda) \varsubsetneq E_8$,
we have $\oknu{\l}{\mu}{\nu}$ as soon as
$\deg(\nu \cap \pos_0) \leq 13$.

By Proposition \ref{prop-chevalley}, $\g^1$ is a good generator.
Let $\mu,\nu$ such that $\notoknu{\g^4}{\mu}{\nu}$. By the above we
have $\deg(\nu) \geq 14$. In particular $\ok{\g^4}{\s^4}$. By recursion
with respect to $\s^4$ we deduce that $\nu$ must contain $(\alpha_2,2)$.
Thus $\deg(\nu) \geq 16$. Thus $\ok{\g^4}{\s^{11}}$. By recursion with
respect to $\tau^7,\s^8$ and $\s^{11}$ we get that $\mu$
does not contain these elements. Since moreover $\mu$ must have degree
at least 12 it follows that $\mu$ is one of the two elements
$\scal{(\a_5,3)}$, $\scal{(\a_4,3),(\a_6,2)}$, of degree respectively
13,12. Thus we can conclude by Lemma \ref{tous-sauf-1} that $\g^4$
is a good generator.

Let us show that $\g^6$ is a good generator.
Let $\mu,\nu$ such that $\notoknu{\g^6}{\mu}{\nu}$.
We know that $\ok{\g^6}{\s}$ for $\deg(\s) \leq 7$ thus for
$\s \in \{ \s^4 , \tau^7 \}$. By recursion with respect to these elements
we deduce that $\mu$ cannot contain $(\a_8,1)$, thus $(\a_2,1) \in \mu$,
and
$\nu$ must contain $(\a_2,2)$. Thus $\deg(\nu) \geq 16$ and
$\ok{\g^6}{\s}$ if $\deg(\s) \leq 9$. The number of Schubert classes
of degree 9 resp. 10,11,12,13,14,15,16 not bigger than $\s^8$ and $\tau^7$
is 3 resp. 3,3,2,2,1,1,0, and moreover the map
induced by the multiplication by $h$ is surjective on this sets. Thus
we conclude thanks to Lemma \ref{lefs}$(\imath\imath)$.

Let $\mu,\nu$ such that $\notoknu{\g^7}{\mu}{\nu}$. Assume first that
$\deg(\mu) = 7$. If $\mu \supset \s^4$ then by recursion with respect
to $\s^4$ it follows that $(\a_2,2) \in \nu$ and $\deg(\nu) \geq 16$,
contradicting $\deg(\mu) = 7$. Since there is only one cell of degree
7 which is not bigger than $\s^4$ (namely $\tau^7$), it follows that
$\ok{\g^7}{\s}$ if $\deg(\s) = 7$. By recursion with respect to $\tau^7$
we also have this property for any $\mu\supset\lt_7$. Thus, again
by recursion with respect to $\s^4$,
$\ok{\g^7}{\s}$ if $\deg(\s) \leq 8$. By recursion with respect to
$\s^8$ we have $(\a_1,2) \not \in \mu$. Since
$h^8(E_8/P_1) = h^9(E_8/P_1) = 5$, we deduce from Lemma \ref{lefs}
that $\ok{\g^7}{\s}$ if $\deg(\s)=9$. Then we can argue as for $\g^6$.

For $\g^{10}$ we already know that $\ok{\g^{10}}{\s}$ if
$(\s)$ is one of the $\g^i$'s or $\s=\s^4$ or $\s=\tau^7$.
By recursion with respect to $\s^4$ and $\tau^7$ we deduce
$\ok{\g^{10}}{\s^\l}$ if $\l \in I(\pos) - I(\pos_0)$, and since
the $\g^i$'s for $i \leq 7$ generate $H^9(\pos_0)$ we have
$\ok{\g^{10}}{\s}$ for $\deg(\s)=9$.
Then we can argue as for $\g^6$ and $\g^7$. 
\end{proo}

We now consider the case of $E_8/P_2$.
Thus we consider the system of $\varpi_2$-minuscule
$E_8$-colored posets $\pos_0$ given by only one quiver $\pos_0$: 

$$
\begin{array}{c}
%auto-ignore

\psset{unit=6mm}

\begin{pspicture*}(1.5,0.5)(-7.5,-6.5)

\psline(-4,-5.84)(-4,-5.16)
\psline(-3.88,-4.88)(-3.12,-4.12)
\psline(-4.12,-4.88)(-4.88,-4.12)
\psline(-2.88,-3.88)(-2.12,-3.12)
\psline(-3.12,-3.88)(-3.88,-3.12)
\psline(-4.88,-3.88)(-4.12,-3.12)
\psline(-5.12,-3.88)(-5.88,-3.12)
\psline(-1.88,-2.88)(-1.12,-2.12)
\psline(-2.12,-2.88)(-2.88,-2.12)
\psline(-3.88,-2.88)(-3.12,-2.12)
\psline(-4,-2.84)(-4,-2.16)
\psline(-4.12,-2.88)(-4.88,-2.12)
\psline(-5.88,-2.88)(-5.12,-2.12)
\psline(-0.88,-1.88)(-0.12,-1.12)
\psline(-1.12,-1.88)(-1.88,-1.12)
\psline(-2.88,-1.88)(-2.12,-1.12)
\psline(-0.12,-0.88)(-0.88,-0.12)
\psline(-1.88,-0.88)(-1.12,-0.12)

\psdiamond[fillstyle=solid,fillcolor=black](-0,-1)(0.256,0.256)   \put(0.374,-1.125){$\sigma^6$} 
\pscircle*(-1,-2){0.16}   %\put(-0.626,-2.125){$\sigma_5$}
\pscircle(-1,-0){0.16}  
%\psline(-0.84,-0)(-1.16,-0)  \psline(-1,0.16)(-1,-0.16)
\psline(-0.888,-0.112)(-1.112,0.112)  \psline(-1.112,-0.112)(-0.888,0.112)
\pscircle*(-2,-3){0.16}   %\put(-1.626,-3.125){$\sigma_4$}
\pscircle(-2,-1){0.16}  \psline(-1.88,-0.88)(-2.12,-1.12)  \psline(-1.88,-1.12)(-2.12,-0.88)
\pscircle*(-3,-4){0.16}   %\put(-2.626,-4.125){$\sigma_3$}
\pscircle(-3,-2){0.16}  \psline(-2.88,-1.88)(-3.12,-2.12)  \psline(-2.88,-2.12)(-3.12,-1.88)
\pscircle*(-4,-6){0.16}
\pscircle*(-4,-5){0.16}
\pscircle(-4,-3){0.16}  \psline(-3.88,-2.88)(-4.12,-3.12)  \psline(-3.88,-3.12)(-4.12,-2.88)
\pscircle(-4,-2){0.16}  \psline(-3.88,-1.88)(-4.12,-2.12)  \psline(-3.88,-2.12)(-4.12,-1.88)   \put(-4.208,-1.609){$\tau^6$}
\pscircle*(-5,-4){0.16}
\pscircle(-5,-2){0.16}  \psline(-4.88,-1.88)(-5.12,-2.12)  \psline(-4.88,-2.12)(-5.12,-1.88)   %\put(-5.208,-1.609){$\sigma_7$}
\psdiamond[fillstyle=solid,fillcolor=black](-6,-3)(0.256,0.256)   \put(-6.958,-3.125){$\tau^4$}

\end{pspicture*}
\end{array}
$$

Set $S=S(D_0)$: we have $S=\{1,8\}$. For $i \in \{1,8\}$
let $(D_i,d_i)$ be a marked Dynkin diagram and $P_i$ be any $d_i$-minuscule
$D_i$-colored poset. Set $\pos = \pos_{\pos_0,(P_1,P_8)}$
with notation \ref{nota-system}.

\begin{lemm}
\label{lemm-e8p2}
With the above notation, assume that Conjecture \ref{main_conj}
holds for $P_1$,$P_8$, and any poset $\pos'$ with
$D_0(\pos') \varsubsetneq E_8$. Then Conjecture \ref{main_conj}
holds for $\pos$.
\end{lemm}
\begin{proo}
Let us define $\l_6=\scal{(\a_8,1)}$ and $\mu_4 = \scal{(\a_1,1)}$
 and define $\s^6 = \s^{\l_6}$ and $\tau^4 = \sigma^{\mu_4}$ which are
 classes of degree 6 and 4 respectively. We also consider $\tau^6$
 which corresponds to the ideal $\mu_6=\scal{(\a_2,2)}$. Let
 $\{\g^1,\g^3,\g^4,\g^5,\g^6,\g^7\}$ be a set of generators of
 $H^*(E_8/P_2)$, with $\deg(\g^i)=i$.

Let $u,v,w \in W$ correspond to ideals $\l,\mu,\nu \in I(\pos)$.
Since Conjecture \ref{main_conj} holds if $D_0(w) \varsubsetneq E_8$,
we may have $\notoknu{\l}{\mu}{\nu}$ only if $\nu \supset \pos_0$.

By Proposition \ref{prop-chevalley} $\g^1$ is a good generator.
For $\g^3$ we therefore have $\ok{\g^3}{\s}$ if
$\deg(\s) \leq 10$. In particular $\ok{\g^3}{\tau^4}$,
$\ok{\g^3}{\s^6}$ and $\ok{\g^3}{\tau^6}$.
By recursion we deduce that $\ok{\g^3}{\s}$
if the ideal $\lt$ associated to $\s$ satisfies $\lt\supset\lt_6$, 
$\lt \supset \mu_4$, or $\lt\supset \mu_6$.
If not, 
then $\deg(\s) \leq 9$. Thus $\g^3$ is a good
generator.
The same argument works for $\g^4$.

For $\g^5$ the same argument says that we have $\ok{\g^5}{\s^\l}$
except possibly for the degree 9 ideal $\l = \scal{(\alpha_6,2)}$.

But then we can use Lemma \ref{tous-sauf-1}.
Thus $\g^5$ is a good generator. For $\g^6$ the same argument also works because
there is also only one element of degree 8 not bigger than
$\s^6,\tau^6,\tau^4$, namely $\scal{(\a_5,2),(\a_7,1)}$.

For $\g^7$ we observe that we have already shown that
$\ok{\g^7}{\s}$ for $\s$ of degree at most 6 or $\s \geq \s^6$ or
$\s \geq \tau^4$.
Since $H^7(\pos,\Q)$ is generated as a $\Q$-vector space by 
$h \cdot H^6(\pos,\Q)$ and the elements $\s$ greater than $\s^6$ or $\tau^4$,
$\g^7$ is a good generator.
\end{proo}

\vskip .8cm

We finally deal with $E_8/P_8$.
Thus we consider the system of $\varpi_8$-minuscule
$E_8$-colored posets $\pos_0$ given by the seven
following maximal elements and their obvious intersections:

$$
\begin{array}{cccc}
%auto-ignore

\psset{unit=5mm}

\begin{pspicture*}(-1.5,-0.5)(7.5,15.5)

\psline(1.12,14.88)(1.88,14.12)
\psline(2,13.84)(2,13.16)
\psline(2.12,13.88)(2.88,13.12)
\psline(2,12.84)(2,12.16)
\psline(2.88,12.88)(2.12,12.12)
\psline(3.12,12.88)(3.88,12.12)
\psline(1.88,11.88)(1.12,11.12)
\psline(2.12,11.88)(2.88,11.12)
\psline(3.88,11.88)(3.12,11.12)
\psline(4.12,11.88)(4.88,11.12)
\psline(1.12,10.88)(1.88,10.12)
\psline(2.88,10.88)(2.12,10.12)
\psline(3.12,10.88)(3.88,10.12)
\psline(4.88,10.88)(4.12,10.12)
\psline(5.12,10.88)(5.88,10.12)
\psline(2,9.84)(2,9.16)
\psline(2.12,9.88)(2.88,9.12)
\psline(3.88,9.88)(3.12,9.12)
\psline(4.12,9.88)(4.88,9.12)
\psline(5.88,9.88)(5.12,9.12)
\psline(2,8.84)(2,8.16)
\psline(2.88,8.88)(2.12,8.12)
\psline(3.12,8.88)(3.88,8.12)
\psline(4.88,8.88)(4.12,8.12)
\psline(1.88,7.88)(1.12,7.12)
\psline(2.12,7.88)(2.88,7.12)
\psline(3.88,7.88)(3.12,7.12)
\psline(0.88,6.88)(0.12,6.12)
\psline(1.12,6.88)(1.88,6.12)
\psline(2.88,6.88)(2.12,6.12)
\psline(0.12,5.88)(0.88,5.12)
\psline(1.88,5.88)(1.12,5.12)
\psline(2,5.84)(2,5.16)
\psline(1.12,4.88)(1.88,4.12)
\psline(2,4.84)(2,4.16)
\psline(2.12,3.88)(2.88,3.12)
\psline(3.12,2.88)(3.88,2.12)
\psline(4.12,1.88)(4.88,1.12)
\psline(5.12,0.88)(5.88,0.12)

\psdiamond[fillstyle=solid,fillcolor=black](0,6)(0.256,0.256)
\put(-1.008,6.03){$\tau^7$} 
\pscircle(1,15){0.16}
\pscircle(1,11){0.16}
\pscircle(1,7){0.16}  \psline(0.88,6.88)(1.12,7.12)
\psline(0.88,7.12)(1.12,6.88)   %\put(-0.408,7.03){$\sigma'_{10}$} 
\pscircle*(1,5){0.16}
\pscircle(2,14){0.16}
\pscircle(2,13){0.16}
\pscircle(2,12){0.16}
\pscircle(2,10){0.16}
\pscircle(2,9){0.16}   \put(0.592,9.03){$\sigma^{13}$}
\pscircle(2,8){0.16}
\pscircle(2,6){0.16}  \psline(1.88,5.88)(2.12,6.12)
\psline(1.88,6.12)(2.12,5.88) 
\psdiamond[fillstyle=solid,fillcolor=black](2,5)(0.256,0.256)
\put(2.368,4.83){$\sigma^{6}$} 
\pscircle*(2,4){0.16}
\pscircle(3,13){0.16}
\pscircle(3,11){0.16}
\pscircle(3,9){0.16}   %\put(3.288,8.83){$\sigma_{14}$}
\pscircle(3,7){0.16}  \psline(2.88,6.88)(3.12,7.12)
\psline(2.88,7.12)(3.12,6.88) 
\pscircle*(3,3){0.16}
\pscircle(4,12){0.16}
\pscircle(4,10){0.16}
\pscircle(4,8){0.16}  \psline(3.88,7.88)(4.12,8.12)
\psline(3.88,8.12)(4.12,7.88)   %\put(4.368,7.83){$\sigma_{10}$} 
\pscircle*(4,2){0.16}
\pscircle(5,11){0.16}  
%\psline(4.88,10.88)(5.12,11.12)
%\psline(4.88,11.12)(5.12,10.88) 
\pscircle(5,9){0.16}  \psline(4.88,8.88)(5.12,9.12)
\psline(4.88,9.12)(5.12,8.88)   %\put(5.368,8.83){$\sigma_{11}$} 
\pscircle*(5,1){0.16}
\pscircle(6,10){0.16}  \psline(5.88,9.88)(6.12,10.12)
\psline(5.88,10.12)(6.12,9.88)   \put(6.368,9.83){$\sigma^{12}$} 
\pscircle*(6,0){0.16}

\end{pspicture*} & %auto-ignore

\psset{unit=5mm}

\begin{pspicture*}(0.5,0.5)(-6.5,-14.5)

\psline(-0.12,-13.88)(-0.88,-13.12)
\psline(-1.12,-12.88)(-1.88,-12.12)
\psline(-2.12,-11.88)(-2.88,-11.12)
\psline(-3.12,-10.88)(-3.88,-10.12)
\psline(-4,-9.84)(-4,-9.16)
\psline(-4.12,-9.88)(-4.88,-9.12)
\psline(-4,-8.84)(-4,-8.16)
\psline(-4.88,-8.88)(-4.12,-8.12)
\psline(-5.12,-8.88)(-5.88,-8.12)
\psline(-3.88,-7.88)(-3.12,-7.12)
\psline(-4.12,-7.88)(-4.88,-7.12)
\psline(-5.88,-7.88)(-5.12,-7.12)
\psline(-2.88,-6.88)(-2.12,-6.12)
\psline(-3.12,-6.88)(-3.88,-6.12)
\psline(-4.88,-6.88)(-4.12,-6.12)
\psline(-1.88,-5.88)(-1.12,-5.12)
\psline(-2.12,-5.88)(-2.88,-5.12)
\psline(-3.88,-5.88)(-3.12,-5.12)
\psline(-4,-5.84)(-4,-5.16)
\psline(-0.88,-4.88)(-0.12,-4.12)
\psline(-1.12,-4.88)(-1.88,-4.12)
\psline(-2.88,-4.88)(-2.12,-4.12)
\psline(-3.12,-4.88)(-3.88,-4.12)
\psline(-4,-4.84)(-4,-4.16)
\psline(-0.12,-3.88)(-0.88,-3.12)
\psline(-1.88,-3.88)(-1.12,-3.12)
\psline(-2.12,-3.88)(-2.88,-3.12)
\psline(-3.88,-3.88)(-3.12,-3.12)
\psline(-4.12,-3.88)(-4.88,-3.12)
\psline(-1.12,-2.88)(-1.88,-2.12)
\psline(-2.88,-2.88)(-2.12,-2.12)
\psline(-3.12,-2.88)(-3.88,-2.12)
\psline(-4.88,-2.88)(-4.12,-2.12)
\psline(-5.12,-2.88)(-5.88,-2.12)
\psline(-2.12,-1.88)(-2.88,-1.12)
\psline(-3.88,-1.88)(-3.12,-1.12)
\psline(-4.12,-1.88)(-4.88,-1.12)
\psline(-5.88,-1.88)(-5.12,-1.12)
\psline(-3.12,-0.88)(-3.88,-0.12)
\psline(-4.88,-0.88)(-4.12,-0.12)

\pscircle*(-0,-14){0.16}
\pscircle(-0,-4){0.16}  \psline(0.12,-3.88)(-0.12,-4.12)  \psline(0.12,-4.12)(-0.12,-3.88)
\pscircle*(-1,-13){0.16}
\pscircle(-1,-5){0.16}  \psline(-0.88,-4.88)(-1.12,-5.12)  \psline(-0.88,-5.12)(-1.12,-4.88)
\pscircle(-1,-3){0.16}
\pscircle*(-2,-12){0.16}
\pscircle(-2,-6){0.16}  \psline(-1.88,-5.88)(-2.12,-6.12)  \psline(-1.88,-6.12)(-2.12,-5.88)
\pscircle(-2,-4){0.16}
\pscircle(-2,-2){0.16}
\pscircle*(-3,-11){0.16}
\pscircle(-3,-7){0.16}  \psline(-2.88,-6.88)(-3.12,-7.12)  \psline(-2.88,-7.12)(-3.12,-6.88)
\pscircle(-3,-5){0.16}
\pscircle(-3,-3){0.16}
\pscircle(-3,-1){0.16}
\pscircle*(-4,-10){0.16}
\psdiamond[fillstyle=solid,fillcolor=black](-4,-9)(0.256,0.256)
\pscircle(-4,-8){0.16}  \psline(-3.88,-7.88)(-4.12,-8.12)  \psline(-3.88,-8.12)(-4.12,-7.88)
\pscircle(-4,-6){0.16}
\pscircle(-4,-5){0.16}
\pscircle(-4,-4){0.16}
\pscircle(-4,-2){0.16}
\pscircle(-4,-0){0.16}
\pscircle*(-5,-9){0.16}
\pscircle(-5,-7){0.16}  \psline(-4.88,-6.88)(-5.12,-7.12)  \psline(-4.88,-7.12)(-5.12,-6.88)
\pscircle(-5,-3){0.16}
\pscircle(-5,-1){0.16}
\psdiamond[fillstyle=solid,fillcolor=black](-6,-8)(0.256,0.256)
\pscircle(-6,-2){0.16}

\end{pspicture*} &
%auto-ignore

\psset{unit=5mm}

\begin{pspicture*}(0.5,0.5)(-6.5,-14.5)

\psline(-0.12,-13.88)(-0.88,-13.12)
\psline(-1.12,-12.88)(-1.88,-12.12)
\psline(-2.12,-11.88)(-2.88,-11.12)
\psline(-3.12,-10.88)(-3.88,-10.12)
\psline(-4,-9.84)(-4,-9.16)
\psline(-4.12,-9.88)(-4.88,-9.12)
\psline(-4,-8.84)(-4,-8.16)
\psline(-4.88,-8.88)(-4.12,-8.12)
\psline(-5.12,-8.88)(-5.88,-8.12)
\psline(-3.88,-7.88)(-3.12,-7.12)
\psline(-4.12,-7.88)(-4.88,-7.12)
\psline(-5.88,-7.88)(-5.12,-7.12)
\psline(-2.88,-6.88)(-2.12,-6.12)
\psline(-3.12,-6.88)(-3.88,-6.12)
\psline(-4.88,-6.88)(-4.12,-6.12)
\psline(-1.88,-5.88)(-1.12,-5.12)
\psline(-2.12,-5.88)(-2.88,-5.12)
\psline(-3.88,-5.88)(-3.12,-5.12)
\psline(-4,-5.84)(-4,-5.16)
\psline(-0.88,-4.88)(-0.12,-4.12)
\psline(-1.12,-4.88)(-1.88,-4.12)
\psline(-2.88,-4.88)(-2.12,-4.12)
\psline(-3.12,-4.88)(-3.88,-4.12)
\psline(-4,-4.84)(-4,-4.16)
\psline(-0.12,-3.88)(-0.88,-3.12)
\psline(-1.88,-3.88)(-1.12,-3.12)
\psline(-2.12,-3.88)(-2.88,-3.12)
\psline(-3.88,-3.88)(-3.12,-3.12)
\psline(-4.12,-3.88)(-4.88,-3.12)
\psline(-1.12,-2.88)(-1.88,-2.12)
\psline(-2.88,-2.88)(-2.12,-2.12)
\psline(-3.12,-2.88)(-3.88,-2.12)
\psline(-4.88,-2.88)(-4.12,-2.12)
\psline(-5.12,-2.88)(-5.88,-2.12)
\psline(-2.12,-1.88)(-2.88,-1.12)
\psline(-3.88,-1.88)(-3.12,-1.12)
\psline(-4,-1.84)(-4,-1.16)
\psline(-3.12,-0.88)(-3.88,-0.12)
\psline(-4,-0.84)(-4,-0.16)

\pscircle*(-0,-14){0.16}
\pscircle(-0,-4){0.16}  \psline(0.12,-3.88)(-0.12,-4.12)  \psline(0.12,-4.12)(-0.12,-3.88)
\pscircle*(-1,-13){0.16}
\pscircle(-1,-5){0.16}  \psline(-0.88,-4.88)(-1.12,-5.12)  \psline(-0.88,-5.12)(-1.12,-4.88)
\pscircle(-1,-3){0.16}
\pscircle*(-2,-12){0.16}
\pscircle(-2,-6){0.16}  \psline(-1.88,-5.88)(-2.12,-6.12)  \psline(-1.88,-6.12)(-2.12,-5.88)
\pscircle(-2,-4){0.16}
\pscircle(-2,-2){0.16}
\pscircle*(-3,-11){0.16}
\pscircle(-3,-7){0.16}  \psline(-2.88,-6.88)(-3.12,-7.12)  \psline(-2.88,-7.12)(-3.12,-6.88)
\pscircle(-3,-5){0.16}
\pscircle(-3,-3){0.16}
\pscircle(-3,-1){0.16}
\pscircle*(-4,-10){0.16}
\psdiamond[fillstyle=solid,fillcolor=black](-4,-9)(0.256,0.256)
\pscircle(-4,-8){0.16}  \psline(-3.88,-7.88)(-4.12,-8.12)  \psline(-3.88,-8.12)(-4.12,-7.88)
\pscircle(-4,-6){0.16}
\pscircle(-4,-5){0.16}
\pscircle(-4,-4){0.16}
\pscircle(-4,-2){0.16}
\pscircle(-4,-1){0.16}
\pscircle(-4,-0){0.16}
\pscircle*(-5,-9){0.16}
\pscircle(-5,-7){0.16}  \psline(-4.88,-6.88)(-5.12,-7.12)  \psline(-4.88,-7.12)(-5.12,-6.88)
\pscircle(-5,-3){0.16}
\psdiamond[fillstyle=solid,fillcolor=black](-6,-8)(0.256,0.256)
\pscircle(-6,-2){0.16}

\end{pspicture*} & %auto-ignore

\psset{unit=5mm}

\begin{pspicture*}(0.5,0.5)(-6.5,-13.5)

\psline(-0.12,-12.88)(-0.88,-12.12)
\psline(-1.12,-11.88)(-1.88,-11.12)
\psline(-2.12,-10.88)(-2.88,-10.12)
\psline(-3.12,-9.88)(-3.88,-9.12)
\psline(-4,-8.84)(-4,-8.16)
\psline(-4.12,-8.88)(-4.88,-8.12)
\psline(-4,-7.84)(-4,-7.16)
\psline(-4.88,-7.88)(-4.12,-7.12)
\psline(-5.12,-7.88)(-5.88,-7.12)
\psline(-3.88,-6.88)(-3.12,-6.12)
\psline(-4.12,-6.88)(-4.88,-6.12)
\psline(-5.88,-6.88)(-5.12,-6.12)
\psline(-2.88,-5.88)(-2.12,-5.12)
\psline(-3.12,-5.88)(-3.88,-5.12)
\psline(-4.88,-5.88)(-4.12,-5.12)
\psline(-1.88,-4.88)(-1.12,-4.12)
\psline(-2.12,-4.88)(-2.88,-4.12)
\psline(-3.88,-4.88)(-3.12,-4.12)
\psline(-4,-4.84)(-4,-4.16)
\psline(-0.88,-3.88)(-0.12,-3.12)
\psline(-1.12,-3.88)(-1.88,-3.12)
\psline(-2.88,-3.88)(-2.12,-3.12)
\psline(-3.12,-3.88)(-3.88,-3.12)
\psline(-4,-3.84)(-4,-3.16)
\psline(-0.12,-2.88)(-0.88,-2.12)
\psline(-1.88,-2.88)(-1.12,-2.12)
\psline(-2.12,-2.88)(-2.88,-2.12)
\psline(-3.88,-2.88)(-3.12,-2.12)
\psline(-4.12,-2.88)(-4.88,-2.12)
\psline(-1.12,-1.88)(-1.88,-1.12)
\psline(-2.88,-1.88)(-2.12,-1.12)
\psline(-3.12,-1.88)(-3.88,-1.12)
\psline(-4.88,-1.88)(-4.12,-1.12)
\psline(-5.12,-1.88)(-5.88,-1.12)
\psline(-2.12,-0.88)(-2.88,-0.12)
\psline(-3.88,-0.88)(-3.12,-0.12)
\psline(-4,-0.84)(-4,-0.16)
\psline(-4.12,-0.88)(-4.88,-0.12)
\psline(-5.88,-0.88)(-5.12,-0.12)

\pscircle*(-0,-13){0.16}
\pscircle(-0,-3){0.16}  \psline(0.12,-2.88)(-0.12,-3.12)  \psline(0.12,-3.12)(-0.12,-2.88)
\pscircle*(-1,-12){0.16}
\pscircle(-1,-4){0.16}  \psline(-0.88,-3.88)(-1.12,-4.12)  \psline(-0.88,-4.12)(-1.12,-3.88)
\pscircle(-1,-2){0.16}
\pscircle*(-2,-11){0.16}
\pscircle(-2,-5){0.16}  \psline(-1.88,-4.88)(-2.12,-5.12)  \psline(-1.88,-5.12)(-2.12,-4.88)
\pscircle(-2,-3){0.16}
\pscircle(-2,-1){0.16}
\pscircle*(-3,-10){0.16}
\pscircle(-3,-6){0.16}  \psline(-2.88,-5.88)(-3.12,-6.12)  \psline(-2.88,-6.12)(-3.12,-5.88)
\pscircle(-3,-4){0.16}
\pscircle(-3,-2){0.16}
\pscircle(-3,-0){0.16}
\pscircle*(-4,-9){0.16}
\psdiamond[fillstyle=solid,fillcolor=black](-4,-8)(0.256,0.256)
\pscircle(-4,-7){0.16}  \psline(-3.88,-6.88)(-4.12,-7.12)  \psline(-3.88,-7.12)(-4.12,-6.88)
\pscircle(-4,-5){0.16}
\pscircle(-4,-4){0.16}
\pscircle(-4,-3){0.16}
\pscircle(-4,-1){0.16}
\pscircle(-4,-0){0.16}
\pscircle*(-5,-8){0.16}
\pscircle(-5,-6){0.16}  \psline(-4.88,-5.88)(-5.12,-6.12)  \psline(-4.88,-6.12)(-5.12,-5.88)
\pscircle(-5,-2){0.16}
\pscircle(-5,-0){0.16}
\psdiamond[fillstyle=solid,fillcolor=black](-6,-7)(0.256,0.256)
\pscircle(-6,-1){0.16}

\end{pspicture*} \\
\end{array}$$
$$\begin{array}{cccc}
%auto-ignore

\psset{unit=5mm}

\begin{pspicture*}(0.5,0.5)(-6.5,-14.5)

\psline(-0.12,-13.88)(-0.88,-13.12)
\psline(-1.12,-12.88)(-1.88,-12.12)
\psline(-2.12,-11.88)(-2.88,-11.12)
\psline(-3.12,-10.88)(-3.88,-10.12)
\psline(-4,-9.84)(-4,-9.16)
\psline(-4.12,-9.88)(-4.88,-9.12)
\psline(-4,-8.84)(-4,-8.16)
\psline(-4.88,-8.88)(-4.12,-8.12)
\psline(-5.12,-8.88)(-5.88,-8.12)
\psline(-3.88,-7.88)(-3.12,-7.12)
\psline(-4.12,-7.88)(-4.88,-7.12)
\psline(-5.88,-7.88)(-5.12,-7.12)
\psline(-2.88,-6.88)(-2.12,-6.12)
\psline(-3.12,-6.88)(-3.88,-6.12)
\psline(-4.88,-6.88)(-4.12,-6.12)
\psline(-1.88,-5.88)(-1.12,-5.12)
\psline(-2.12,-5.88)(-2.88,-5.12)
\psline(-3.88,-5.88)(-3.12,-5.12)
\psline(-4,-5.84)(-4,-5.16)
\psline(-0.88,-4.88)(-0.12,-4.12)
\psline(-1.12,-4.88)(-1.88,-4.12)
\psline(-2.88,-4.88)(-2.12,-4.12)
\psline(-3.12,-4.88)(-3.88,-4.12)
\psline(-4,-4.84)(-4,-4.16)
\psline(-0.12,-3.88)(-0.88,-3.12)
\psline(-1.88,-3.88)(-1.12,-3.12)
\psline(-2.12,-3.88)(-2.88,-3.12)
\psline(-3.88,-3.88)(-3.12,-3.12)
\psline(-4.12,-3.88)(-4.88,-3.12)
\psline(-1.12,-2.88)(-1.88,-2.12)
\psline(-2.88,-2.88)(-2.12,-2.12)
\psline(-3.12,-2.88)(-3.88,-2.12)
\psline(-4.88,-2.88)(-4.12,-2.12)
\psline(-5.12,-2.88)(-5.88,-2.12)
\psline(-4,-1.84)(-4,-1.16)
\psline(-4.12,-1.88)(-4.88,-1.12)
\psline(-5.88,-1.88)(-5.12,-1.12)
\psline(-4,-0.84)(-4,-0.16)
\psline(-4.88,-0.88)(-4.12,-0.12)

\pscircle*(-0,-14){0.16}
\pscircle(-0,-4){0.16}  \psline(0.12,-3.88)(-0.12,-4.12)  \psline(0.12,-4.12)(-0.12,-3.88)
\pscircle*(-1,-13){0.16}
\pscircle(-1,-5){0.16}  \psline(-0.88,-4.88)(-1.12,-5.12)  \psline(-0.88,-5.12)(-1.12,-4.88)
\pscircle(-1,-3){0.16}
\pscircle*(-2,-12){0.16}
\pscircle(-2,-6){0.16}  \psline(-1.88,-5.88)(-2.12,-6.12)  \psline(-1.88,-6.12)(-2.12,-5.88)
\pscircle(-2,-4){0.16}
\pscircle(-2,-2){0.16}
\pscircle*(-3,-11){0.16}
\pscircle(-3,-7){0.16}  \psline(-2.88,-6.88)(-3.12,-7.12)  \psline(-2.88,-7.12)(-3.12,-6.88)
\pscircle(-3,-5){0.16}
\pscircle(-3,-3){0.16}
\pscircle*(-4,-10){0.16}
\psdiamond[fillstyle=solid,fillcolor=black](-4,-9)(0.256,0.256)
\pscircle(-4,-8){0.16}  \psline(-3.88,-7.88)(-4.12,-8.12)  \psline(-3.88,-8.12)(-4.12,-7.88)
\pscircle(-4,-6){0.16}
\pscircle(-4,-5){0.16}
\pscircle(-4,-4){0.16}
\pscircle(-4,-2){0.16}
\pscircle(-4,-1){0.16}
\pscircle(-4,-0){0.16}
\pscircle*(-5,-9){0.16}
\pscircle(-5,-7){0.16}  \psline(-4.88,-6.88)(-5.12,-7.12)  \psline(-4.88,-7.12)(-5.12,-6.88)
\pscircle(-5,-3){0.16}
\pscircle(-5,-1){0.16}
\psdiamond[fillstyle=solid,fillcolor=black](-6,-8)(0.256,0.256)
\pscircle(-6,-2){0.16}

\end{pspicture*} & %auto-ignore

\psset{unit=5mm}

\begin{pspicture*}(0.5,0.5)(-6.5,-15.5)

\psline(-0.12,-14.88)(-0.88,-14.12)
\psline(-1.12,-13.88)(-1.88,-13.12)
\psline(-2.12,-12.88)(-2.88,-12.12)
\psline(-3.12,-11.88)(-3.88,-11.12)
\psline(-4,-10.84)(-4,-10.16)
\psline(-4.12,-10.88)(-4.88,-10.12)
\psline(-4,-9.84)(-4,-9.16)
\psline(-4.88,-9.88)(-4.12,-9.12)
\psline(-5.12,-9.88)(-5.88,-9.12)
\psline(-3.88,-8.88)(-3.12,-8.12)
\psline(-4.12,-8.88)(-4.88,-8.12)
\psline(-5.88,-8.88)(-5.12,-8.12)
\psline(-2.88,-7.88)(-2.12,-7.12)
\psline(-3.12,-7.88)(-3.88,-7.12)
\psline(-4.88,-7.88)(-4.12,-7.12)
\psline(-1.88,-6.88)(-1.12,-6.12)
\psline(-2.12,-6.88)(-2.88,-6.12)
\psline(-3.88,-6.88)(-3.12,-6.12)
\psline(-4,-6.84)(-4,-6.16)
\psline(-0.88,-5.88)(-0.12,-5.12)
\psline(-1.12,-5.88)(-1.88,-5.12)
\psline(-2.88,-5.88)(-2.12,-5.12)
\psline(-3.12,-5.88)(-3.88,-5.12)
\psline(-4,-5.84)(-4,-5.16)
\psline(-0.12,-4.88)(-0.88,-4.12)
\psline(-1.88,-4.88)(-1.12,-4.12)
\psline(-2.12,-4.88)(-2.88,-4.12)
\psline(-3.88,-4.88)(-3.12,-4.12)
\psline(-4.12,-4.88)(-4.88,-4.12)
\psline(-3.12,-3.88)(-3.88,-3.12)
\psline(-4.88,-3.88)(-4.12,-3.12)
\psline(-5.12,-3.88)(-5.88,-3.12)
\psline(-4,-2.84)(-4,-2.16)
\psline(-4.12,-2.88)(-4.88,-2.12)
\psline(-5.88,-2.88)(-5.12,-2.12)
\psline(-4,-1.84)(-4,-1.16)
\psline(-4.88,-1.88)(-4.12,-1.12)
\psline(-3.88,-0.88)(-3.12,-0.12)

\pscircle*(-0,-15){0.16}
\pscircle(-0,-5){0.16}  \psline(0.12,-4.88)(-0.12,-5.12)  \psline(0.12,-5.12)(-0.12,-4.88)
\pscircle*(-1,-14){0.16}
\pscircle(-1,-6){0.16}  \psline(-0.88,-5.88)(-1.12,-6.12)  \psline(-0.88,-6.12)(-1.12,-5.88)
\pscircle(-1,-4){0.16}
\pscircle*(-2,-13){0.16}
\pscircle(-2,-7){0.16}  \psline(-1.88,-6.88)(-2.12,-7.12)  \psline(-1.88,-7.12)(-2.12,-6.88)
\pscircle(-2,-5){0.16}
\pscircle*(-3,-12){0.16}
\pscircle(-3,-8){0.16}  \psline(-2.88,-7.88)(-3.12,-8.12)  \psline(-2.88,-8.12)(-3.12,-7.88)
\pscircle(-3,-6){0.16}
\pscircle(-3,-4){0.16}
\pscircle(-3,-0){0.16}
\pscircle*(-4,-11){0.16}
\psdiamond[fillstyle=solid,fillcolor=black](-4,-10)(0.256,0.256)
\pscircle(-4,-9){0.16}  \psline(-3.88,-8.88)(-4.12,-9.12)  \psline(-3.88,-9.12)(-4.12,-8.88)
\pscircle(-4,-7){0.16}
\pscircle(-4,-6){0.16}
\pscircle(-4,-5){0.16}
\pscircle(-4,-3){0.16}
\pscircle(-4,-2){0.16}
\pscircle(-4,-1){0.16}
\pscircle*(-5,-10){0.16}
\pscircle(-5,-8){0.16}  \psline(-4.88,-7.88)(-5.12,-8.12)  \psline(-4.88,-8.12)(-5.12,-7.88)
\pscircle(-5,-4){0.16}
\pscircle(-5,-2){0.16}
\psdiamond[fillstyle=solid,fillcolor=black](-6,-9)(0.256,0.256)
\pscircle(-6,-3){0.16}

\end{pspicture*} &
%auto-ignore

\psset{unit=5mm}

\begin{pspicture*}(0.5,0.5)(-6.5,-16.5)

\psline(-0.12,-15.88)(-0.88,-15.12)
\psline(-1.12,-14.88)(-1.88,-14.12)
\psline(-2.12,-13.88)(-2.88,-13.12)
\psline(-3.12,-12.88)(-3.88,-12.12)
\psline(-4,-11.84)(-4,-11.16)
\psline(-4.12,-11.88)(-4.88,-11.12)
\psline(-4,-10.84)(-4,-10.16)
\psline(-4.88,-10.88)(-4.12,-10.12)
\psline(-5.12,-10.88)(-5.88,-10.12)
\psline(-3.88,-9.88)(-3.12,-9.12)
\psline(-4.12,-9.88)(-4.88,-9.12)
\psline(-5.88,-9.88)(-5.12,-9.12)
\psline(-2.88,-8.88)(-2.12,-8.12)
\psline(-3.12,-8.88)(-3.88,-8.12)
\psline(-4.88,-8.88)(-4.12,-8.12)
\psline(-1.88,-7.88)(-1.12,-7.12)
\psline(-2.12,-7.88)(-2.88,-7.12)
\psline(-3.88,-7.88)(-3.12,-7.12)
\psline(-4,-7.84)(-4,-7.16)
\psline(-0.88,-6.88)(-0.12,-6.12)
\psline(-1.12,-6.88)(-1.88,-6.12)
\psline(-2.88,-6.88)(-2.12,-6.12)
\psline(-3.12,-6.88)(-3.88,-6.12)
\psline(-4,-6.84)(-4,-6.16)
\psline(-2.12,-5.88)(-2.88,-5.12)
\psline(-3.88,-5.88)(-3.12,-5.12)
\psline(-4.12,-5.88)(-4.88,-5.12)
\psline(-3.12,-4.88)(-3.88,-4.12)
\psline(-4.88,-4.88)(-4.12,-4.12)
\psline(-5.12,-4.88)(-5.88,-4.12)
\psline(-4,-3.84)(-4,-3.16)
\psline(-4.12,-3.88)(-4.88,-3.12)
\psline(-5.88,-3.88)(-5.12,-3.12)
\psline(-4,-2.84)(-4,-2.16)
\psline(-4.88,-2.88)(-4.12,-2.12)
\psline(-3.88,-1.88)(-3.12,-1.12)
\psline(-2.88,-0.88)(-2.12,-0.12)

\pscircle*(-0,-16){0.16}
\pscircle(-0,-6){0.16}  \psline(0.12,-5.88)(-0.12,-6.12)  \psline(0.12,-6.12)(-0.12,-5.88)
\pscircle*(-1,-15){0.16}
\pscircle(-1,-7){0.16}  \psline(-0.88,-6.88)(-1.12,-7.12)  \psline(-0.88,-7.12)(-1.12,-6.88)
\pscircle*(-2,-14){0.16}
\pscircle(-2,-8){0.16}  \psline(-1.88,-7.88)(-2.12,-8.12)  \psline(-1.88,-8.12)(-2.12,-7.88)
\pscircle(-2,-6){0.16}
\pscircle(-2,-0){0.16}
\pscircle*(-3,-13){0.16}
\pscircle(-3,-9){0.16}  \psline(-2.88,-8.88)(-3.12,-9.12)  \psline(-2.88,-9.12)(-3.12,-8.88)
\pscircle(-3,-7){0.16}
\pscircle(-3,-5){0.16}
\pscircle(-3,-1){0.16}
\pscircle*(-4,-12){0.16}
\psdiamond[fillstyle=solid,fillcolor=black](-4,-11)(0.256,0.256)
\pscircle(-4,-10){0.16}  \psline(-3.88,-9.88)(-4.12,-10.12)  \psline(-3.88,-10.12)(-4.12,-9.88)
\pscircle(-4,-8){0.16}
\pscircle(-4,-7){0.16}
\pscircle(-4,-6){0.16}
\pscircle(-4,-4){0.16}
\pscircle(-4,-3){0.16}
\pscircle(-4,-2){0.16}
\pscircle*(-5,-11){0.16}
\pscircle(-5,-9){0.16}  \psline(-4.88,-8.88)(-5.12,-9.12)  \psline(-4.88,-9.12)(-5.12,-8.88)
\pscircle(-5,-5){0.16}
\pscircle(-5,-3){0.16}
\psdiamond[fillstyle=solid,fillcolor=black](-6,-10)(0.256,0.256)
\pscircle(-6,-4){0.16}

\end{pspicture*} & \\
\end{array}
$$

We have $S(D_0) = \{1,2\}$. For $i \in \{1,2\}$
let $(D_i,d_i)$ be a marked Dynkin diagram and $P_i$ be any $d_i$-minuscule
$D_i$-colored poset. Set $\pos = \pos_{\pos_0,(P_1,P_2)}$
with notation \ref{nota-system}.

\begin{lemm}
\label{lemm-e8p8}
With the above notation, assume that Conjecture \ref{main_conj}
holds for $P_1$,$P_2$, and any $\lambda$ in $I(\pos)$ with
$D_0(\lambda) \varsubsetneq E_8$. Then Conjecture \ref{main_conj}
holds for $\pos$.
\end{lemm}
\begin{proo}
Set $\l_6 = \scal{(\a_2,1)}$, set $\mu_7 = \scal{(\a_1,1)}$, 
set $\l_{12} = \scal{(\alpha_{8},2)}$ and set $\l_{13} =
\scal{(\a_2,2)}$.
Set $\s^i = \s^{\l_i}$ 
and $\tau^7 = \s^{\mu_7}$.
Let $\{ \g^1,\g^6,\g^{10} \}$ be a set of generators of
$H^*(E_8/P_8)$, with $\deg(\g^i)=i$.
By the hypothesis of the lemma
we know that $\oknu{\l}{\mu}{\nu}$ if $\deg(\nu) \leq 13$.
Let us also give the dimensions of the graded parts of the cohomology:

\vskip 0.2 cm

\centerline{\begin{tabular}{c|lllllllllllllllllllllllllllll}
\hline
$d$&0&1&2&3&4&5&6&7&8&9&10&11&12&13&14\\
\hline
$\dim\ H^d(E_8/P_8)$&1&1&1&1&1&1&2&2&2&2&3&3&4&4&4\\
\hline
\hline
$d$&15&16&17&18&19&20&21&22 
&23&24&25&26&27&28\\ 
\hline
$\dim\ H^d(E_8/P_8)$&4&5&5&6&6&6&6&7&7&7&7&7&7&8\\
\hline
\end{tabular}}

\vskip 0.2 cm

By Proposition \ref{prop-chevalley} we know that $\g^1$ is a good generator.
For $\g^6$ we have $\ok{\g^6}{\s}$ if $\s \leq 7$.
In particular we get $\ok{\g^6}{\tau^7}$ and
$\ok{\g^6}{\s^6}$. Let $\mu,\nu$ such that
$\notoknu{\g^6}{\mu}{\nu}$.
By recursion with respect to $\s^6$ and $\tau^7$ we deduce
that $\ok{\g^6}{\s^\l}$
if $\l \in I(\pos)-I(\pos_0)$
Thus $\mu \in I(\pos_0)$.

Assume that
$\deg(\mu) \leq 13$. By recursion with respect to
$\tau^7$ it follows that if $(\a_1,1) \in \mu$ then $(\a_1,2) \in \nu$.
But $\nu$ must also contain $(\a_8,2)$, thus $\deg(\nu) \geq 20$ and
this contradicts $\deg(\mu) \leq 13$. Thus $(\a_1,1) \not \in \mu$.
Since there is exactly one possible $\mu$ with $7\leq \deg(\nu) \leq 12$
and none with $\deg(\mu) > 12$, 
by Lemma \ref{tous-sauf-1} it follows that $\ok{\g^6}{\s}$ if
$\deg(\s) \leq 13$. 
By Lemma \ref{lefs} it follows that $\ok{\g^6}{\s}$ if $\deg(\s) \leq 15$.

Let us assume that $\deg(\mu)=16$ and $(\a_2,2) \in \mu$.
By recursion with respect to $\s^{13}$ we deduce that
$(\a_2,3) \in \nu$. Since $(\a_1,2) \in \nu$ also it follows that
$\deg(\nu) \geq 23$, and we get a contradiction. Since moreover there
is only one class $\mu$ of degree 16 such that $(\a_8,2) \not \in \mu$
and $(\a_2,2) \not \in \mu$, we conclude that $\deg(\mu) > 16$ by
Lemma \ref{tous-sauf-1}.
By Lemma \ref{lefs} it follows that $\ok{\g^6}{\s}$ if $\deg(\s) \leq 17$.

There are three classes $\s^\l$ of degree 17 resp. 18 such that
$(\a_8,2) \not \in \l$, namely
$\scal{(\a_2,2),(\a_6,3)}$, $\scal{(\a_4,4),(\a_7,2)}$, $\scal{(\a_3,2)}$
resp.
$\scal{(\a_1,2)}$, $\scal{(\a_3,3)$, $(\a_7,2)},\scal{(\a_4,4),(\a_6,3)}$,
and the corresponding map given by multiplication
by $h$ is surjective; thus we conclude thanks to Lemma 
\ref{lefs}$(\imath\imath)$
that $\ok{\g^6}{\s}$ if $\deg(\s) = 18$. Then Lemma \ref{lefs}$(\imath)$
gives the same identity for $\deg(\s) \leq 21$.

We finish showing that $\g^6$ is a good generator thanks again to
Lemma \ref{lefs}$(\imath\imath)$, because there are exactly two classes 
in each degree
21 and 22 which are not bigger than $\s^{12}$.

\vskip .5 cm

We now consider $\g^{10}$. By Lemma \ref{lemm-sous-algebre}
we have already proved that
$\ok{\g^{10}}{\s}$ for $\deg(\s) \leq 9$. 
For $\deg(\s) = 10$, 
we must consider some degrees and we shall assume $\g^{10}=\s^{10}$.
Remark first that
all classes $(\tau^{20,i})_{i\in[1,6]}$ of degree 20 in $H^*(\pos_0)$
do not contain one of the vertices $(\a_1,2)$ or $(\a_8,2)$
except for $\tau^{20,4}=\scal{(\a_1,2),(\a_8,2)}$. In particular, we
have the equalities $\oknu{\g^{10}}{\s}{\s'}$ for all
degree 10 classes $\s$ and all degree 20 classes $\s'\neq\tau^{20,4}$.
Let us be more precise here, define the following ideals
$$\begin{array}{lll}
\lt_{20,1}=\scal{(\a_4,4),(\a_7,3)} &
\lt_{20,2}=\scal{(\a_5,4),(\a_8,2)} & 
\lt_{20,3}=\scal{(\a_3,3),(\a_6,3),(\a_8,2)} \\  
\lt_{20,4}=\scal{(\a_1,2),(\a_8,2)} &
\lt_{20,5}=\scal{(\a_1,2),(\a_6,3)}&
\lt_{20,6}=\scal{(\a_3,3),(\a_5,4)}\\
\end{array}$$
and the cohomology classes $\tau^{20,i}=\s^{\lt_{20,i}}$ for
$i\in[1,6]$.

As we have seen, the algebras $H^*(\pos)$ and
$H^*_t(X)$ coincide in degree 20 except maybe for that class
$\tau^{20,4}$. Using jeu de taquin, we have the equality
$$\begin{array}{l}
\s^{10} \pp \s^{10} = 16 \tau^{20,1} + 8 \tau ^{20,2} + 14 \tau^{20,3} + 
7 \tau^{20,4} + 4 \tau^{20,5} + 2 \tau^{20,6}.
\end{array}$$
and because of the coincidence of the two algebras we get
$$\begin{array}{l}
\s^{10} \cdot \s^{10} = 16 \tau^{20,1} + 8 \tau^{20,2} + 14 \tau^{20,3} + 
x \tau^{20,4} + 4 \tau^{20,5} + 2 \tau^{20,6}
\end{array}$$
for some non negative integer $x$.
However, we are able to compute the coefficient
$x=c_{\s^{10},\s^{10}}^{\tau^{20,4}}$ using the degree of
these classes. To compute $x$ we use the 
Hasse diagram to obtain the following degrees:
$$\begin{array}{lll}
\deg(\tau^{20,1})=4322859480&\deg(\tau^{20,2})=6717795480& 
\deg(\tau^{20,3})=8298453240\\
\deg(\tau^{20,4})=1560699960&\deg(\tau^{20,5})=3789366840&
\deg(\tau^{20,6})=10269733320\\ 
\deg(\s^{10} \cdot \s^{10})=285708294600.&\\
\end{array}$$
Remark that here we made the computation in the cohomology of the
finite dimensional homogeneous space $E_8/P_8$ and used Poincar{\'e}
duality over this space. The degree is linear and we get the value $x=7$. 
Thus the two algebras also coincide in degre 20.

Let us remark here that these computation where made with the help of a 
computer. It is a quite easy computation for the Hasse diagram. For the 
jeu de taquin, we made an (easy) adaptation of the computer program writen 
by H. Thomas and A. Yong for the cominuscule jeu de taquin.

\vskip .2cm

Since any class of degree at most 19 can be expressed as
$P(\g^1,\g^6) + \g^{10} \cdot Q(\g^1,\g^6)$ we have
$\ok{\g^{10}}{\s}$ if $\deg(\s) \leq 19$ by Lemma
\ref{lemm-petite-puissance}. For higher degree classes, we conclude as
for $\g^6$.
\end{proo}

%%% Local Variables: 
%%% mode: latex
%%% TeX-master: t
%%% End: 

\def \ii {\overline \imath}

\section{Non simply-laced case}

\label{section-pliage}

\subsection{General results for the push-forward of a minuscule class}

We will now explain how it is possible to obtain Theorem \ref{main-theo}
in the non simply-laced cases using folding. First we deal with the
minuscule case. More precisely, let $(D_0,d_0),(E,e)$ be
marked Dynkin diagrams, $p$ an integer, and $x_0 \in D_0$.
We consider the disjoint union $\amalg_{1\leq i\leq p} D_0^i$ of
$p$ copies of $D_0$ denoted $D_0^i$ and an
automorphism $\t$ of $\amalg_i D_0^i$
induced by a cyclic permutation of order $p$
of $[1,p]$. In each $D_0^i$ we denote $x_0^i$ the element
corresponding to $x_0$. We consider the Dynkin diagram obtained from the
disjoint union of $E$ and $\amalg_i D_0^i$ connecting each $x_0^i$ with $e$.
We still denote $\t$ the automorphism of $D$ extending $\theta$ by setting
$\t(x)=x$ for $x \in E$. Moreover we denote $d=d_0^1 \in D$.

Thus $D$ defines a Kac-Moody algebra $\fg$, and
$\t$ an automorphism of $\fg$. We denote $\fg^\t$ the subalgebra of
invariant elements, with Dynkin diagram $D^\t$ indexed by the
equivalence classes of elements in $D$ modulo $\t$,
$G^\t$ the corresponding subgroup of $G$, and $W^\t$ the Weyl group
of $D^\t$. For $i \in D$ let $\overline \imath \in D^\t$ denote
its natural projection. Denote $\overline D_0$ resp. $\overline E$
the image of $D_0^1$ resp. $E$ under this projection.
We denote $\overline x_0$ the element
$\overline{x_0^i}$ for any $i\in[1,p]$.

\centerline{\begin{pspicture*}(0,0)(12.0,7.2000003)

\psellipse[fillstyle=solid,fillcolor=black](1.2,1.2)(0.096,0.096)
\psellipse(1.8,5.4)(0.096,0.096)
\psline(1.732,5.332)(1.868,5.468)
\psline(1.732,5.468)(1.868,5.332)
\psellipse(1.8,4.2)(0.096,0.096)
\psline(1.732,4.132)(1.868,4.268)
\psline(1.732,4.268)(1.868,4.132)
\psellipse[fillstyle=solid,fillcolor=black](4.8,1.2)(0.096,0.096)
\psellipse(5.4,4.2)(0.096,0.096)
\psline(5.332,4.132)(5.468,4.268)
\psline(5.332,4.268)(5.468,4.132)
\psellipse(7.8,5.4)(0.096,0.096)
\psline(7.732,5.332)(7.868,5.468)
\psline(7.732,5.468)(7.868,5.332)
\psellipse(7.8,4.2)(0.096,0.096)
\psline(7.732,4.132)(7.868,4.268)
\psline(7.732,4.268)(7.868,4.132)
\psellipse(10.2,4.2)(0.096,0.096)
\psline(10.132,4.132)(10.268,4.268)
\psline(10.132,4.268)(10.268,4.132)
\psline(1.692,5.4)(1.692,4.2)
\psline(1.764,5.311)(1.764,4.289)
\psline(1.836,5.311)(1.836,4.289)
\psline(1.908,5.4)(1.908,4.2)
%\psline(1.8,4.74)(1.877,4.832)
\psline(1.8,4.74)(2,5)
\psline(1.8,4.74)(1.6,5)
\psline(7.714,5.357)(5.486,4.243)
\psline(7.8,5.304)(7.8,4.296)
\psline(7.886,5.357)(10.114,4.243)
\uput[u](1.8,1.8){$\overline{D}_0$}
\uput[u](5.4,1.8){$D_0^1$}
\uput[u](7.8,1.8){$D_0^2$}
\uput[u](10.2,1.8){$D_0^p$}
\uput[l](1.8,4.8){$p$}
\uput[l](1.8,4.2){$\overline{x}_0$}
\uput[l](5.4,4.2){${x}_0^1$}
\uput[l](7.8,4.2){${x}_0^2$}
\uput[r](10.2,4.2){${x}_0^p$}
\put(1.65,6.0){$\overline{E}$}
\put(7.65,6.0){$E$}
\put(7.2,5.4){$e$}
\put(1.2,5.4){$\overline{e}$}
\put(9.0,4.2){$\ldots$}
\uput[dl](1.2,1.2){$\overline{d}$}
\uput[dl](4.8,1.2){${d}$}
\psbezier[linewidth=2pt](1.8,4.2)(-1.1279981,-0.21600123)(4.7280016,-0.21599893)(1.8,4.2)
\psbezier[linewidth=2pt](5.4,4.2)(2.472002,-0.21600123)(8.328002,-0.21599893)(5.4,4.2)
\psbezier[linewidth=2pt](7.8,4.2)(4.872002,-0.21600123)(10.728003,-0.21599893)(7.8,4.2)
\psbezier[linewidth=2pt](10.2,4.2)(7.2720017,-0.21600123)(13.128002,-0.21599893)(10.2,4.2)
\psbezier[linewidth=2pt](1.8,5.4)(4.3920007,7.031999)(-0.79199934,7.032001)(1.8,5.4)
\psbezier[linewidth=2pt](7.8,5.4)(10.392,7.031999)(5.2080007,7.032001)(7.8,5.4)

\end{pspicture*}}

Let $P$ resp. $P^\t$ be the parabolic subgroup of $G$
resp. $G^\t$ corresponding to $d$ resp. $\overline d$; we have injections
$i : W^\t \to W$ and $\iota : G^\t/P^\t \to G/P$. Denoting with $t_m$ the 
simple reflections in $W^\t$ and with $s_j$ the simple reflections in $W$, 
note that we have $i(t_m) = \prod_{j:\overline \jmath=m} s_j \in W$.
The idea to prove Conjecture \ref{main_conj} in this situation is to use 
the fact that $\iota^* : H^*(G/P) \to H^*(G^\t/P^\t)$
and $\iota_* : H^*(G^\t/P^\t) \to H^*(G/P)$ are adjoint and to 
to compare Littlewood-Richardson coefficients on $G/P$ with those
on $G^\t/P^\t$. For this it is usefull to show that minuscule Schubert
cells are mapped to minuscule Schubert cells by $\iota_*$.

We first show that if $p \geq 3$ then the situation is quite simple
because there are very few $\overline d$-minuscule elements.

\begin{lemm}
\label{lemm-n3}
If $p \geq 3$ then any $\overline d$-minuscule element is either
in $W(\overline D_0)$ or can be written
as $vu$ with $u\in W(\overline D_0)$ a
$\overline d$-minuscule element and $v \in W(\overline E)$
an $\overline e$-minuscule element.
\end{lemm}
\begin{proo}
If the reflexion with respect to $\overline e$
does not appear in a reduced expression of $w$, then clearly
$w$ belongs to $W(\overline D_0)$. Assuming now that
there exists an integer $k$ such that $m_k = \overline e$, since we have
$\scal{t_{m_{k+1}} \cdots t_{m_l} (\overline \Lambda) ,
\beta_{\overline x_0}^\vee} \geq -1$, we deduce
$\scal{t_{m_{k}} \cdots t_{m_l} (\overline \Lambda) ,
\beta_{\overline x_0}^\vee} \geq p-1 \geq 2$, so that
for all $k'\leq k$ we have
$\scal{t_{m_{k'}} \cdots t_{m_l} (\overline \Lambda) ,
\beta_{\overline x_0}^\vee} \geq 2$ and
$m_{k'} \not = \overline x_0$. Therefore up to using some
commutation relations we may write $w$ as $vu$ with $v\in W(\overline{E})$ 
and $u\in W(\overline{D}_0)$. Since any reduced expression of $w$ satisfies 
the conditions of Definition \ref{defi-minuscule}, $v$ is 
$\overline{e}$-minuscule and $u$ is $\overline{d}$-minuscule.
\end{proo}

\begin{lemm}
\label{lemm-folding-schubert}
Let $w \in W^\t$ be $\overline d$-minuscule. Then the class of $i(w)$
in $W/W_P$ can be represented by a unique $d$-minuscule element $u$.
This element satisfies $l(u) = l(w)$ and we have the equality
$\iota(\overline{B^\t wP^\t/P^\t}) = \overline{BuP/P}$.
\end{lemm}
\begin{proo}
Let $\Lambda$ resp. $\overline \Lambda$
be the weight corresponding to $d$ resp. $\overline d$. Then
$\overline \Lambda$ is the restriction of $\Lambda$ to $\fh^\t$. Let
$\a_j, j\in D$ resp. $\beta_m,m\in D^\t$ denote the simple roots
of $G$ resp. $G^\t$.
Let us denote $t_m \in W^\t$ the reflexion corresponding to $m \in D^\t$.
Let $w \in W^\t$ be $\overline \Lambda$-minuscule and
let $w = t_{m_1} \cdots t_{m_l}$ be a reduced decomposition of $w$.
Since $w$ is $\overline \Lambda$-minuscule, we have
$\scal{t_{m_{2}} \cdots t_{m_l} (\overline \Lambda) , \beta^\vee_{m_1}} = 1$.
We have
$$
\beta^\vee_{m_1} = \sum_{j:\overline \jmath = m_1} \a_j^\vee,
$$
thus we get that
$\sum_j \scal{\iota(t_{m_{2}} \cdots t_{m_l}) (\Lambda) , \a_j^\vee } = 1$.

\vskip .4cm

We claim that if $\overline \jmath = m_1$ then
$\scal{\iota(t_{m_{2}} \cdots t_{m_l}) (\Lambda) , \a_j^\vee }$ is 
nonnegative. In case $p>2$ the claim is easily verified using Lemma 
\ref{lemm-n3}.

Let us now assume that $p=2$ and let us
choose $j$ with $\overline \jmath = m_1$ and 
$\scal{\iota(t_{m_{2}} \cdots t_{m_l}) (\Lambda) , \a_j^\vee } > 0$.
If $j$ is the unique element $k$ such that $\overline k = m_1$,
then the claim is true, so we can assume that $\t(j) \not = j$, so that
$\{j,\t(j)\} = \{k : \overline k = m_1 \}$.
For $w \in W$ let $l_P(w)$ denote the length of its minimal length
representative in $W/W_P$. 
Since $Bi(t_{m_{1}} \cdots t_{m_l})P/P$ contains
$\iota(B^\t t_{m_{1}} \cdots t_{m_l} P^\t/P^\t)$, of dimension 
$l=l_{P^\t}(w)$, we have
$l_P(i(t_{m_{1}} \cdots t_{m_l})) \geq l$.
By induction we may assume that 
$l_P(i(t_{m_{2}} \cdots t_{m_l})) = l-1$.
If $\scal{\iota(t_{m_{2}} \cdots t_{m_l}) (\Lambda) , \a_{\t(j)}^\vee } < 0$,
then we would have
$l_P(s_{\t(j)} \cdot \iota(t_{m_{2}} \cdots t_{m_l})) < l-1$ and thus
$l_P(\iota(t_{m_{1}} \cdots t_{m_l})) \leq l-1$. We have already seen that
this does not occur.

\vskip .4cm

Thus the claim is proved
and the class of $\iota(w)$ in $W/W_P$ is equal to the class of
$s_j \cdot t_{m_2} \cdots t_{m_l}$, and thus $l_P(\iota(w)) = l$.
Moreover if $u_2$ is a $d$-minuscule element which represents the
class of $i(t_{m_2} \cdots t_{m_l})$ in $W/W_P$, then
$s_j \cdot u_2$ represents $i(w)$. Finally,
$\iota$ restricts to an inclusion
$B^\t w P^\t/P^\t \to BuP/P$ of $l$-dimensional irreducible varieties,
so we have the equality $\iota(\overline{B^\t w P^\t/P^\t}) = 
\overline{BuP/P}$.
\end{proo}

\begin{nota}
Let $w \in W^\t$ be $\overline d$-minuscule. We denote $\overline \imath(w)$
the unique $d$-minuscule element in $W$ which has the same class as
$i(w)$ modulo $P$. Such an element exists by Lemma \ref{lemm-folding-schubert}.

For $w\in W^\t$ resp. $v \in W$ let $\s_w,\s^w$ resp. $\tau_v,\tau^v$
denote the corresponding homology and cohomology classes.
\end{nota}

\begin{lemm}
\label{lemm-iota*}

(\i) Let $w\in W^\t$ be $\overline d$-minuscule. Then
$\iota_* \s_w = \tau_{\overline \imath(w)}$.

(\i\i) Let $w \in W^\t$ be $\overline d$-minuscule
and assume that all degree $d$
classes in $G^\t/P^\t$ are $d$-minuscule. Then
$\iota^* \tau^{\overline \imath(w)} = \s^w$.

\end{lemm}
\begin{proo}
Point $(\imath)$ follows directly from Lemma \ref{lemm-folding-schubert}.
Let $w \in W^\t$ be as in $(\imath\imath)$.
Then for $w' \in W^\t$ a $d$-minuscule element one computes that
$$
\scal { \iota^* \tau^{\overline \imath(w)} , \s_{w'} } =
\scal { \tau^{\overline \imath(w)} , \iota_* \s_{w'} } =
\scal { \tau^{\overline \imath(w)} , \tau_{\overline \imath(w')} } =
\delta_{w',w} = \scal { \s^w , \s_{w'} },
$$
thus the lemma is proved.
\end{proo}

\subsection{Type $B_n$ case}

In this case, we consider  the system of $\varpi_n$-minuscule 
$B_n$-colored posets $\pos_0$ given by the poset of an isotropic
Grassmannian $\G_Q(n,2n+1)$.
We have $S_0=\{1\}$. Let $(D_1,d_1)$ be a
marked Dynkin diagram and $P_1$ be any $d_1$-minuscule $D_1$-colored
poset. Set $\pos=\pos_{\pos_0,\{P_{1}\}}$. 
We apply the above construction with $D_0$ reduced to one vertex $d_0=d$,
$E$ obtained from a union $D_1 \cup A$, where $A$ is of type $A_{n-1}$, 
attaching $d_1$ on the first node of $A$, $e$ the last element of $A$,
and $p=2$. Thus $D$ resp. $D^\t$ is obtained as a union of $D_1$ and a Dynkin
diagram of type $D_{n+1}$ resp. $B_n$. Moreover we see that
the heaps of $w$ and $\overline \imath(w)$ are isomorphic
for any $w$ corresponding to an ideal in $\pos$
(although they are not isomorphic as colored heaps).

\centerline{\begin{pspicture*}(0,0)(12.0,3.6000001)

\psellipse(2.4,1.8)(0.096,0.096)
\psellipse(3.0,1.8)(0.096,0.096)
\psellipse(4.2,1.8)(0.096,0.096)
\psellipse[fillstyle=solid,fillcolor=black](4.8,1.8)(0.096,0.096)
\psellipse(9.0,1.8)(0.096,0.096)
\psellipse(9.6,1.8)(0.096,0.096)
\psellipse(10.8,1.8)(0.096,0.096)
\psellipse(11.4,2.4)(0.096,0.096)
\psellipse[fillstyle=solid,fillcolor=black](11.4,1.2)(0.096,0.096)
\psline(2.304,1.8)(2.304,1.8)
\psline(2.496,1.8)(2.904,1.8)
\psline(4.289,1.764)(4.711,1.764)
\psline(4.289,1.836)(4.711,1.836)
\psline(4.55,1.8)(4.458,1.877)
\psline(4.55,1.8)(4.458,1.723)
\psline(9.096,1.8)(9.504,1.8)
\psline(10.868,1.868)(11.332,2.332)
\psline(10.868,1.732)(11.332,1.268)
\psline[linestyle=dashed](3.096,1.8)(4.104,1.8)
\psline[linestyle=dashed](9.696,1.8)(10.704,1.8)
\put(1.2,1.8){$D_1$}
\put(7.8,1.8){$D_1$}
\uput[u](9.0,1.8){1}
\uput[u](2.4,1.8){1}
\uput[u](4.2,1.8){$n-1$}
\uput[u](10.8,1.8){$n-1$}
\uput[d](4.8,1.8){$\overline{d}$}
\uput[d](11.4,1.2){$d$}
\psbezier[linewidth=2pt](9.0,1.8)(6.6,4.2000003)(6.6,-0.5999999)(9.0,1.8)
\psbezier[linewidth=2pt](2.4,1.8)(0,4.2000003)(0,-0.5999999)(2.4,1.8)

\end{pspicture*}}

\begin{lemm}
\label{lemm-bnpn}
With the above notation, assume that Conjecture \ref{main_conj} holds
for $P_{1}$. Then Conjecture \ref{main_conj} 
holds for $\pos$.
\end{lemm}
\begin{proo}
Let $\g^1,\ldots,\g^n$ be a set of generators of $H^*(\pos_0)$ with
$\deg(\g^i)=i$. By Lemma \ref{oreille} it is enough to show that
$\ok{\g^i}{\s}$ for any $\s \in H^*(\pos_0)$.

Let $u_i \in W^\t$ denote the element corresponding to $\gamma^i$.
It is enough to show that for any elements $v,w$ in $W^\t$
we have $t_{u_i,v}^w = c_{u_i,v}^w$.
We compute $c_{u_i,v}^w$ as the coefficient of
$\s_v$ in $\s^{u_i} \cap \s_w$.
Since $\deg(\g^i) = i \leq n$, all classes of degree $i$
correspond to ideals in
$\pos_0$ and thus are minuscule. So by Lemma \ref{lemm-iota*} we deduce
that $\iota^* \tau^{\overline \imath(u_i)} = \s^{u_i}$.
Thus by Lemma \ref{lemm-iota*} again we get
$$
\iota_* ( \s^{u_i} \cap \s_w ) =
\iota_* ( \iota^* \tau^{\overline \imath(u_i)} \cap \s_w ) =
\tau^{\overline \imath(u_i)} \cap \iota_*(\s_w) =
\tau^{\ii(u_i)} \cap \tau_{\ii(w)}.
$$
Thus the coefficient of $\s_v$ in $\s^{u_i} \cap \s_w$
is the same as the coefficient
of $\tau_{\ii(v)}$ in the cap product $\tau^{\ii(u_i)} \cap \tau_{\ii(w)}$.
In other words $c_{u_i,v}^w = c_{\ii(u_i),\ii(v)}^{\ii(w)}$. Now by
Lemma \ref{lemm-dnpn} we know that the latter equals
$t_{\ii(u_i),\ii(v)}^{\ii(w)}$.
Since the heaps of $w$ and $\ii(w)$
are isomorphic, we deduce
$t_{\ii(u_i),\ii(v)}^{\ii(w)} = t_{u_i,v}^w$. Therefore
$c_{u_i,v}^{w} = t_{u_i,v}^w$, which is exactly what we wanted to prove.
\end{proo}

\subsection{Type $F_4$ minuscule case}

In this case, we consider  the system of $\varpi_4$-minuscule 
$F_4$-colored posets $\pos_0$ given by the following picture:

\centerline{\begin{pspicture*}(0,0)(3.0,3.6000001)

\pspolygon[fillstyle=solid,fillcolor=black](0.72,2.4)(0.6,2.52)(0.48,2.4)(0.6,2.28)
\psellipse(1.2,3.0)(0.096,0.096)
\psline(1.132,2.932)(1.268,3.068)
\psline(1.132,3.068)(1.268,2.932)
\psellipse[fillstyle=solid,fillcolor=black](1.2,1.8)(0.096,0.096)
\psellipse(1.8,2.4)(0.096,0.096)
\psline(1.732,2.332)(1.868,2.468)
\psline(1.732,2.468)(1.868,2.332)
\psellipse[fillstyle=solid,fillcolor=black](1.8,1.2)(0.096,0.096)
\psellipse(2.4,3.0)(0.096,0.096)
\psline(2.332,2.932)(2.468,3.068)
\psline(2.332,3.068)(2.468,2.932)
\psellipse[fillstyle=solid,fillcolor=black](2.4,0.6)(0.096,0.096)
\psline(2.332,0.668)(1.868,1.132)
\psline(1.237,1.712)(1.712,1.237)
\psline(1.288,1.763)(1.763,1.288)
\psline(1.55,1.45)(1.54,1.57)
\psline(1.55,1.45)(1.43,1.46)
\psline(1.288,1.837)(1.763,2.312)
\psline(1.237,1.888)(1.712,2.363)
\psline(1.55,2.15)(1.43,2.14)
\psline(1.55,2.15)(1.54,2.03)
\psline(1.237,2.912)(1.712,2.437)
\psline(1.288,2.963)(1.763,2.488)
\psline(1.55,2.65)(1.54,2.77)
\psline(1.55,2.65)(1.43,2.66)
\psline(1.868,2.468)(2.332,2.932)
\psline(0.66,2.34)(1.132,1.868)
\psline(0.66,2.46)(1.132,2.932)

\end{pspicture*}}

We have $S_0=\{1\}$. Let $(D_1,d_1)$ be a
marked Dynkin diagram and $P_1$ be any $d_1$-minuscule $D_1$-colored
poset. Set $\pos=\pos_{\pos_0,\{P_{1}\}}$.

\centerline{\begin{pspicture*}(0,0)(11.400001,3.6000001)

\psellipse(3.0,1.8)(0.096,0.096)
\psellipse(3.6,1.8)(0.096,0.096)
\psellipse(4.2,1.8)(0.096,0.096)
\psellipse[fillstyle=solid,fillcolor=black](4.8,1.8)(0.096,0.096)
\psellipse(8.4,1.8)(0.096,0.096)
\psellipse(9.0,1.8)(0.096,0.096)
\psellipse(9.6,2.4)(0.096,0.096)
\psellipse(9.6,1.2)(0.096,0.096)
\psellipse(10.2,2.4)(0.096,0.096)
\psellipse[fillstyle=solid,fillcolor=black](10.2,1.2)(0.096,0.096)
\psline(3.689,1.764)(4.111,1.764)
\psline(3.689,1.836)(4.111,1.836)
\psline(3.95,1.8)(3.858,1.877)
\psline(3.95,1.8)(3.858,1.723)
\psline(3.096,1.8)(3.504,1.8)
\psline(4.296,1.8)(4.704,1.8)
\psline(8.496,1.8)(8.904,1.8)
\psline(9.068,1.868)(9.532,2.332)
\psline(9.696,2.4)(10.104,2.4)
\psline(9.068,1.732)(9.532,1.268)
\psline(9.696,1.2)(10.104,1.2)
\uput[d](10.2,1.2){$d$}
\uput[d](4.8,1.8){$\overline{d}$}
\put(1.8,1.8){$D_1$}
\put(7.2,1.8){$D_1$}
\psbezier[linewidth=2pt](3.0,1.8)(0.5999999,4.2000003)(0.5999999,-0.5999999)(3.0,1.8)
\psbezier[linewidth=2pt](8.4,1.8)(5.9999995,4.2000003)(5.9999995,-0.5999999)(8.4,1.8)

\end{pspicture*}}

\begin{lemm}
\label{lemm-f4p4}
With the above notation, assume that Conjecture \ref{main_conj} holds
for $P_{1}$. Then Conjecture \ref{main_conj} holds for $\pos$.
\end{lemm}
\begin{proo}
Let $\g^1,\g^4$ be a set of generators of $H^*(\pos_0)$ with
$\deg(\g^i)=i$. By Lemma \ref{oreille} it is enough to show that
$\ok{\g^i}{\s}$ for any $\s \in H^*(\pos_0)$. For $\g^1$ this is
already known by Proposition \ref{prop-chevalley}. Moreover
by Lemma \ref{lefs} $(\imath \imath)$ it is enough to show that
$\ok{\g^4}{\g^4}$.

To prove this we consider the above construction with
$D_0$ of type $A_2$ and $d_0$ the first node of $A_2$ and $d$ the last 
node, $E$ obtained as a connected union of $D_1$ and again a Dynkin
diagram of type $A_2$, and $p=2$. Here $D$ resp. $D^\t$ is a connected
union of $D_1$ and a Dynkin diagram of type $E_6$ resp. $F_4$ (cf. the 
above picture).
Again we see that
the heaps of $w$ and $\overline \imath(w)$ are isomorphic
for any $w$ corresponding to an ideal in $\pos$.

The rest of the proof of the lemma is the same as for
Lemma \ref{lemm-bnpn}, using the fact that any class of degree 4
corresponding to an ideal in $\pos$
is minuscule and Lemma \ref{lemm-e6p1}.
\end{proo}

%%% Local Variables: 
%%% mode: latex
%%% TeX-master: t
%%% End: 

\def \bb {{{\underline{\beta}}}}

\subsection{General result for $\Lambda$-cominuscule classes}

Let $(\cD_0,d)$ be a marked Dynkin diagram, let $G_0$ be the associated 
Kac-Moody group and $P_0$ the corresponding parabolic subgroup. 
Let $W_{G_0}$ the Weyl group of $G_0$ and let $w\in W_{G_0}^{P_0}$ 
(the set of minimal length representatives for $P_0$). We shall assume 
that $\cD_0$ is the support of $w$. 
Choose a simple root $\a$ or equivalently a vertex of $\cD_0$ 
(still denoted $\a$)
and a Dynkin diagram $\cD'$ containing $\cD_0$ and one more root $\beta$
only connected to $\a$ in $\cD'$. If $\scal{\beta,\a^\vee}=p$ we also 
define a
Dynkin diagram $\cD$ containing $\cD_0$ and $p$ more vertices
labelled $(\beta_i)_{i\in[1,p]}$ all only connected to $\a$ with a simple 
edge. In the following we depicted $\cD'$ on the left and $\cD$ on
the right.

\centerline{\begin{pspicture*}(0,0)(9.0,3.6000001)

\psellipse[fillstyle=solid,fillcolor=black](1.2,1.8)(0.096,0.096)
\psellipse[fillstyle=solid,fillcolor=black](3.0,1.8)(0.096,0.096)
\psellipse[fillstyle=solid,fillcolor=black](3.6,1.8)(0.096,0.096)
\psellipse[fillstyle=solid,fillcolor=black](5.4,1.8)(0.096,0.096)
\psellipse[fillstyle=solid,fillcolor=black](7.2,1.8)(0.096,0.096)
\psellipse[fillstyle=solid,fillcolor=black](7.8,3.0)(0.096,0.096)
\psellipse[fillstyle=solid,fillcolor=black](7.8,2.4)(0.096,0.096)
\psellipse[fillstyle=solid,fillcolor=black](7.8,1.2)(0.096,0.096)
\psellipse[fillstyle=solid,fillcolor=black](7.8,0.6)(0.096,0.096)
\psline(3.0,1.692)(3.6,1.692)
\psline(3.089,1.764)(3.511,1.764)
\psline(3.089,1.836)(3.511,1.836)
\psline(3.0,1.908)(3.6,1.908)
%\psline(3.33,1.8)(3.238,1.877)
%\psline(3.33,1.8)(3.238,1.723)
\psline(3.36,1.8)(3.21,2)
\psline(3.36,1.8)(3.21,1.6)
%
%\psline(3.36,1.8)(3.268,1.877)
%\psline(3.36,1.8)(3.268,1.723)
%
%\psline(3.39,1.8)(3.298,1.877)
%\psline(3.39,1.8)(3.298,1.723)
\psline(7.243,1.886)(7.757,2.914)
\psline(7.268,1.868)(7.732,2.332)
\psline(7.268,1.732)(7.732,1.268)
\psline(7.243,1.714)(7.757,0.686)
\uput[l](1.2,1.8){$d$}
\uput[l](5.4,1.8){$d$}
\uput[l](7.2,1.8){$\a$}
\uput[l](7.2,0.5){$\cD$}
\uput[l](3.0,1.8){$\a$}
\uput[l](3.0,0.5){$\cD'$}
\uput[dr](3.0,1.6){$p$}
\uput[r](3.6,1.8){$\beta$}
\uput[r](7.8,3.0){$\beta_1$}
\uput[r](7.8,2.4){$\beta_2$}
\uput[r](7.8,0.6){$\beta_p$}
\put(7.8,1.8){$\vdots$}
\psbezier[linewidth=2pt](3.0,1.8)(0.5999999,4.2000003)(0.5999999,-0.5999999)(3.0,1.8)
\psbezier[linewidth=2pt](7.2,1.8)(4.7999997,4.2000003)(4.7999997,-0.5999999)(7.2,1.8)

\end{pspicture*}}

Let us denote with $G'$ resp. $G$ the group whose Dynkin diagram is $\cD'$
resp. $\cD$ and with $P'$ resp. $P$ the maximal parabolic subgroup of $G'$
resp. $G$ corresponding to the marked node $d$. We have a commutative
diagram: 
$$\xymatrix{G_0/P_0\ar@{=}[r]\ar@{^{(}->}[d]&G_0/P_0\ar@{^{(}->}[d]\\
G'/P'\ar@{^{(}->}[r]^\iota&G/P.}$$

We may define extended elements $w'$ and $(w_i)_{i\in[1,n]}$ of $w$ 
in $W_{G'}^{P'}$ and $W_G^P$ by $w'=s_\beta w$ and $w_i=s_{\beta_i} w$. 
Their length is $l(w)+1$.
For example let us consider 
$w=s_{\a_1}s_{\a_3}s_{\a_2}
s_{\a_4}s_{\a_3}s_{\a_2}s_{\a_1}$ in the Weyl group of $F_4$ 
(with notation as in \cite{bou}). This is a $\varpi_1$-cominuscule 
element. The elements $w'$ and $(w_i)_{i\in[1,n]}$ will also be 
$\varpi_1$-cominuscule. We depict here their heaps (in the following 
diagrams we depicted with crossed nodes the added vertices of $w'$ 
and $(w_i)_{i\in[1,n]}$).

\centerline{\begin{pspicture*}(0,0)(12.0,5.4)

\psellipse[fillstyle=solid,fillcolor=black](0.6,4.2)(0.096,0.096)
\psellipse[fillstyle=solid,fillcolor=black](0.6,1.8)(0.096,0.096)
\psellipse[fillstyle=solid,fillcolor=black](1.2,3.6)(0.096,0.096)
\psellipse[fillstyle=solid,fillcolor=black](1.2,2.4)(0.096,0.096)
\psellipse[fillstyle=solid,fillcolor=black](1.8,4.2)(0.096,0.096)
\psellipse[fillstyle=solid,fillcolor=black](1.8,3.0)(0.096,0.096)
\psellipse[fillstyle=solid,fillcolor=black](2.4,3.6)(0.096,0.096)
\psellipse[fillstyle=solid,fillcolor=black](4.2,4.2)(0.096,0.096)
\psellipse[fillstyle=solid,fillcolor=black](4.2,1.8)(0.096,0.096)
\psellipse[fillstyle=solid,fillcolor=black](4.8,3.6)(0.096,0.096)
\psellipse[fillstyle=solid,fillcolor=black](4.8,2.4)(0.096,0.096)
\psellipse[fillstyle=solid,fillcolor=black](5.4,4.2)(0.096,0.096)
\psellipse[fillstyle=solid,fillcolor=black](5.4,3.0)(0.096,0.096)
\psellipse[fillstyle=solid,fillcolor=black](6.0,3.6)(0.096,0.096)
\psellipse(6.6,4.2)(0.096,0.096)
\psline(6.532,4.132)(6.668,4.268)
\psline(6.532,4.268)(6.668,4.132)
\psellipse[fillstyle=solid,fillcolor=black](8.4,4.2)(0.096,0.096)
\psellipse[fillstyle=solid,fillcolor=black](8.4,1.8)(0.096,0.096)
\psellipse[fillstyle=solid,fillcolor=black](9.0,3.6)(0.096,0.096)
\psellipse[fillstyle=solid,fillcolor=black](9.0,2.4)(0.096,0.096)
\psellipse[fillstyle=solid,fillcolor=black](9.6,4.2)(0.096,0.096)
\psellipse[fillstyle=solid,fillcolor=black](9.6,3.0)(0.096,0.096)
\psellipse[fillstyle=solid,fillcolor=black](10.2,3.6)(0.096,0.096)
\psellipse(10.8,4.8)(0.096,0.096)
\psline(10.732,4.732)(10.868,4.868)
\psline(10.732,4.868)(10.868,4.732)
\psellipse(10.8,4.2)(0.096,0.096)
\psline(10.732,4.132)(10.868,4.268)
\psline(10.732,4.268)(10.868,4.132)
\psellipse(10.8,3.0)(0.096,0.096)
\psline(10.732,2.932)(10.868,3.068)
\psline(10.732,3.068)(10.868,2.932)
\psellipse(10.8,2.4)(0.096,0.096)
\psline(10.732,2.332)(10.868,2.468)
\psline(10.732,2.468)(10.868,2.332)
\psline(0.668,1.868)(1.132,2.332)
\psline(1.868,3.068)(2.332,3.532)
\psline(2.332,3.668)(1.868,4.132)
\psline(1.132,3.668)(0.668,4.132)
\psline(1.288,3.637)(1.763,4.112)
\psline(1.237,3.688)(1.712,4.163)
\psline(1.55,3.95)(1.43,3.94)
\psline(1.55,3.95)(1.54,3.83)
\psline(1.237,3.512)(1.712,3.037)
\psline(1.288,3.563)(1.763,3.088)
\psline(1.55,3.25)(1.54,3.37)
\psline(1.55,3.25)(1.43,3.26)
\psline(1.288,2.437)(1.763,2.912)
\psline(1.237,2.488)(1.712,2.963)
\psline(1.55,2.75)(1.43,2.74)
\psline(1.55,2.75)(1.54,2.63)
\psline(4.268,1.868)(4.732,2.332)
\psline(5.468,3.068)(5.932,3.532)
\psline(5.932,3.668)(5.468,4.132)
\psline(4.732,3.668)(4.268,4.132)
\psline(4.888,2.437)(5.363,2.912)
\psline(4.837,2.488)(5.312,2.963)
\psline(5.15,2.75)(5.03,2.74)
\psline(5.15,2.75)(5.14,2.63)
\psline(4.837,3.512)(5.312,3.037)
\psline(4.888,3.563)(5.363,3.088)
\psline(5.15,3.25)(5.14,3.37)
\psline(5.15,3.25)(5.03,3.26)
\psline(4.888,3.637)(5.363,4.112)
\psline(4.837,3.688)(5.312,4.163)
\psline(5.15,3.95)(5.03,3.94)
\psline(5.15,3.95)(5.14,3.83)
\psline(6.076,3.524)(6.676,4.124)
\psline(6.088,3.637)(6.563,4.112)
\psline(6.037,3.688)(6.512,4.163)
\psline(5.924,3.676)(6.524,4.276)
\psline(6.33,3.93)(6.21,3.92)
\psline(6.33,3.93)(6.32,3.81)
\psline(8.468,1.868)(8.932,2.332)
\psline(9.668,3.068)(10.132,3.532)
\psline(10.132,3.668)(9.668,4.132)
\psline(8.932,3.668)(8.468,4.132)
\psline(9.088,3.637)(9.563,4.112)
\psline(9.037,3.688)(9.512,4.163)
\psline(9.35,3.95)(9.23,3.94)
\psline(9.35,3.95)(9.34,3.83)
\psline(9.037,3.512)(9.512,3.037)
\psline(9.088,3.563)(9.563,3.088)
\psline(9.35,3.25)(9.34,3.37)
\psline(9.35,3.25)(9.23,3.26)
\psline(9.088,2.437)(9.563,2.912)
\psline(9.037,2.488)(9.512,2.963)
\psline(9.35,2.75)(9.23,2.74)
\psline(9.35,2.75)(9.34,2.63)
\psline(10.243,3.686)(10.757,4.714)
\psline(10.268,3.668)(10.732,4.132)
\psline(10.268,3.532)(10.732,3.068)
\psline(10.243,3.514)(10.757,2.486)
\uput[d](10.8,4.2){$\vdots$}
\put(6.6,3.6){$p$}
\uput[dr](1.4,1.2){$w$}
\uput[dr](5,1.2){$w'$}
\uput[dr](9,1.2){$(w_i)_{i\in[1,p]}$}
\uput[r](10.8,4.8){$w_1$}
\uput[r](10.8,2.4){$w_{p}$}

\end{pspicture*}}

For $w$ as above, we define $\s_w$ the corresponding homology
class in $G_0/P_0$ and also in $G'/P'$. We denote with $\tau_w$ the same class
in $H_*(G/P)$. We denote with $\s_{w'}$ the homology class in $G'/P'$
corresponding to $w'$ and with $\tau_{w_i}$ the homology class in $G/P$
corresponding to $w_i$ for $i\in[1,p]$.

\begin{prop}
\label{prop-iota_*}
  We have the equality $\displaystyle{\iota_*\s_{w'}=\sum_{i=1}^p
    \tau_{w_i}}$.
\end{prop}

\begin{proo}
We proceed by induction on the length of $w$. Let us write 
$$\iota_*\s_{w'}=\sum_{x\in W_G^P:\ l(x)=l(w)+1} b_x \tau_x.$$
Let us first of all prove that the only classes appearing in this sum
are the classes $(\tau_{w_i})_{i\in[1,p]}$.

\begin{lemm}
Let $x\in W^P$ with $b_x>0$, then we have
$x=w_i$ for some $i\in[1,p]$. 
\end{lemm}

\begin{proo}
Let us introduce some notation. Let us denote with $\delta$ the simple root
associated to the vertex $d$. We denote with $P'_{\b,\delta}$ 
and $P'_\b$ (resp. $P_{\bb,\delta}$ and $P_\bb$) the parabolic subgroups
of $G'$ (resp. $G$) associated to the set of simple roots
$\{\beta,\delta\}$ and $\{\beta\}$ (resp. $\{(\b_i)_{i\in[1,p]},\delta\}$ and
$\{(\b_i)_{i\in[1,p]}\}$). We also denote, for $u\in W_{G'}$
(resp. $v\in W_G$), with $X_{\beta,\delta}(u)$ and $X_\beta(u)$
(resp. $X_{\bb,\delta}(v)$ and $X_\bb(v)$) the associated Schubert
varieties  in $G'/P'_{\beta,\delta}$ and $G'/P'_{\beta}$
(resp. $G/P_{\bb,\delta}$ and $G/P_\bb$). Finally we introduce the
projections  $p':G'/P'_{\beta,\delta}\to G'/P'$ and
$q':G'/P'_{\beta,\delta}\to G'/P'_{\beta}$ (resp. $p:G/P_{\bb,\delta}\to
G/P$ and $q:G/P_{\bb,\delta}\to G/P_{\bb}$).

Choose a reduced expression $s_{\a_1}\cdots s_{\a_l}$ for $w$ with
$\a_i$ simple roots of $G_0$. We must have the equality $\a_l=\delta$. We deduce a reduced expression 
$w'=s_\beta s_{\a_1}\cdots s_{\a_l}$. Let us consider the unipotent
subgroup $U_{w}=U_{\a_1}\cdots U_{\a_l}$ of $G_0$ and the unipotent
subgroup $U_{w'}=U_\beta U_w$ of $G'$. We have an inclusion 
$U_{w'}\subset U_{\b_1}\cdots U_{\b_p}U_w$. This induces the following inclusions of 
Schubert varieties $\iota: X(w')\subset X(s_{\b_1}\cdots s_{\b_p}w)$, $\iota_\beta:
X_\b(s_\b)\subset X_\bb(s_{\b_1}\cdots s_{\b_p})$ 
and $\iota_{\beta,\delta}: X_{\beta,\delta}(w')\subset
X_{\bb,\delta}(s_{\b_1}\cdots s_{\b_p}w)$. We have the commutative diagram:
$$\xymatrix{&X_{\beta}(s_\beta)\ar@{^{(}->}[r]^{\iota_\beta}&X_{\bb}(s_{\b_1}\cdots
  s_{\b_p}) \\
X_{\beta,\delta}(w')\ar[ur]^{q'}\ar@{^{(}->}[r]^{\iota_{\beta,\delta}\ \ \ }\ar[d]^{p'}
 &
  X_{\bb,\delta}(s_{\b_1}\cdots s_{\b_p}w)\ar[ur]^{q}\ar[d]^{p} &\\
X(w')\ar@{^{(}->}[r]^{\iota\ \ \ }&
  X(s_{\b_1}\cdots s_{\b_p}w). &\\}$$
Remark that the Schubert variety $X_\beta(s_\beta)$ is isomorphic to the
projective line $\pu$
  while the Schubert variety $X_\bb(s_{\b_1}\cdots s_{\b_p})$ is
  isomorphic to $(\pu)^p$ the map $\iota_\beta$ being given by the diagonal
  embedding.

Let $\tau_x$ be a class with $b_x>0$. We thus have $x\leq
s_{\b_1}\cdots s_{\b_p}w$. In particular, as any reduced expression
for $s_{\b_1}\cdots s_{\b_p}w$ is obtained by multiplying on the left
with $s_{\b_1}\cdots s_{\b_p}$ a reduced expression for $w$, we
obtain (using the characterisation of Bruhat order described in
\cite[Section 3 Proposition 5]{demazure}) that  
\begin{equation}
\label{formule-int}
x=\prod_{k\in A}s_{\b_k}y
\end{equation}
with $A\subset[1,p]$ and $y\leq w$. The same argument gives that if we
write 
$$(\iota_{\beta,\delta})_*[X_{\beta,\delta}(w')] = \sum_{t\in W_G^{P_{\bb,\delta}}
  : \ l(t)=l(w)+1} c_t \cdot [X_{\bb,\delta}(t)],$$
then $c_t>0$ implies 
\begin{equation}
\label{formule-int-bis}
t=\prod_{k\in B}s_{\b_k}u
\end{equation}
with $B\subset[1,p]$ and $u\leq w$. We now prove that $A$ has
at most one element and for this, we prove that $B$ has at most one
element. 

Let $[X_{\bb,\delta}(t)]$ be a class with $c_t>0$ and assume that in the
expression (\ref{formule-int-bis}) the set $B$ contains at least two
elements say $i$ and $j$ in $[1,p]$. Let us consider the two
  degree one cohomology classes $h_i$ and $h_j$ of $(\pu)^p$
  corresponding to the factors $i$ and $j$. We have $(h_i\cup h_j)\cap
  [X_{\bb,\delta}(t)]\neq0$ by Chevalley formula and because
  $c_t>0$ we get $(h_i\cup h_j)\cap
  (\iota_{\beta,\delta})_*[X_{\beta,\delta}(w')] \neq0$. By projection formula we get
  ${\iota_{\b,\delta}}_*(\iota_{\b,\delta}^*(h_i\cup h_j)\cap 
  [X_{\beta,\delta}(w')])\neq 0$. On the other hand we have
  ${\iota_{\b,\delta}}_*(\iota_{\b,\delta}^*(h_i\cup h_j)\cap
  [X_{\beta,\delta}(w')])=0$ because $\iota_{\b,\delta}^*(h_i\cup
  h_j)=\iota_{\b,\delta}^*q^*(h_i\cup h_j)={q'}^*\iota_{\b}^*(h_i\cup h_j)$
  and $\iota_{\b}^*(h_i\cup h_j)$ vanishes as a degree 2 class on $\pu$, a
  contradiction. Thus $B$ has at most one element.

Because the maps $p'$ is birational, we have
$p'_*[X_{\beta,\delta}(w')]=[X(w')]=\s_{w'}$ and thus the
equality 
$$\iota_*\s_{w'}=p_*(\iota_{\beta,\delta})_*[X_{\beta,\delta}(w')]=\sum_{t\in
  W_G^{P_{\bb,\delta}} : \ l(t)=l(w)+1} c_t \cdot p_*[X_{\bb,\delta}(t)].$$
Now for $t\in W_G^{P_{\bb,\delta}}$, we have 
$$p_*[X_{\bb,\delta}(t)]=\left\{\begin{array}{ll}
\tau_t& \textrm{for $t\in W_G^P$}\\
0 & \textrm{otherwise.}\\
\end{array}\right.
$$
We deduce that $A$ has only one element. Now from (\ref{formule-int})
and the fact that, $b_x>0$ implies that $l(x)=l(w)+1$, the result
follows.
\end{proo}

We deduce that there is an
integer $i$ such that $b_{w_i} > 0$. Furthermore, the group $G'$ is
obtained from $G$ by taking the subgroup invariant by an automorphism
of order $p$ of the Dynkin diagram $\cD$: the permutation of the $p$
vertices we added to $\cD_0$. In particular the class $\iota_*\s_{w'}$
is invariant under this permutation thus we have the equalities
$b_{w_i}=b_{w_j}$ for $i$ and $j$ in $[1,p]$. We may therefore set
$b=b_{w_i}$ for any $i\in[1,p]$, we have $b>0$ and 
$$\iota_*\s_{w'}=b\sum_{i=1}^p\tau_{w_i}.$$
Computing the coefficient of $\tau_w$ in
$h \cap \iota_*\s_{w'} = \iota_* ( h \cap \s_{w'} )$, we get the
equality 
$$bp\frac{(\a,\a)}{(\delta,\delta)} 
=p\frac{(\a,\a)}{(\delta,\delta)},$$
thus $b=1$. 
\end{proo}

\subsection{Type $C_n $ case}

In this case, we consider  the system of $\varpi_n$-cominuscule 
$C_n$-colored posets $\pos_0$ given by the posets of a Lagrangian
Grassmannian $\G_\omega(n,2n)$.
We have $S_0=\{1\}$. Let $(D_1,d_1)$ be a
marked Dynkin diagram and $P_1$ be any  $d_1$-minuscule $D_1$-colored
poset. Set $\pos=\pos_{\pos_0,\{P_{1}\}}$. 
The heap $\pos$ for type $C_6$ is the same as the one for type $D_7$ 
except for the colors. It was described in (\ref{equa-systeme-d7}).

\begin{lemm}
\label{lemm-cnpn}
With the above notation, assume that Conjecture \ref{main_conj} holds
for $P_{1}$ and any $\lt$ in $I(\pos)$ with
$D_0(\lt)\varsubsetneq C_n$. Then Conjecture \ref{main_conj} holds for
$\pos$. 
\end{lemm}

\begin{proo}
Let us define the degree $i$ ideals $\lt_{i}=\scal{(\a_{n+1-i},1)}$
for $i\in[1,n]$ and set $\s^{i}=s^{\lt_i}$. Take 
a set of generators $\{\gamma^1,\cdots,\gamma^{n}\}$ with $\deg(\gamma^i)=i$. 
We start to prove that the generators $(\gamma^i)_{i\in[1,n-1]}$ are good 
generators and shall prove at the end that $\gamma^n$ is also a good 
generator.

Since by assumption the conjecture holds for any $\lt\in I(\pos)$ with
$D_0(\lt)\varsubsetneq C_n$, we have $c_{\lt,\mu}^\nu=t_{\lt,\mu}^\nu$
as soon as $\deg(\nu\cap\pos_0)\leq2n-2$. In particular if 
$\deg(\lt)\leq n$ and $i\leq n-2$ we have $\deg(\lt)+i\leq 2n-2$ and
$\gamma^i\cdot\s^\lt=\gamma^i\pp\s^\lt$. Furthermore, for $i=n-1$ there 
is a unique 
ideal $\nu$ (namely $\nu=\scal{(\a_2,2)}$) of degree $2n-1$ for which we 
cannot compute 
$c_{\gamma^{n-1},\lt}^\nu$. By Lemma \ref{poincare} we conclude that 
$\gamma^{n-1}\cdot\s^\lt = \gamma^{n-1}\pp\s^\lt$. In particular we have 
$\gamma^i\cdot\s^j=\gamma^i\pp\s^j$ for $i\in[1,n-1]$ and $j\in[1,n]$.

If $\lt\supset\lt_{n}$, then by recursion with respect to $\s^{n}$
we have $\gamma^i\cdot\s^\lt=\gamma^i\pp\s^\lt$. 

If $\lt\not\supset\lt_{n}$, then we first consider the case where $\lt$ is
an ideal of the form $\scal{(\a_k,l)}$ for some simple root
$\a_k$ and some integer $l$. We prove the equality
$\gamma^i\cdot\s^\lt=\gamma^i\pp\s^\lt$ by induction on $\deg(\lt)$ 
in that case. We may of course assume that $\lt$ is distinct from all
the $\lt_i$. 
We 
consider the two
subideals $\lt'$ and $\lt''$ 
in $\lt$ described by $\scal{(\a_{k-1},l')}$ and
$\scal{(\a_{k+1},l'')}$ (if $k=n$ we consider only $\lt'$)
where $l'=\max\{a\ /\ (\a_{k-1},a)\in \lt\}$ and 
$l''=\max\{a\ /\ (\a_{k+1},a)\in \lt\}$.
By recursion with respect to
$\lt'$ or $\lt''$, we have
$c_{\gamma^i,\s^\lt}^{\s^\nu}=t_{\gamma^i,\s^\lt}^{\s^\nu}$ for any $\nu$ 
not containing
$(\a_{k-1},l'+1)$ or $(\a_{k+1},l''+1)$ (the last condition is empty for 
$k=n$).
By induction
on $\pos_0$ it is also true if $\nu$ does not contain $(\a_1,1)$. For
an ideal $\nu$ in $\pos$ containing all these elements of $\pos_0$, we have
$\deg(\nu\cap \pos_0)\geq \deg(\lt)+n-1$. For such a $\nu$ and $i\leq n-2$, 
we have
$c_{\gamma^i,\s^\lt}^{\s^\nu}=0=t_{\gamma^i,\s^\lt}^{\s^\nu}$
for degree reasons. For $i=n-1$ however, the equality
$c_{\gamma^i,\s^\lt}^{\s^\nu}=t_{\gamma^i,\s^\lt}^{\s^\nu}$
holds for all $\lt=\scal{(\a_k,l)}$ and $\nu\neq\scal{(\a_k,l+1)}$. 
We conclude 
by Lemma \ref{poincare}.

We finish by dealing with $\lt$ not of the previous form. Let us
consider the set $M(\lt)$ of maximal elements in $\lt$. For $(\a_k,l)\in
M(\lt)$, define the ideal $\lt(\a_k,l)=\scal{(\a_k,l)}$. We have
$\gamma^i\cdot \s^{\lt(\a_k,l)}=\gamma^i\pp\s^{\lt(\a_k,l)}$. 
In particular we can
use recursion with respect to $\lt(\a_k,l)$ and we deduce that
$c_{\gamma^i,\s^\lt}^{\s^\nu}=t_{\gamma^i,\s^\lt}^{\s^\nu}$
for any $\nu$ not containing
$(\a_{k},l+1)$. By induction on $\pos_0$ it is also true if $\nu$ does
not contain $(\a_1,1)$. For an ideal $\nu$ in $\pos$ containing all the 
elements
$(\a_k,l+1)$ for $(\a_k,l)\in M(\lt)$ as well as $(\a_1,1)$, we have
$\deg(\nu\cap \pos_0)\geq \deg(\lt)+n$. For such a $\nu$ and $i\leq n-1$, 
we have
$c_{\gamma^i,\s^\lt}^{\s^\nu}=0=t_{\gamma^i,\s^\lt}^{\s^\nu}$
for degree reasons. 

To finish the proof, we need to deal with $\gamma^n$. The first formula we 
need to verify is the equality $\ok{\g^n}{\g^n}$. This will be the most 
difficult one. Indeed, assume this formula holds, then
$\ok{\g^n}{\s^n}$ and by recursion $\ok{\g^n}{\s^\lt}$ 
for $\lt\supset\lt_n$. Now take $\lt\not\supset\lt_n$, then in the
cohomology of $G_\omega(n-1,2(n-1))$ 
we may write $\s^\lt=P(\g^1,\cdots,\g^{n-1})$
where $P$ is a polynomial in $n-1$ variables. If we consider the class
$P(\g^1,\cdots,\g^{n-1})$ in $H^*(\pos)$ then its pull-back to
$H^*(G_\omega(n-1,2(n-1)))$ is $\s^\lt$ thus
$P(\g^1,\cdots,\g^{n-1})=\s^\lt+A$ where $A$ is a linear combination
of classes $\s^\mu$ with $\mu\supset\lt_n$. We thus have
$\ok{\g^n}{A}$. Furthermore, by Lemma
\ref{lemm-petite-puissance} we have
$\ok{\g^n}{P(\g^1,\cdots,\g^{n-1})}$ and the result follows.

To prove $\ok{\g^n}{\g^n}$, we remark that there are two ideals $\nu$
of degree $2n$ for which we do not know that
$\oknu{\g^n}{\g^n}{\s^\nu}$. These ideals are
$\nu=\scal{(\a_2,2),(\a_n,3)}$ and $\nu'=\scal{(\a_0,1),(\a_2,2)}$
where we denote with $\a_0$ the simple root corresponding to the vertex
$d_1$ in $D_1$. Since $\nu$ is contained in $\pos_0$ and is the
only class in that degree in $\pos_0$, we may apply Lemma
\ref{poincare} to get $\oknu{\g^n}{\g^n}{\s^\nu}$. For $\nu'$ however,
we may not apply Lemma \ref{poincare} since $D_0(\nu')$ is not the
Dynkin diagram of a finite group. However if the edge between $\a_0$
and $\a_1$ is simple, then $D_0(\nu')=C_{n+1}$ is of finite type and
$\oknu{\g^n}{\g^n}{\s^{\nu'}}$ by Lemma \ref{poincare}. If the edge between
$\a_0$ and $\a_1$ is a $p$-tuple edge (i.e. $\scal{\a_1,\a_0^\vee}=p$),
then by Proposition \ref{prop-iota_*} we have
$$\iota_*\s_{\nu'}=\sum_{i=1}^p\tau_{\nu_i}$$ 
with notation as in Proposition \ref{prop-iota_*}.
We then have, because $\iota^*\g^n=\g^n$, the equality 
$$\iota_*(\g^n\cap\s_{\nu'})=\sum_{i=1}^p\g^n\cap\tau_{\nu_i}$$
and it follows that
$\displaystyle{c_{\g^n,\g^n}^{\s^{\nu'}}= 
\sum_{i=1}^pc_{\g^n,\g^n}^{\tau^{\nu_i}}}$
and the result follows.
\end{proo}

\subsection{Type $F_4$ cominuscule case}

As for type $C_n$ we shall need to use Proposition \ref{prop-iota_*}
and foldings to get the result. However, we need here one more
step. Indeed, with the notation of Proposition \ref{prop-iota_*}, if
$D=C_n$ then $D'=C_{n+1}$ is still of finite type with quite well
understood cohomology, for $D=F_4$, then $D'=\widetilde{F}_4^2$ which
is a twisted affine Dynkin diagram (see \cite{kac}). To compute some 
intersections in its cohomology we will use a folding of 
$\widetilde{E}_7^1$ to $\widetilde{F}_4^2$ and compute direct images 
by hand (this is done in Lemma \ref{lemm-f_4} and in Proposition 
\ref{prop-f_4}).

\subsubsection{Foldings with $F_4$}
 
We start with notation and set up. Let us denote with $\iota$ the
inclusion of the group $F_4$ in the group $E_6$ given by folding of
the Dynkin diagram. We also denote with $\iota$ the inclusion of
$F_4/P_1$ in $E_6/P_2$. We want to describe the map
$\iota_*:H_*(F_4/P_1)\to H_*(E_6/P_2)$. For this we introduce some 
notation to
describe the classes in these homology groups. Let $\Lambda_F$ and
$\Lambda_E$ the fundamental weights corresponding to $F_4/P_1$ and $E_6/P_2$
respectively. Any element of length at most 7 in $(W_{F_4})^{P_1}$ is
$\Lambda_F$-cominuscule. The two heaps of size 7 are as follows:

\centerline{\begin{pspicture*}(0,0)(6.6000004,5.4)

\psellipse[fillstyle=solid,fillcolor=black](0.6,4.2)(0.096,0.096)
\psellipse[fillstyle=solid,fillcolor=black](0.6,1.8)(0.096,0.096)
\psellipse[fillstyle=solid,fillcolor=black](1.2,3.6)(0.096,0.096)
\psellipse[fillstyle=solid,fillcolor=black](1.2,2.4)(0.096,0.096)
\psellipse[fillstyle=solid,fillcolor=black](1.8,4.2)(0.096,0.096)
\psellipse[fillstyle=solid,fillcolor=black](1.8,3.0)(0.096,0.096)
\psellipse[fillstyle=solid,fillcolor=black](2.4,3.6)(0.096,0.096)
\psellipse[fillstyle=solid,fillcolor=black](4.2,1.8)(0.096,0.096)
\psellipse[fillstyle=solid,fillcolor=black](4.8,4.8)(0.096,0.096)
\psellipse[fillstyle=solid,fillcolor=black](4.8,3.6)(0.096,0.096)
\psellipse[fillstyle=solid,fillcolor=black](4.8,2.4)(0.096,0.096)
\psellipse[fillstyle=solid,fillcolor=black](5.4,4.2)(0.096,0.096)
\psellipse[fillstyle=solid,fillcolor=black](5.4,3.0)(0.096,0.096)
\psellipse[fillstyle=solid,fillcolor=black](6.0,3.6)(0.096,0.096)
\psline(0.668,1.868)(1.132,2.332)
\psline(1.868,3.068)(2.332,3.532)
\psline(2.332,3.668)(1.868,4.132)
\psline(4.268,1.868)(4.732,2.332)
\psline(1.132,3.668)(0.668,4.132)
\psline(5.468,3.068)(5.932,3.532)
\psline(5.932,3.668)(5.468,4.132)
\psline(4.888,2.437)(5.363,2.912)
\psline(4.837,2.488)(5.312,2.963)
\psline(5.15,2.75)(5.03,2.74)
\psline(5.15,2.75)(5.14,2.63)
\psline(4.837,3.512)(5.312,3.037)
\psline(4.888,3.563)(5.363,3.088)
\psline(5.15,3.25)(5.14,3.37)
\psline(5.15,3.25)(5.03,3.26)
\psline(4.888,3.637)(5.363,4.112)
\psline(4.837,3.688)(5.312,4.163)
\psline(5.15,3.95)(5.03,3.94)
\psline(5.15,3.95)(5.14,3.83)
\psline(4.837,4.712)(5.312,4.237)
\psline(4.888,4.763)(5.363,4.288)
\psline(5.15,4.45)(5.14,4.57)
\psline(5.15,4.45)(5.03,4.46)
\psline(1.288,2.437)(1.763,2.912)
\psline(1.237,2.488)(1.712,2.963)
\psline(1.55,2.75)(1.43,2.74)
\psline(1.55,2.75)(1.54,2.63)
\psline(1.237,3.512)(1.712,3.037)
\psline(1.288,3.563)(1.763,3.088)
\psline(1.55,3.25)(1.54,3.37)
\psline(1.55,3.25)(1.43,3.26)
\psline(1.288,3.637)(1.763,4.112)
\psline(1.237,3.688)(1.712,4.163)
\psline(1.55,3.95)(1.43,3.94)
\psline(1.55,3.95)(1.54,3.83)
\uput[d](4.8,1.2){$\s_{7,2}$}
\uput[d](1.2,1.2){$\s_{7,1}$}

\end{pspicture*}}

To fix notation we define the following homology classes in
$H_*(F_4/P_1)$. These are all classes of degree $d\in[4,7]$. By 
convention the notation $\s_{a,b}$ or $\tau_{a,b}$ (resp. 
$\s^{a,b}$ or $\tau^{a,b}$) denote homology (resp. cohomology) classes 
of degree a. The set of all indices $b$ is an index set of (co)homology 
classes of that degree. For example, as the following array shows, 
there are two homology classes of degree 4 denoted $\s_{4,1}$ and 
$\s_{4,2}$. 
$$\begin{array}{llll}
  \s_{4,1}=\scal{(\a_2,2)}&
  \s_{5,1}=\scal{(\a_1,2)} &
  \s_{6,1}=\scal{(\a_1,2),(\a_4,1)} &
  \s_{7,1}=\scal{(\a_1,2),(\a_3,2)}
  \\
\s_{4,2}=\scal{(\a_4,1)} &
\s_{5,2}=\scal{(\a_2,2),(\a_4,1)}&
\s_{6,2}=\scal{(\a_3,2)} &
\s_{7,2}=\scal{(\a_2,3)} \\
\end{array}$$

The two elements of length 8 are fully commutative. The heaps of these
length 8 elements are as follows:

\centerline{\begin{pspicture*}(0,0)(6.6000004,6.0)

\psellipse[fillstyle=solid,fillcolor=black](0.6,4.2)(0.096,0.096)
\psellipse[fillstyle=solid,fillcolor=black](0.6,1.8)(0.096,0.096)
\psellipse[fillstyle=solid,fillcolor=black](1.2,4.8)(0.096,0.096)
\psellipse[fillstyle=solid,fillcolor=black](1.2,3.6)(0.096,0.096)
\psellipse[fillstyle=solid,fillcolor=black](1.2,2.4)(0.096,0.096)
\psellipse[fillstyle=solid,fillcolor=black](1.8,4.2)(0.096,0.096)
\psellipse[fillstyle=solid,fillcolor=black](1.8,3.0)(0.096,0.096)
\psellipse[fillstyle=solid,fillcolor=black](2.4,3.6)(0.096,0.096)
\psellipse[fillstyle=solid,fillcolor=black](4.2,5.4)(0.096,0.096)
\psellipse[fillstyle=solid,fillcolor=black](4.2,1.8)(0.096,0.096)
\psellipse[fillstyle=solid,fillcolor=black](4.8,4.8)(0.096,0.096)
\psellipse[fillstyle=solid,fillcolor=black](4.8,3.6)(0.096,0.096)
\psellipse[fillstyle=solid,fillcolor=black](4.8,2.4)(0.096,0.096)
\psellipse[fillstyle=solid,fillcolor=black](5.4,4.2)(0.096,0.096)
\psellipse[fillstyle=solid,fillcolor=black](5.4,3.0)(0.096,0.096)
\psellipse[fillstyle=solid,fillcolor=black](6.0,3.6)(0.096,0.096)
\psline(0.668,1.868)(1.132,2.332)
\psline(1.868,3.068)(2.332,3.532)
\psline(2.332,3.668)(1.868,4.132)
\psline(1.132,3.668)(0.668,4.132)
\psline(0.668,4.268)(1.132,4.732)
\psline(1.288,2.437)(1.763,2.912)
\psline(1.237,2.488)(1.712,2.963)
\psline(1.55,2.75)(1.43,2.74)
\psline(1.55,2.75)(1.54,2.63)
\psline(1.237,3.512)(1.712,3.037)
\psline(1.288,3.563)(1.763,3.088)
\psline(1.55,3.25)(1.54,3.37)
\psline(1.55,3.25)(1.43,3.26)
\psline(1.288,3.637)(1.763,4.112)
\psline(1.237,3.688)(1.712,4.163)
\psline(1.55,3.95)(1.43,3.94)
\psline(1.55,3.95)(1.54,3.83)
\psline(1.237,4.712)(1.712,4.237)
\psline(1.288,4.763)(1.763,4.288)
\psline(1.55,4.45)(1.54,4.57)
\psline(1.55,4.45)(1.43,4.46)
\psline(4.268,1.868)(4.732,2.332)
\psline(5.468,3.068)(5.932,3.532)
\psline(5.932,3.668)(5.468,4.132)
\psline(4.732,4.868)(4.268,5.332)
\psline(4.837,4.712)(5.312,4.237)
\psline(4.888,4.763)(5.363,4.288)
\psline(5.15,4.45)(5.14,4.57)
\psline(5.15,4.45)(5.03,4.46)
\psline(4.888,3.637)(5.363,4.112)
\psline(4.837,3.688)(5.312,4.163)
\psline(5.15,3.95)(5.03,3.94)
\psline(5.15,3.95)(5.14,3.83)
\psline(4.837,3.512)(5.312,3.037)
\psline(4.888,3.563)(5.363,3.088)
\psline(5.15,3.25)(5.14,3.37)
\psline(5.15,3.25)(5.03,3.26)
\psline(4.888,2.437)(5.363,2.912)
\psline(4.837,2.488)(5.312,2.963)
\psline(5.15,2.75)(5.03,2.74)
\psline(5.15,2.75)(5.14,2.63)
\uput[d](1.2,1.2){$\s_{8,1}$}
\uput[d](4.8,1.2){$\s_{8,2}$}

\end{pspicture*}}

We define $\s_{8,1}$ to be the class associated to the left heap and
$\s_{8,2}$ to be the class associated to the right one. Let us also give 
the Hasse diagram for $F_4/P_1$. In the following picture we decribe on 
the lowest raw the classes $\s_{i,1}$ and on the top raw the classes 
$\s_{i,2}$ with $i$ growing from left to right. We also indicated the 
degree (with respect to the hyperplane classe) of the lower dimension 
classes.

\centerline{\begin{pspicture*}(0,0)(10.200001,3.0)

\psellipse[fillstyle=solid,fillcolor=black](0.6,1.2)(0.096,0.096)
\psellipse[fillstyle=solid,fillcolor=black](1.2,1.2)(0.096,0.096)
\psellipse[fillstyle=solid,fillcolor=black](1.8,1.2)(0.096,0.096)
\psellipse[fillstyle=solid,fillcolor=black](2.4,1.2)(0.096,0.096)
\psellipse[fillstyle=solid,fillcolor=black](3.0,1.8)(0.096,0.096)
\psellipse[fillstyle=solid,fillcolor=black](3.0,1.2)(0.096,0.096)
\psellipse[fillstyle=solid,fillcolor=black](3.6,1.8)(0.096,0.096)
\psellipse[fillstyle=solid,fillcolor=black](3.6,1.2)(0.096,0.096)
\psellipse[fillstyle=solid,fillcolor=black](4.2,1.8)(0.096,0.096)
\psellipse[fillstyle=solid,fillcolor=black](4.2,1.2)(0.096,0.096)
\psellipse[fillstyle=solid,fillcolor=black](4.8,1.8)(0.096,0.096)
\psellipse[fillstyle=solid,fillcolor=black](4.8,1.2)(0.096,0.096)
\psellipse[fillstyle=solid,fillcolor=black](5.4,1.8)(0.096,0.096)
\psellipse[fillstyle=solid,fillcolor=black](5.4,1.2)(0.096,0.096)
\psellipse[fillstyle=solid,fillcolor=black](6.0,1.8)(0.096,0.096)
\psellipse[fillstyle=solid,fillcolor=black](6.0,1.2)(0.096,0.096)
\psellipse[fillstyle=solid,fillcolor=black](6.6,1.8)(0.096,0.096)
\psellipse[fillstyle=solid,fillcolor=black](6.6,1.2)(0.096,0.096)
\psellipse[fillstyle=solid,fillcolor=black](7.2,1.8)(0.096,0.096)
\psellipse[fillstyle=solid,fillcolor=black](7.2,1.2)(0.096,0.096)
\psellipse[fillstyle=solid,fillcolor=black](7.8,1.2)(0.096,0.096)
\psellipse[fillstyle=solid,fillcolor=black](8.4,1.2)(0.096,0.096)
\psellipse[fillstyle=solid,fillcolor=black](9.0,1.2)(0.096,0.096)
\psellipse[fillstyle=solid,fillcolor=black](9.6,1.2)(0.096,0.096)
\psline(0.696,1.2)(1.104,1.2)
\psline(1.296,1.2)(1.704,1.2)
\psline(2.496,1.2)(2.904,1.2)
\psline(3.096,1.2)(3.504,1.2)
\psline(3.096,1.8)(3.504,1.8)
\psline(3.668,1.732)(4.132,1.268)
\psline(4.268,1.732)(4.732,1.268)
\psline(4.296,1.8)(4.704,1.8)
\psline(9.504,1.2)(9.096,1.2)
\psline(8.904,1.2)(8.496,1.2)
\psline(7.704,1.2)(7.296,1.2)
\psline(7.104,1.2)(6.696,1.2)
\psline(1.889,1.164)(2.311,1.164)
\psline(1.889,1.236)(2.311,1.236)
\psline(2.488,1.237)(2.963,1.712)
\psline(2.437,1.288)(2.912,1.763)
\psline(3.088,1.237)(3.563,1.712)
\psline(3.037,1.288)(3.512,1.763)
\psline(8.311,1.236)(7.889,1.236)
\psline(8.311,1.164)(7.889,1.164)
\psline(7.763,1.288)(7.288,1.763)
\psline(7.712,1.237)(7.237,1.712)
\psline(7.163,1.288)(6.688,1.763)
\psline(7.112,1.237)(6.637,1.712)
\psline(3.689,1.764)(4.111,1.764)
\psline(3.689,1.836)(4.111,1.836)
\psline(3.689,1.164)(4.111,1.164)
\psline(3.689,1.236)(4.111,1.236)
\psline(4.289,1.164)(4.711,1.164)
\psline(4.289,1.236)(4.711,1.236)
\psline(4.889,1.164)(5.311,1.164)
\psline(4.889,1.236)(5.311,1.236)
\psline(5.489,1.164)(5.911,1.164)
\psline(5.489,1.236)(5.911,1.236)
\psline(6.089,1.164)(6.511,1.164)
\psline(6.089,1.236)(6.511,1.236)
\psline(4.889,1.764)(5.311,1.764)
\psline(4.889,1.836)(5.311,1.836)
\psline(6.511,1.836)(6.089,1.836)
\psline(6.511,1.764)(6.089,1.764)
\psline(7.104,1.8)(6.696,1.8)
\psline(5.904,1.8)(5.496,1.8)
\psline(5.332,1.732)(4.868,1.268)
\psline(4.868,1.732)(5.332,1.268)
\psline(5.932,1.732)(5.468,1.268)
\psline(6.532,1.732)(6.068,1.268)
\uput[d](0.6,1.2){1}
\uput[d](1.2,1.2){1}
\uput[d](1.8,1.2){1}
\uput[d](2.4,1.2){2}
\uput[d](3.0,1.2){2}
\uput[d](3.6,1.2){2}
\uput[u](3.0,1.8){4}
\uput[u](3.6,1.8){8}
\uput[d](4.2,1.2){12}
\uput[d](4.8,1.2){40}
\uput[u](4.2,1.8){16}
\uput[u](4.8,1.8){16}

\end{pspicture*}}

Let us now
describe some classes in $E_6/P_2$. Recall that we described the
maximal slant-irreducible heap in $E_6/P_2$ in Section
\ref{section-e_6}. To fix notation we define the following homology classes in
$H_*(E_6/P_2)$. These are all classes of degree $d\in[3,8]$. 
$$\begin{array}{lll}
  \tau_{3,1}=\scal{(\beta_3,1)}&  
  \tau_{4,1}=\scal{(\beta_1,1)} &
  \tau_{5,1}=\scal{(\beta_1,1),(\beta_5,1)}\\
  \tau_{3,2}=\scal{(\beta_5,1)}&  
  \tau_{4,2}=\scal{(\beta_3,1),(\beta_5,1)} &
  \tau_{5,2}=\scal{(\beta_4,2)}\\
     &  
  \tau_{4,3}=\scal{(\beta_6,1)} &
  \tau_{5,3}=\scal{(\beta_3,1),(\beta_6,1)}\vspace{1 mm}\\
  \tau_{6,1}=\scal{(\beta_1,1),(\beta_4,2)} &
  \tau_{7,1}=\scal{(\beta_3,2)} &
  \tau_{8,1}=\scal{(\beta_3,2),(\beta_2,2)}  \\
  \tau_{6,2}=\scal{(\beta_2,2)} &
  \tau_{7,2}=\scal{(\beta_1,1),(\beta_2,2)} &
  \tau_{8,2}=\scal{(\beta_3,2),(\beta_6,1)}  \\
  \tau_{6,3}=\scal{(\beta_1,1),(\beta_6,1)} &
  \tau_{7,3}=\scal{(\beta_1,1),(\beta_4,2),(\beta_6,1)} &
  \tau_{8,3}=\scal{(\beta_1,1),(\beta_2,2),(\beta_6,1)}  \\
  \tau_{6,4}=\scal{(\beta_6,1),(\beta_4,2)} &
  \tau_{7,4}=\scal{(\beta_2,2),(\beta_6,1)} &
  \tau_{8,4}=\scal{(\beta_1,1),(\beta_5,2)}  \\
   &
  \tau_{7,5}=\scal{(\beta_5,2)} &
  \tau_{8,5}=\scal{(\beta_5,2),(\beta_2,2)}  \\
\end{array}$$

\begin{lemm}
\label{lemm-f_4}
Let $\iota$ denote the inclusion of $F_4/P_1$ into $E_6/P_2$. We have
$$\begin{array}{lll}
  \iota_*\s_{4,1}=\tau_{4,2}&
  \iota_*\s_{5,1}=\tau_{5,2} &
  \iota_*\s_{6,1}=\tau_{6,1}+\tau_{6,2}+\tau_{6,4} \\
  \iota_*\s_{4,2}=\tau_{4,1}+\tau_{4,2}+\tau_{4,3} &
  \iota_*\s_{5,2}=\tau_{5,1}+\tau_{5,2}+\tau_{5,3}&
  \iota_*\s_{6,2}=\tau_{6,1}+\tau_{6,3}+\tau_{6,4} \\
  \iota_*\s_{7,1}=\tau_{7,1}+\tau_{7,2}+\tau_{7,3}+
\tau_{7,4}+\tau_{7,5} &
  \multicolumn{2}{l}{\iota_*\s_{8,1}=\tau_{8,1}+\tau_{8,2}+\tau_{8,3}+
\tau_{8,4}+\tau_{8,5}}\\
  \iota_*\s_{7,2}=\tau_{7,3} &
  \iota_*\s_{8,2}=\tau_{8,2}+\tau_{8,3}+\tau_{8,4} &\\
\end{array}$$
\end{lemm}
\begin{proo}
We shall denote with $h$ the hyperplane class in $H^*(E_6/P_2)$ and in
$H^*(F_4/P_1)$ by identifying it to its pull-back. Let $g$ be the Weyl
involution of the Lie algebra ${\mathfrak e_6}$. Then $g$ induces an
outer automorphism of $E_6/P_2$, which fixes pointwise
$\iota(F_4/P_1)$. Since $g \circ \iota = \iota$, we have $g_* \iota_*
\s = \iota_* \s$ for $\s \in H_*(F_4/P_1)$. In other words, the
classes in the image of $\iota_*$ are invariant under $g$. 

Thus there exist non negative integers $a,b,c,d$ such that
$$
\left \{
\begin{array}{rcl}
\iota_* \s_{4,1} & = & a(\tau_{4,1}+\tau_{4,3}) + b \tau_{4,2} \\
\iota_* \s_{4,2} & = & c(\tau_{4,1}+\tau_{4,3}) + d \tau_{4,2}.
\end{array}
\right .
$$
By the same argument there exist non negative integers 
$\a,\b,\gamma,\delta,\epsilon,\eta$ such that
$$
\left \{
\begin{array}{rcl}
  \iota_* \s_{8,1} & = & \a(\tau_{8,1}+\tau_{8,5}) + \b (\tau_{8,2} +
  \tau_{8,4}) + \gamma \tau_{8,3}\\
\iota_* \s_{8,2} & = & \delta(\tau_{8,1}+\tau_{8,5}) + \epsilon (\tau_{8,2} +
\tau_{8,4})  + \eta\tau_{8,3}.
\end{array}
\right .
$$
The degree of $\s_{4,1}$ resp. $\s_{4,2},\tau_{4,1},\tau_{4,2},\tau_{4,3}$ is
$2$ resp. $4,1,2,1$ so we have
\begin{equation}
\label{equa-degre-tau-4}
a+b=1 \mbox{ and } c+d=2.
\end{equation}
The degree of
$\s_{8,1}$ resp. $\s_{8,2},\tau_{8,1},\tau_{8,2},\tau_{8,3},
\tau_{8,4},\tau_{8,5}$ is
$96$ resp. $72,12,21,30,21,12$ so we have
\begin{equation}
\label{equa-degre-tau-8}
24\a + 42\b + 30 \gamma = 96 \mbox{ and }
24\delta + 42\epsilon + 30 \eta = 72.
\end{equation}

To get more precise information we use the relation $\s^{4,2} \cup
\s^{4,2} = \s^{8,1} + \s^{8,2}$, which follows from the fact that the
degree of $(\s^{4,1})^2$ resp. $\s^{8,1} , \s^{8,2}$ is $56$
resp. $40,16$ (here we identify via Poincar{\'e} duality the cohomology 
classes $\s^{8,i}$ with the homology classes $\s_{7,i}$ for 
$i\in\{1,2\}$). We deduce the relations $\s^{4,1} \cup \s^{4,2} = 3 \s^{8,1}
+ 2 \s^{8,2}$ and $(\s^{4,1})^2 = 8 \s^{8,1} + 6 \s^{8,2}$. Thus one computes
that $\iota^* \tau^{4,1} \cup \iota^* \tau^{4,1} = (8a^2+6ac+c^2) \s^{8,1}
+ (6a^2+4ac+c^2) \s^{8,2}$.

On the other hand using the jeu de taquin rule
we have $\tau^{4,1} \cup \tau^{4,1} = \tau^{8,2}$ so
$\iota^*(\tau^{4,1} \cup \tau^{4,1}) = \iota^* \tau^{8,2} = \b \s_{8,1} + 
\epsilon \s_{8,2}$.
This implies that $\beta = 8a^2 + 6ac + c^2$
and $\epsilon=6a^2+4ac+c^2$. By (\ref{equa-degre-tau-8})
we have $\beta \leq 2$ so $a=0$ and $\beta=c=1$.
By (\ref{equa-degre-tau-4}) and (\ref{equa-degre-tau-8}) we deduce the
result for $\iota_*$ applied to degree 4 and 8 classes.

To compute $\iota_*$ for classes of degree lower than 8, we use the
projection formula $h\cap\iota_*\s=\iota_*(h\cap\s)$. For example
applying this to $\s_{8,1}$ and $\s_{8,2}$ we get 
$$h\cap(\tau_{8,1}+\tau_{8,2}+\tau_{8,3}+\tau_{8,4}+\tau_{8,5})=
\iota_*(2\s_{7,1}+\s_{7,2})
\ {\rm and}\ h\cap(\tau_{8,2}+\tau_{8,3}+\tau_{8,4})= 
\iota_*(\s_{7,1}+2\s_{7,2}).$$
Resolving this system gives the result in degree 7. The same procedure
gives the result in lower degrees.
\end{proo} 

\begin{rema}
Let us also remark that there is only one class in $H_*(F_4/P_1)$ in
degree 3. We denote this class $\s_3$. We have
$\iota_*\s_3=a\tau_{3,1}+b\tau_{3,2}$ but
$2=\deg(\s_3)=h^3\cap\iota^*\s_3= ah^3\cap\tau_{3,1}+b
h\cap\tau_{3,2}=a+b$ thus $a=b=1$ by symmetry and
$\iota_*\s_3=\tau_{3,1}+\tau_{3,2}$.
\end{rema}

We need to extend the Dynkin diagrams of $F_4$ and $E_6$. We first
consider the Kac-Moody groups $\widetilde{F}_4^2$ and
$\widetilde{E}_7^1$ with the notation of \cite{kac}. Their Dynkin
diagrams are:

\centerline{\begin{pspicture*}(0,0)(9.6,4.2000003)

\psellipse[fillstyle=solid,fillcolor=black](0.6,2.4)(0.096,0.096)
\psellipse[fillstyle=solid,fillcolor=black](1.2,2.4)(0.096,0.096)
\psellipse[fillstyle=solid,fillcolor=black](1.8,2.4)(0.096,0.096)
\psellipse[fillstyle=solid,fillcolor=black](2.4,2.4)(0.096,0.096)
\psellipse[fillstyle=solid,fillcolor=black](3.0,2.4)(0.096,0.096)
\psellipse[fillstyle=solid,fillcolor=black](5.4,2.4)(0.096,0.096)
\psellipse[fillstyle=solid,fillcolor=black](6.0,2.4)(0.096,0.096)
\psellipse[fillstyle=solid,fillcolor=black](6.6,2.4)(0.096,0.096)
\psellipse[fillstyle=solid,fillcolor=black](7.2,3.0)(0.096,0.096)
\psellipse[fillstyle=solid,fillcolor=black](7.2,2.4)(0.096,0.096)
\psellipse[fillstyle=solid,fillcolor=black](7.8,2.4)(0.096,0.096)
\psellipse[fillstyle=solid,fillcolor=black](8.4,2.4)(0.096,0.096)
\psellipse[fillstyle=solid,fillcolor=black](9.0,2.4)(0.096,0.096)
\psline(8.904,2.4)(8.496,2.4)
\psline(8.304,2.4)(7.896,2.4)
\psline(7.704,2.4)(7.296,2.4)
\psline(7.104,2.4)(6.696,2.4)
\psline(6.504,2.4)(6.096,2.4)
\psline(5.904,2.4)(5.496,2.4)
\psline(7.2,2.496)(7.2,2.904)
\psline(2.904,2.4)(2.496,2.4)
\psline(2.304,2.4)(1.896,2.4)
\psline(1.104,2.4)(0.696,2.4)
\psline(1.289,2.364)(1.711,2.364)
\psline(1.289,2.436)(1.711,2.436)
\psline(1.55,2.4)(1.458,2.477)
\psline(1.55,2.4)(1.458,2.323)
\uput[d](0.6,2.4){$\a_1$}
\uput[d](6.0,2.4){$\a_1$}
\uput[d](1.2,2.4){$\a_2$}
\uput[u](7.2,3.0){$\a_2$}
\uput[d](1.8,2.4){$\a_3$}
\uput[d](6.6,2.4){$\a_3$}
\uput[d](2.4,2.4){$\a_4$}
\uput[d](7.2,2.4){$\a_4$}
\uput[d](3.0,2.4){$\a_5$}
\uput[d](7.8,2.4){$\a_5$}
\uput[d](8.4,2.4){$\a_6$}
\uput[d](9.0,2.4){$\a_7$}
\uput[d](5.4,2.4){$\a_0$}
\put(1.8,1.2){$\widetilde{F}_4^2$}
\put(7.2,1.2){$\widetilde{E}_7^1$}

\end{pspicture*}}

Any length 8 element is $\Lambda_F$-cominuscule and there are three
new $\Lambda_F$-cominuscule heaps of length 8 in $\widetilde{F}_4^2$
with heaps as follows:

\centerline{\begin{pspicture*}(0,0)(12.6,5.4)

\psellipse[fillstyle=solid,fillcolor=black](0.6,4.2)(0.096,0.096)
\psellipse[fillstyle=solid,fillcolor=black](0.6,1.8)(0.096,0.096)
\psellipse[fillstyle=solid,fillcolor=black](1.2,3.6)(0.096,0.096)
\psellipse[fillstyle=solid,fillcolor=black](1.2,2.4)(0.096,0.096)
\psellipse[fillstyle=solid,fillcolor=black](1.8,4.2)(0.096,0.096)
\psellipse[fillstyle=solid,fillcolor=black](1.8,3.0)(0.096,0.096)
\psellipse[fillstyle=solid,fillcolor=black](2.4,3.6)(0.096,0.096)
\psellipse[fillstyle=solid,fillcolor=black](3.0,4.2)(0.096,0.096)
\psellipse[fillstyle=solid,fillcolor=black](4.8,1.8)(0.096,0.096)
\psellipse[fillstyle=solid,fillcolor=black](5.4,4.8)(0.096,0.096)
\psellipse[fillstyle=solid,fillcolor=black](5.4,3.6)(0.096,0.096)
\psellipse[fillstyle=solid,fillcolor=black](5.4,2.4)(0.096,0.096)
\psellipse[fillstyle=solid,fillcolor=black](6.0,4.2)(0.096,0.096)
\psellipse[fillstyle=solid,fillcolor=black](6.0,3.0)(0.096,0.096)
\psellipse[fillstyle=solid,fillcolor=black](6.6,3.6)(0.096,0.096)
\psellipse[fillstyle=solid,fillcolor=black](7.2,4.2)(0.096,0.096)
\psellipse[fillstyle=solid,fillcolor=black](9.0,1.8)(0.096,0.096)
\psellipse[fillstyle=solid,fillcolor=black](9.6,3.6)(0.096,0.096)
\psellipse[fillstyle=solid,fillcolor=black](9.6,2.4)(0.096,0.096)
\psellipse[fillstyle=solid,fillcolor=black](10.2,4.2)(0.096,0.096)
\psellipse[fillstyle=solid,fillcolor=black](10.2,3.0)(0.096,0.096)
\psellipse[fillstyle=solid,fillcolor=black](10.8,4.8)(0.096,0.096)
\psellipse[fillstyle=solid,fillcolor=black](10.8,3.6)(0.096,0.096)
\psellipse[fillstyle=solid,fillcolor=black](11.4,4.2)(0.096,0.096)
\psline(0.668,1.868)(1.132,2.332)
\psline(1.868,3.068)(2.332,3.532)
\psline(2.468,3.668)(2.932,4.132)
\psline(2.332,3.668)(1.868,4.132)
\psline(1.132,3.668)(0.668,4.132)
\psline(1.288,3.637)(1.763,4.112)
\psline(1.237,3.688)(1.712,4.163)
\psline(1.55,3.95)(1.43,3.94)
\psline(1.55,3.95)(1.54,3.83)
\psline(1.237,3.512)(1.712,3.037)
\psline(1.288,3.563)(1.763,3.088)
\psline(1.55,3.25)(1.54,3.37)
\psline(1.55,3.25)(1.43,3.26)
\psline(1.288,2.437)(1.763,2.912)
\psline(1.237,2.488)(1.712,2.963)
\psline(1.55,2.75)(1.43,2.74)
\psline(1.55,2.75)(1.54,2.63)
\psline(4.868,1.868)(5.332,2.332)
\psline(6.068,3.068)(6.532,3.532)
\psline(6.668,3.668)(7.132,4.132)
\psline(6.532,3.668)(6.068,4.132)
\psline(5.488,2.437)(5.963,2.912)
\psline(5.437,2.488)(5.912,2.963)
\psline(5.75,2.75)(5.63,2.74)
\psline(5.75,2.75)(5.74,2.63)
\psline(5.437,3.512)(5.912,3.037)
\psline(5.488,3.563)(5.963,3.088)
\psline(5.75,3.25)(5.74,3.37)
\psline(5.75,3.25)(5.63,3.26)
\psline(5.488,3.637)(5.963,4.112)
\psline(5.437,3.688)(5.912,4.163)
\psline(5.75,3.95)(5.63,3.94)
\psline(5.75,3.95)(5.74,3.83)
\psline(5.437,4.712)(5.912,4.237)
\psline(5.488,4.763)(5.963,4.288)
\psline(5.75,4.45)(5.74,4.57)
\psline(5.75,4.45)(5.63,4.46)
\psline(9.068,1.868)(9.532,2.332)
\psline(10.268,3.068)(10.732,3.532)
\psline(10.868,3.668)(11.332,4.132)
\psline(11.332,4.268)(10.868,4.732)
\psline(10.732,4.732)(10.268,4.268)
\psline(10.268,4.132)(10.732,3.668)
\psline(9.688,2.437)(10.163,2.912)
\psline(9.637,2.488)(10.112,2.963)
\psline(9.95,2.75)(9.83,2.74)
\psline(9.95,2.75)(9.94,2.63)
\psline(9.637,3.512)(10.112,3.037)
\psline(9.688,3.563)(10.163,3.088)
\psline(9.95,3.25)(9.94,3.37)
\psline(9.95,3.25)(9.83,3.26)
\psline(9.688,3.637)(10.163,4.112)
\psline(9.637,3.688)(10.112,4.163)
\psline(9.95,3.95)(9.83,3.94)
\psline(9.95,3.95)(9.94,3.83)
\uput[d](1.8,1.2){$\s_{8,3}$}
\uput[d](6.0,1.2){$\s_{8,4}$}
\uput[d](10.2,1.2){$\s_{8,5}$}

\end{pspicture*}}

We shall also consider the following heap in $\widetilde{E}_7^1/P_2$:

\centerline{\begin{pspicture*}(0,0)(4.8,3.6000001)

\psellipse[fillstyle=solid,fillcolor=black](0.6,3.0)(0.096,0.096)
\psellipse[fillstyle=solid,fillcolor=black](1.2,2.4)(0.096,0.096)
\psellipse[fillstyle=solid,fillcolor=black](1.8,3.0)(0.096,0.096)
\psellipse[fillstyle=solid,fillcolor=black](1.8,1.8)(0.096,0.096)
\psellipse[fillstyle=solid,fillcolor=black](2.4,3.0)(0.096,0.096)
\psellipse[fillstyle=solid,fillcolor=black](2.4,2.4)(0.096,0.096)
\psellipse[fillstyle=solid,fillcolor=black](2.4,1.2)(0.096,0.096)
\psellipse[fillstyle=solid,fillcolor=black](2.4,0.6)(0.096,0.096)
\psellipse[fillstyle=solid,fillcolor=black](3.0,3.0)(0.096,0.096)
\psellipse[fillstyle=solid,fillcolor=black](3.0,1.8)(0.096,0.096)
\psellipse[fillstyle=solid,fillcolor=black](3.6,2.4)(0.096,0.096)
\psellipse[fillstyle=solid,fillcolor=black](4.2,3.0)(0.096,0.096)
\psline(0.668,2.932)(1.132,2.468)
\psline(1.268,2.332)(1.732,1.868)
\psline(1.868,1.732)(2.332,1.268)
\psline(2.4,1.104)(2.4,0.696)
\psline(4.132,2.932)(3.668,2.468)
\psline(3.532,2.332)(3.068,1.868)
\psline(2.932,1.732)(2.468,1.268)
\psline(2.4,2.904)(2.4,2.496)
\psline(2.468,2.332)(2.932,1.868)
\psline(2.332,2.332)(1.868,1.868)
\psline(1.732,2.932)(1.268,2.468)
\psline(1.868,2.932)(2.332,2.468)
\psline(2.932,2.932)(2.468,2.468)
\psline(3.068,2.932)(3.532,2.468)

\end{pspicture*}}

We complete our notation and define homology classes in
$H_*(\widetilde{F}_4^2/P_1)$. 
The previous classes are again classes
and there are few more classes to obtain all classes of degree
$d\in[4,8]$.
$$\begin{array}{llll}
  \s_{5,3}=\scal{(\a_5,1)}&
  \s_{6,3}=\scal{(\a_2,2),(\a_5,1)} &
  \s_{7,3}=\scal{(\a_1,2),(\a_5,1)} &
  \s_{8,3}=\scal{(\a_1,2),(\a_3,2),(\a_5,1)}
  \\
 &
 &
\s_{7,4}=\scal{(\a_3,2),(\a_5,1)} &
\s_{8,4}=\scal{(\a_2,3),(\a_5,1)} \\
&
&
&
\s_{8,5}=\scal{(\a_4,2)}\\
\end{array}$$

In the same way, we complete our notation and define homology classes in
$H_*(\widetilde{E}_7^1/P_2)$. 
The previous classes are again classes
and there are few more classes to obtain all classes of degree
$d\in[4,8]$.
We define
$$\begin{array}{llll}
  \tau_{5,4}=\scal{(\beta_0,1)}&
  \tau_{6,5}=\scal{(\beta_0,1),(\beta_5,1)} &
  \tau_{7,6}=\scal{(\beta_0,1),(\beta_4,2)}&
  \tau_{8,6}=\scal{(\beta_0,1),(\beta_3,2)} \\
  \tau_{5,5}=\scal{(\beta_7,1)}&
  \tau_{6,6}=\scal{(\beta_3,1),(\beta_7,1)} &
  \tau_{7,7}=\scal{(\beta_0,1),(\beta_6,1)}&
  \tau_{8,7}=\scal{(\beta_0,1),(\beta_2,2)} \\
&&
  \tau_{7,8}=\scal{(\beta_1,1),(\beta_7,1)}&
  \tau_{8,8}=\scal{(\beta_0,1),(\beta_4,2),(\beta_6,1)}   \\
&&
  \tau_{7,7}=\scal{(\beta_4,2),(\beta_7,1)}&
  \tau_{8,9}=\scal{(\beta_0,1),(\beta_7,1)} \\
&&&  \tau_{8,{10}}=\scal{(\beta_1,1),(\beta_4,2),(\beta_7,1)} \\
&&&  \tau_{8,{11}}=\scal{(\beta_2,2),(\beta_7,1)} \\
&&&  \tau_{8,{12}}=\scal{(\beta_5,2),(\beta_7,1)} \\
\end{array}$$

We prove the following 

\begin{prop}
\label{prop-f_4}
  We have the formula
$$\tau_{4,1}\cap \iota_*\s_{8,3}=4\tau_{4,1}+12\tau_{4,2}+
4\tau_{4,3}.$$
\end{prop}

Before going into the proof of this proposition, which is a long but
simple computation we prove

\begin{coro}
\label{coro-f_4}
  We have the equalities
$$c_{\s_{4,2},\s_{4,2}}^{\s_{8,3}}=4=m_{\s_{4,2},\s_{4,2}}^{\s_{8,3}}\cdot
t_{\s_{4,2},\s_{4,2}}^{\s_{8,3}} 
 \ \ {\rm and}\ \
 c_{\s_{4,1},\s_{4,2}}^{\s_{8,3}}=8=m_{\s_{4,1},\s_{4,2}}^{\s_{8,3}}\cdot
 t_{\s_{4,1},\s_{4,2}}^{\s_{8,3}}.$$
\end{coro}

\begin{proo}
  By Lemma \ref{lemm-f_4}, we have in $F_4/P_1$ the equality 
$\scal{\iota^*\tau^{4,1},\s_{4,i}}=\scal{\tau^{4,1},\iota_*\s_{4,i}}
=\delta_{i,2}$.
In particular, this implies the equality $\iota^*\tau^{4,1}=\s^{4,2}$. 
On the other hand, Lemma \ref{lemm-f_4} and the previous Proposition
imply the equality
$\tau^{4,1}\cap\iota_*\s_{8,3}=\iota_*(8\s_{4,1}+4\s_{4,2})$. 
We compute 
$$ \begin{array}{ll}
\iota_*(\s^{4,2}\cap\s_{8,3})&=\iota_*(\iota^*\tau^{4,1}\cap\s_{8,3})\\
&=\tau^{4,1}\cap\iota_*\s_{8,3}\\
&=4\tau_{4,1}+12\tau_{4,2}+4\tau_{4,3}\\
&=\iota_*(8\s_{4,1}+4\s_{4,2}).\\
\end{array}$$
The result follows by injectivity of $\iota_*$.
\end{proo}

\begin{proo-prop}
  The main tool here will be the fact that the pull-back by $\iota$ of
  an hyperplane section is again an hyperplane section. We will write
  this as $\iota^*h=h$ and use it with projection formula to obtain
  
  \begin{equation}
    \label{projection}
h\cap\iota_*\s=\iota_*(h\cap\s)
  \end{equation}
where $\s\in H_*(\widetilde{F}_4^2/P_1)$. We shall also use the
following observation: for $\s\in H_*(\widetilde{F}_4^2/P_1)$ and
$\tau\in H^*(\widetilde{E}_7^1/P_2)$, the cap product
$\tau\cap\iota_*\s$ is symmetric with respect to the folding. Indeed,
we have $\tau\cap\iota_*\s=\iota_*(\iota^*\tau\cap\s)$. 
We shall in particular need the following cap products (we compute them using
the product $\pp$ which is valid for all degree 8 classes $\s_\lt$ in
$H^*(\widetilde{E}_7^1/P_2)$ because $D_0(\lt)$ is of finite type and because
  we have already proved the simply laced case).

\vskip 0.2 cm

$$\begin{array}{l|cccccc}
& \tau_{6,1} & \tau_{6,2} & \tau_{6,3} & \tau_{6,4} & \tau_{6,5} & 
\tau_{6,6}  \\
\hline
  \tau^{3,1}\cap\bullet &2\tau_{3,1}+\tau_{3,2}&\tau_{3,2}
  &\tau_{3,1}+2\tau_{3,2} &\tau_{3,1}+\tau_{3,2}
  &2\tau_{3,1}+\tau_{3,2}&\tau_{3,2} \\
\hline
\hline
& \tau_{8,1} & \tau_{8,2} & \tau_{8,3} & \tau_{8,4} & \tau_{8,5} & 
\tau_{8,6} \\
\hline
  \tau^{4,1}\cap\bullet &\tau_{4,2}&\tau_{4,1}+\tau_{4,2}
  &\tau_{4,2}+\tau_{4,3} &\tau_{4,2}
  &0&2\tau_{4,1}+\tau_{4,2} \\
\hline
\hline
& \tau_{8,7} & \tau_{8,8} & \tau_{8,9} & \tau_{8,10} & \tau_{8,11} & 
\tau_{8,12}\\  
\hline
  \tau^{4,1}\cap\bullet &2\tau_{4,2}&\tau_{4,1}+3\tau_{4,2}+\tau_{4,3}
  &\tau_{4,2}+2\tau_{4,3} &\tau_{4,2}+\tau_{4,3}
  &0&0 \\
\end{array}$$

\vskip 0.2 cm

We will not
explicitly commute the direct image $\iota_*\s_{8,3}$ (we will have four
possible solutions) but this will be enough to get the result.

Write
$\iota_*\s_{5,3}=a(\tau_{5,1}+\tau_{5,3})+b\tau_{5,2}+c(\tau_{5,4}+
\tau_{5,5})$ with $(a,b,c)$ non negative integers.
By equation (\ref{projection}), we get
$$2(\tau_{4,1}+\tau_{4,2}+\tau_{4,3})=\iota_*(2\s_{4,2})
=\iota_*(h\cap\s_{5,3}) = h\cap\iota_*(a(\tau_{5,1}+\tau_{5,3})+
b\tau_{5,2}+c(\tau_{5,4}+\tau_{5,5}))
$$
and the equalities $2a+b=2=a+c$. The only solutions are
$(a,b,c)=(1,0,1)$ or $(0,2,2)$. 

Now write
$\iota_*\s_{6,3}=\a(\tau_{6,1}+\tau_{6,4})+\beta\tau_{6,2}
+\gamma\tau_{6,3}+ \delta(\tau_{6,5}+\tau_{6,6})$ with 
$(\a,\beta,\gamma,\delta)$ non
negative integers. As before, we get the equalities 
$\delta=c$, $\a+\gamma+\delta=a+2$, $2\a+\beta=b+2$. 
If $(a,b,c)=(1,0,1)$ then $(\a,\beta,\gamma,\delta)=(0,2,2,1)$ or
$(1,0,1,1)$ and if $(a,b,c)=(0,2,2)$ then 
$(\a,\beta,\gamma,\delta)=(0,4,0,2)$. Computing the cap product
$\tau^{3,1}\cap\iota_*\s_{6,3}$ we see that the only solution for
$(\a,\beta,\gamma,\delta)$ such that $\tau^{3,1}\cap\iota_*\s_{6,3}$ 
is symmetric with respect to the folding is $(1,0,1,1)$ and we deduce
that $(a,b,c)=(1,0,1)$.

Let us now write
$\iota_*\s_{7,3}=x(\tau_{7,1}+\tau_{7,5})+y(\tau_{7,2}+\tau_{7,4})
+z\tau_{7,3}+ t(\tau_{7,6}+\tau_{7,9})+ u(\tau_{7,7}+\tau_{7,8})$ with 
$(x,y,z,t,u)$ non
negative integers. As before, we get the equalities 
$x+y+z+t=3$, $2y=2$, $z+2u=1$, $t+u=1$. The only solution is
$(x,y,z,t,u)=(0,1,1,1,0)$.

Write
$\iota_*\s_{7,4}=x'(\tau_{7,1}+\tau_{7,5})+y'(\tau_{7,2}+\tau_{7,4})
+z'\tau_{7,3}+ t'(\tau_{7,6}+\tau_{7,9})+ u'(\tau_{7,7}+\tau_{7,8})$ 
with $(x',y',z',t',u')$
non negative integers. As before, we get the equalities 
$x'+y'+z'+t'=4$, $2y'=0$, $z'+2u'=4$, $t'+u'=2$. The only solutions are
$(x',y',z',t',u')=(1,0,2,1,1)$ and $(4,0,0,0,2)$.

Write
$\iota_*\s_{8,3}=A(\tau_{8,1}+\tau_{8,5})+B(\tau_{8,2}+\tau_{8,4})
+C\tau_{8,3}+ D(\tau_{8,6}+\tau_{8,12})+E(\tau_{8,7}+\tau_{8,11})
+ F(\tau_{8,8}+\tau_{8,10})+G\tau_{8,9}$ with $(A,B,C,D,E,F,G)$
non negative integers. We get the equalities 
$A+B+D=x'+2$, $A+C+E=y'+4$, $2B+C+2F=z'+4$, $D+E+F=t'+2$, $F+G=u'$. If
$(x',y',z',t',u')=(1,0,2,1,1)$ then 
$(A,B,C,D,E,F,G)=(0,2,2,1,2,0,1)$ or $(1,1,2,1,1,1,0)$ and if
$(x',y',z',t',u')=(4,0,0,0,2)$ then $(A,B,C,D,E,F)=(4,1,0,1,0,1,1)$ or
$(3,2,0,1,1,0,2)$. We now compute for all these solution the cap product with
$\tau^{4,1}$. It gives in all cases $\tau^{4,1}\cap
\iota_*\s_{8,3}=4\tau_{4,1}+12\tau_{4,2}+4\tau_{4,3}.$ 
\end{proo-prop}

\subsubsection{Proof for $F_4$}

We consider  the system of $\varpi_1$-cominuscule $F_4$-colored posets
$\pos_0$ given by the unique following poset:

\centerline{\begin{pspicture*}(0,0)(3.6000001,3.6000001)

\psellipse(0.6,3.0)(0.096,0.096)
\psline(0.532,2.932)(0.668,3.068)
\psline(0.532,3.068)(0.668,2.932)
\psellipse[fillstyle=solid,fillcolor=black](0.6,0.6)(0.096,0.096)
\psellipse(1.2,2.4)(0.096,0.096)
\psline(1.132,2.332)(1.268,2.468)
\psline(1.132,2.468)(1.268,2.332)
\psellipse[fillstyle=solid,fillcolor=black](1.2,1.2)(0.096,0.096)
\psellipse(1.8,3.0)(0.096,0.096)
\psline(1.732,2.932)(1.868,3.068)
\psline(1.732,3.068)(1.868,2.932)
\psellipse[fillstyle=solid,fillcolor=black](1.8,1.8)(0.096,0.096)
\pspolygon[fillstyle=solid,fillcolor=black](2.52,2.4)(2.4,2.52)(2.28,2.4)(2.4,2.28)
\psline(1.704,3.0)(1.704,3.0)
\psline(1.868,2.932)(2.34,2.46)
\psline(2.34,2.34)(1.868,1.868)
\psline(1.132,1.132)(0.668,0.668)
\psline(0.668,2.932)(1.132,2.468)
\psline(1.288,2.437)(1.763,2.912)
\psline(1.237,2.488)(1.712,2.963)
\psline(1.55,2.75)(1.43,2.74)
\psline(1.55,2.75)(1.54,2.63)
\psline(1.237,2.312)(1.712,1.837)
\psline(1.288,2.363)(1.763,1.888)
\psline(1.55,2.05)(1.54,2.17)
\psline(1.55,2.05)(1.43,2.06)
\psline(1.288,1.237)(1.763,1.712)
\psline(1.237,1.288)(1.712,1.763)
\psline(1.55,1.55)(1.43,1.54)
\psline(1.55,1.55)(1.54,1.43)
\uput[r](2.4,2.4){$\s^{4,2}$}

\end{pspicture*}}

We have $S_0=\{4\}$. Let $(D_4,d_4)$ be a
marked Dynkin diagram and $P_4$ be any $d_4$-minuscule $D_4$-colored
poset. Set $\pos=\pos_{\pos_0,\{P_4\}}$. 

\begin{lemm}
\label{lemm-f4p1}
With the above notation, assume that Conjecture \ref{main_conj} holds
for $P_4$ and any $\lt$ in $I(\pos)$ with
$D_0(\lt)\varsubsetneq F_4$. Then Conjecture \ref{main_conj} holds for
$\pos$. 
\end{lemm}

\begin{proo}
  Choose some generators $\g^1$ and $\g^4$ of degree 1 and 4 of
  $H^*(\pos_0)$. It is easy to see that we may choose $\g^4=\s^{4,2}$ 
with the notation of the
  previous section. The variety $F_4/P_1$ has dimension 15 and the
dimensions of $H^d(F_4/P_1)$ are

\vskip 0.2 cm

\centerline{\begin{tabular}{c|lllllllllllllllllllllll}
\hline
$d$&0&1&2&3&4&5&6&7&8\\
\hline
$\dim H^d(F_4/P_1)$&1&1&1&1&2&2&2&2&2\\
\hline
\end{tabular}}

\vskip 0.2 cm

In particular by Lemma \ref{lefs}, we only need to prove
$\ok{\g^4}{\g^4}$.
Since by assumption the conjecture holds for any
$\lambda \in I(\pos)$ with $D_0(\lambda) \varsubsetneq F_4$,
we have $\oknu{\g^4}{\g^4}{\s^\nu}$ as soon as the ideal $\nu$, of degree
8 satisfies $(\a_1,2)\not\in\nu$ or $(\a_3,2)\not\in\nu$. There is a
unique such ideal $\nu$ in $\pos$. We denote it by $\nu'$.

We first deal with the case
$D_0(\nu')=\widetilde{F}_4^2$. In that case, the class $\s^{\nu'}$ is
$\s^{8,3}$ in the notation of the previous section. In particular, we
have $c_{\g^4,\g^4}^{\s^{8,3}}=c_{\s^{4,2},\s^{4,2}}^{\s^{8,3}}=4$ by
Corollary \ref{coro-f_4}. 

Now we deal with the general case where $D_0(\nu')$ is obtained from
$F_4$ by adding one vertex with $n$-tuple edge linking it to the
simple root $\a_4$. By Proposition \ref{prop-iota_*} and with the
notation of that proposition, we have
$\iota_*\s_{\nu'}=\sum_{i=1}^n\tau_{\nu_i}$. 
We then have, because $\iota^*\g^4=\g^4$, the equality 
$$\g^4\cap\s_{\nu'}=\sum_{i=1}^n\g^4\cap\tau_{{\nu_i}}$$
and it follows that
$c_{\g^4,\g^4}^{\s^{\nu'}}=\sum_{i=1}^nc_{\g^4,\g^4}^{\tau^{\nu_i}}$
and the result follows.
\end{proo}

%%% Local Variables: 
%%% mode: latex
%%% TeX-master: t
%%% End: 

\bigskip\noindent
Pierre-Emmanuel {\sc Chaput}, \\
{\it Laboratoire de Math{\'e}matiques Jean Leray,} 
UMR 6629 du CNRS, UFR Sciences et Techniques,  2 rue de la
Houssini{\`e}re, BP 92208, 44322 Nantes cedex 03, France and \\
{\it Max Planck Institut f{\"u}r Mathematik, }
Vivatsgasse 7, 
53111 Bonn, Germany.

\noindent {\it email}: \texttt{pierre-emmanuel.chaput@math.univ-nantes.fr}.

\medskip\noindent
Nicolas {\sc Perrin}, \\
{\it Hausdorff Center for Mathematics,}
Universit{\"a}t Bonn, Villa Maria, Endenicher
Allee 62, 
53115 Bonn, Germany and \\
{\it Institut de Math{\'e}matiques de Jussieu,} 
Universit{\'e} Pierre et Marie Curie, Case 247, 4 place
Jussieu, 75252 Paris Cedex 05, France.

\noindent {\it email}: \texttt{nicolas.perrin@hcm.uni-bonn.de}.


\begin{thebibliography}{BeGeGe73}

\bibitem[BeGeGe73]{BGG} Bern{\u s}te{\u \i}n, I. N., Gel'fand, I. M., 
Gel'fand, S. I., {\it Schubert cells and the cohomology of a flag space.}
Funkcional. Anal. i Prilo{\u z}en. {\bf 7} (1973), no. 1.

\bibitem[Bor53]{borel} Borel, A., {\it Sur la cohomologie des espaces 
fibr{\'e}s principaux et des espaces homog{\`e}nes de groupes de Lie 
compacts.} Ann. of Math. (2) {\bf 57}, (1953).

\bibitem[Bou54]{bou} Bourbaki, N., {\it Groupes et alg{\`e}bres de Lie}.
  Hermann 1954.

\bibitem[ChMaPe08]{CMP} Chaput, P.-E., Manivel, L., Perrin, N., {\it Quantum
    cohomology of minuscule homogeneous spaces.}  Transform. Groups  {\bf 13} 
    (2008), no. 1, 47--89.

\bibitem[ChPe09]{CP} Chaput, P.-E., Perrin, N., {\it Quantum cohomology
    of adjoint homogeneous spaces}. In preparation.

\bibitem[Dem74]{demazure} Demazure, M., {\it D{\'e}singularisation 
des vari{\'e}t{\'e}s de Schubert g{\'e}n{\'e}ralis{\'e}es.} Collection 
of articles dedicated to Henri Cartan on the occasion of his 70th birthday, I.
Ann. Sci. {\'E}cole Norm. Sup. (4) {\bf 7} (1974).

\bibitem[Dua05]{duan} Duan, H., {\it Multiplicative rule of Schubert class}.  
Invent. Math.  {\bf 159}  (2005),  no. 2.

\bibitem[Kac90]{kac} Kac, V. G., {\it Infinite-dimensional Lie algebras.}
 Third edition. Cambridge University Press, Cambridge, 1990.

\bibitem[Kum02]{kumar} Kumar, S., {\it Kac-Moody groups, their flag varieties 
and representation theory.} Progress in Mathematics, 204. Birkh{\"a}user 
Boston, 
Inc., Boston, MA, 2002.

\bibitem[Laz04]{lazarsfeld} Lazarsfeld, R., {\it Positivity in algebraic
    geometry. I. Classical setting: line bundles and linear series.}
  Ergebnisse der Mathematik und ihrer Grenzgebiete. 3. Folge. A Series
  of Modern Surveys in Mathematics, 48. Springer-Verlag, Berlin, 2004.

\bibitem[LiRi34]{LR} Littlewood, D.E., Richardson, A.R., {\it Group 
characters and  algebras}. Phil. Trans. R. Soc. London Ser. A {\bf 233} 
(1934).

\bibitem[Per07]{small} Perrin, N., {\it Small resolutions of minuscule 
Schubert varieties.} Compos. Math.  {\bf 143}  (2007),  no. 5.

\bibitem[Pra91]{pragacz} Pragacz, P,. {\it Algebro-geometric applications 
of Schur $S$- and $Q$-polynomials}. Topics in invariant theory 
(Paris, 1989/1990), Lecture Notes in Math., 1478, Springer, Berlin, 1991.

\bibitem[Pro99a]{proctor_minuscule} Proctor, R. A., {\it Minuscule elements 
of Weyl groups, the numbers game, and $d$-complete posets.}
J. Algebra  {\bf 213}  (1999),  no. 1.

\bibitem[Pro99b]{proctor-classification} Proctor, R. A., {\it Dynkin diagram 
classification of $\lambda$-minuscule Bruhat lattices and of 
$d$-complete posets.}  J. Algebraic Combin.  {\bf 9}  (1999),  no. 1.

\bibitem[Pro04]{proctor_taquin} Proctor, R. A., {\it d-Complete posets 
generalize Young diagrams for the jeu de taquin property}. preprint 2004, 
available at \texttt{http://www.math.unc.edu/Faculty/rap}.

\bibitem[Sch77]{schutz} Sch{\"u}tzenberger, M.-P., {\it La corres\-pondance 
de Robinson}. Combinatoire et repr{\'e}\-sentation du groupe sym{\'e}trique 
(Actes Table Ronde CNRS, Univ. Louis-Pasteur Strasbourg, Strasbourg, 1976). 
Lecture Notes in Math., Vol. 579, Springer, Berlin, 1977. 

\bibitem[Ste96]{stembridge} Stembridge, J. R., {\it On the fully 
commutative
elements of Coxeter groups}.  J. Algebraic Combin.  {\bf 5}  (1996),  
no. 4.

\bibitem[Ste01]{stembridge-classification} Stembridge, J. R.,
{\it Minuscule elements of Weyl groups.}  J. Algebra  {\bf 235}  (2001),
no. 2.

\bibitem[ThYo08]{TY} Thomas, H., Yong, A., {\it A combinatorial rule for
(co)minuscule Schubert calculus}. To appear in Advances in Math.

\bibitem[VLe01]{vanlee} Van Leeuwen, M. A., {\it The Littlewood-Richardson 
rule, and related combinatorics.} Interaction of combinatorics and 
representation theory, MSJ Mem., {\bf 11}, Math. Soc. Japan, Tokyo, 2001.

\bibitem[Vie86]{viennot} Viennot, G. X., {\it Heaps of pieces. I. 
Basic definitions and combinatorial lemmas.} Combinatoire {\'e}num{\'e}rative 
(Montreal, Que., 1985/Quebec, Que., 1985), Lecture Notes 
in Math., {\bf 1234}, Springer, Berlin, 1986.

\bibitem[Wor84]{worley} Worley, D. R., {\it A theory of shifted Young 
tableaux}. Thesis, MIT, 1984. Available at 
\texttt{http://hdl.handle.net/1721.1/15599}.
\end{thebibliography}
\end{document}